\documentclass[final,3p]{elsarticle}
\usepackage[english]{babel}
\usepackage{amsmath,amssymb,mathtools}
\usepackage{amsthm}
\usepackage{algorithm2e}
\usepackage{enumitem}
\usepackage{bm}
\usepackage{empheq}
\usepackage{booktabs}
\usepackage{graphicx}
\usepackage{subcaption}
\usepackage{caption}
\usepackage{subcaption}

\usepackage{float}
\usepackage{placeins}
\usepackage{comment}

\usepackage{pdflscape}
\usepackage{xcolor}
\usepackage[colorlinks=true, allcolors=blue]{hyperref}
\usepackage[capitalize]{cleveref}
\allowdisplaybreaks

\numberwithin{equation}{section}
\newtheorem{theorem}{Theorem}[section]
\newtheorem{lemma}{Lemma}[section]

\newtheorem{remark}{Remark}[section]

\newtheorem{assumption}{Assumption}[section]
\newcommand{\lemref}[1]{\hyperref[#1]{Lemma~\ref*{#1}}}
\newcommand{\thmref}[1]{\hyperref[#1]{Theorem~\ref*{#1}}}
\newcommand{\defref}[1]{\hyperref[#1]{Definition~\ref*{#1}}}
\newcommand{\remref}[1]{\hyperref[#1]{Remark~\ref*{#1}}}
\newcommand{\probref}[1]{\hyperref[#1]{Problem~\ref*{#1}}}
\newcommand{\appref}[1]{\hyperref[#1]{Appendix~\ref*{#1}}}
\newcommand{\assref}[1]{\hyperref[#1]{Appendix~\ref*{#1}}}
\journal{Computer Methods in Applied Mechanics and Engineering}
\begin{document}
\begin{frontmatter}
\title{A Rigorous Fredholm Framework and Stabilized Equal-Order Finite Element Methods for the Time-Harmonic Biot Model} 

	\author[MS]{Mingchao Cai\corref{cor1}}
		\cortext[cor1]{Corresponding author.}\ead{Mingchao.Cai@morgan.edu}
        
       \author[SuTech]{Yiming Fu} 
	    \ead{12532024@mail.sustech.edu.cn}
        
		\author[SuTech]{Jingzhi Li} 
        \ead{lijz@sustech.edu.cn} 


        \affiliation[MS]{organization={Department of Mathematics},
			addressline={Morgan State University}, 
			city={Baltimore},
			state={MD 21251},
			country={USA}}
        
		\affiliation[SuTech]{organization={Department of Mathematics},
			addressline={Southern University of Science and Technology}, 
			city={Shenzhen},
            state={Guangdong},
			country={PR China}}
		
		
		
		
		

\begin{abstract}
We investigate the well-posedness and finite element approximation of a three-field formulation of the time-harmonic Biot model, motivated by frequency-domain magnetic resonance elastography. The underlying complex-valued mixed variational problem is noncoercive because of the inertial elasticity term and the coupled three-field structure. Continuous well-posedness is established under a non-resonance condition by decomposing the variational operator into an invertible principal part and a compact perturbation and applying the Fredholm alternative. We then develop stabilized equal-order finite element methods based on momentum-residual stabilization and local projection stabilization of the pore pressure. For the discrete problem, we establish mesh-uniform stability and well-posedness on sufficiently fine meshes via an a priori stability estimate and a Schatz-type compactness argument, without invoking discrete $T$-coercivity. Optimal-order a priori error estimates are derived in the natural energy norm. We also develop systematic selection principles for the stabilization parameters. Numerical experiments confirm the predicted convergence rates and demonstrate stable approximations for heterogeneous low-permeability media and anatomically realistic brain geometries relevant to magnetic resonance elastography.
\end{abstract}

\end{frontmatter}

\section{Introduction}
Biot's theory of poroelasticity provides a classical framework for describing the coupled deformation of fluid-saturated porous media, with applications in geomechanics, petroleum engineering, biomechanics, and biomedical imaging. The present work is motivated primarily by magnetic resonance elastography (MRE) \cite{hirsch2017magnetic,sack2023magnetic,galarce2023displacement}, a non-invasive imaging modality for estimating the mechanical properties of biological tissues. In MRE, tissue is subjected to harmonic mechanical excitation, typically in the low frequency range of $1$--$100\,\mathrm{Hz}$, and the induced displacement field is measured using phase-contrast magnetic resonance imaging. These measurements can then be used in inverse problems to estimate elastic moduli, hydraulic permeability, and other poroelastic parameters. Related partitioned methods for time-domain Biot systems have been analyzed in \cite{bukac2015analysis}. Such applications naturally motivate frequency-domain Biot models obtained through a single-frequency ansatz.
The resulting time-harmonic Biot model is a complex-valued coupled system. In contrast to the quasi-static case, the solid momentum equation contains an inertial term that leads to a Helmholtz-type indefinite elasticity operator \cite{cummings2006sharp}, so the associated variational formulation is generally noncoercive. Its analysis must therefore account simultaneously for operator indefiniteness, mixed saddle-point coupling, and complex-valued coefficients. These features also complicate the construction and analysis of stable finite element discretizations.

A useful framework is the three-field total-pressure formulation \cite{lee2017parameter, oyarzua2016locking}, in which displacement, pore pressure, and total pressure are treated as independent unknowns. The total-pressure variable separates the volumetric response from the deviatoric deformation and helps alleviate volumetric locking in the nearly incompressible regime. Cárcamo, Caiazzo, Galarce, and Mura~\cite{carcamo2024stabilized} proposed a continuous equal-order finite element discretization of this formulation, combining an elementwise momentum-residual stabilization with a Brezzi--Pitk\"aranta-type pressure-gradient stabilization \cite{brezzi1984stabilization} to address the lack of discrete inf--sup stability and pressure instability in low-permeability regimes.

The mathematical treatment of the time-harmonic Biot problem requires particular care. The inertial elasticity contribution cannot, in general, be handled by standard coercivity arguments, even under a non-resonance condition. A rigorous well-posedness analysis must instead identify an invertible principal operator, isolate the compact lower-order contributions, and apply the Fredholm alternative \cite{kress1977fredholm} to the full mixed system. At the discrete level, stability and convergence must be established consistently with the noncoercive continuous problem and the introduced stabilization mechanisms.

In this work, we develop a rigorous Fredholm framework for the three-field time-harmonic Biot problem, together with stability and convergence analyses for stabilized equal-order finite element approximations. Our continuous analysis begins with the operator decomposition
\[
\mathcal A=\mathcal B+\mathcal C,
\]
where $\mathcal B$ contains the coercive elasticity, hydraulic diffusion, total-pressure, and skew-Hermitian divergence-coupling contributions, while $\mathcal C$ collects the inertial and lower-order pressure--total-pressure couplings. We prove that $\mathcal B$ is an isomorphism and $\mathcal C$ is compact, thereby establishing continuous well-posedness under a non-resonance condition through the Fredholm alternative. This approach avoids imposing coercivity on the indefinite elasticity operator.

For the spatial discretization, we develop stabilized equal-order finite element methods for the displacement, pore pressure, and total pressure. An elementwise momentum-residual stabilization addresses the missing displacement--total-pressure inf--sup condition, while a pressure local projection stabilization \cite{ganesan2008local} controls the unresolved pore-pressure gradient modes. The latter reduces to a Brezzi--Pitk\"aranta-type stabilization \cite{brezzi1984stabilization} for linear elements and limits unnecessary artificial diffusion for higher-order elements. We further derive systematic selection principles for the stabilization parameters based on the scaling of the dynamic momentum operator and a reduced pressure Schur-complement criterion.

For the discrete analysis, we establish mesh-uniform stability and well-posedness on sufficiently fine meshes without relying on discrete $T$-coercivity \cite{ciarlet2012t}. Instead, we derive an a priori stability estimate, establish asymptotic consistency of the stabilization terms, and employ a Schatz-type compactness argument \cite{schatz1974observation}. By combining continuous injectivity with a contradiction argument, we prove a mesh-uniform discrete inf--sup condition, linking the continuous Fredholm theory to the stabilized finite element formulation.

Under standard regularity assumptions, we also establish optimal-order a priori error estimates for all three variables in the natural energy norm, consistent with the noncoercive structure of the time-harmonic problem. The dependence of the stability and approximation constants on the physical and stabilization parameters is also identified.

The theoretical results are supported by numerical experiments in two and three spatial dimensions. Benchmark problems with exact solutions confirm the predicted convergence rates, while layered heterogeneous media demonstrate stable pressure approximations in the presence of strongly discontinuous permeability coefficients. Finally, simulations on anatomically realistic brain geometries illustrate the applicability of the proposed method to magnetic resonance elastography.


\section{Model problem}

\subsection{A linear poroelastic model in the frequency domain}
We consider the following Biot model with a dynamic solid momentum balance coupled to quasi-static fluid mass conservation \cite{biot1941general,showalter2000diffusion,mella2023hybrid} posed on a bounded Lipschitz domain $\Omega\subset\mathbb{R}^d$ $(d=2,3)$:
\begin{empheq}[left=\empheqlbrace]{equation}
\begin{aligned}
\rho\bm{u}_{tt} - \mathrm{div}\,\boldsymbol{\sigma} + \alpha\nabla p &= \bm{f}
&& \text{in } \Omega \times (0,T],\\
S_\epsilon p_t + \alpha\,\nabla\!\cdot\!\bm{u}_t
-\frac{\kappa}{\mu_f}\Delta p &= g
&& \text{in } \Omega \times (0,T],
\end{aligned}
\label{eq:biot_time_domain}
\end{empheq}
where $\bm{u}$ and $p$ denote the solid displacement and pore-fluid pressure, respectively. The Cauchy stress tensor is given by
\begin{equation*}
\boldsymbol{\sigma}
=
2\mu_e\varepsilon(\bm{u})
+\lambda(\nabla\!\cdot\!\bm{u})\mathbb{I},
\qquad
\varepsilon(\bm{u})
=
\frac12\bigl(\nabla\bm{u}+(\nabla\bm{u})^\top\bigr),
\label{eq:cauchy_stress}
\end{equation*}
where
\[
\mu_e=\frac{E}{2(1+\nu)},
\qquad
\lambda=\frac{E\nu}{(1+\nu)(1-2\nu)}
\]
are the Lam\'e parameters determined by Young's modulus $E$ and Poisson's ratio $\nu$. The remaining physical parameters are the permeability $\kappa$, fluid viscosity $\mu_f$, bulk density $\rho$, Biot--Willis coefficient $\alpha$, body force $\bm{f}$, fluid source term $g$, and the constrained storage coefficient
\[
S_\epsilon
=
\frac{3\alpha(1-\alpha B)(1-2\nu)}{BE},
\]
where $B$ denotes Skempton's coefficient. We assume throughout that the physical parameters satisfy
\[
E>0,\qquad
0<\nu<\frac12,\qquad
\rho>0,\qquad
\mu_f>0,\qquad
\kappa>0,\qquad
\omega>0,\qquad
0<\alpha\leq 1,
\]
and $B>0$ with $\alpha B<1$, so that $S_\epsilon>0$. To avoid Poisson locking \cite{gu2023iterative,lee2017parameter,oyarzua2016locking}, system \eqref{eq:biot_time_domain} is reformulated in a three-field form by introducing the total pressure 
\begin{equation}
\phi=\alpha p-\lambda\operatorname{tr}(\varepsilon(\bm{u})).
\label{phi}
\end{equation}

We assume a time-harmonic solution of the form
\[
\bm u(x,t) = \bm u(x)e^{i\omega t},\qquad
p(x,t) = p(x)e^{i\omega t},
\]
where $\omega=2\pi f$ denotes the angular frequency and $f$ is the excitation frequency. Consequently, we consider the following frequency-domain Biot system corresponding to a prescribed angular frequency $\omega$:

\begin{empheq}[left=\empheqlbrace]{equation}
\begin{aligned}
-\omega^2\rho\bm{u} - \operatorname{div}(2\mu_e\varepsilon(\bm{u}) - \phi\mathbb{I}) &= \bm{f}, \\
i\omega\Big(S_\epsilon + \frac{\alpha^2}{\lambda}\Big)p
- i\omega\frac{\alpha}{\lambda}\phi
- \frac{\kappa}{\mu_f}\Delta p &= g, \\
\phi - \alpha p + \lambda\,\operatorname{tr}(\varepsilon(\bm{u})) &= 0,
\end{aligned}
\label{eq:frequency_domain}
\end{empheq}
For simplicity, the $\omega$-Fourier modes of the displacement, pore pressure, and total pressure are again denoted by $\bm{u}$, $p$, and $\phi$, respectively, and $i=\sqrt{-1}$ denotes the imaginary unit.

The system \eqref{eq:frequency_domain} shall be complemented by appropriate boundary conditions on the displacement and pressure fields. Throughout the analysis, we assume that $\Omega\subset\mathbb{R}^d$ is a
bounded Lipschitz domain. Let $\Gamma_{\bm{u}}$ and $\Gamma_p$ be relatively open
subsets of $\partial\Omega$ such that
\[
\Gamma_{\bm{u}}\cap\Gamma_p=\varnothing,
\qquad
\overline{\Gamma_{\bm{u}}}\cup\overline{\Gamma_p}=\partial\Omega,
\]
and assume that both boundary portions have positive surface measure, namely,
\[
\mathcal{H}^{d-1}(\Gamma_{\bm{u}})>0,
\qquad
\mathcal{H}^{d-1}(\Gamma_p)>0,
\]
where $\mathcal{H}^{d-1}$ denotes the $(d-1)$-dimensional Hausdorff measure. We also assume
$\mathcal H^{d-1}\!\left(\partial\Omega\setminus
(\Gamma_{\bm u}\cup\Gamma_p)\right)=0$.
All boundary descriptions and boundary conditions below are understood up to sets of zero surface measure.
Denoting $\bm{n}$ as the outward normal vector to the boundary, we consider the following mixed-type boundary conditions
\begin{empheq}[left=\empheqlbrace]{equation}
\begin{aligned}
\bm{u} &= \textbf{0} \quad \text{on } \Gamma_{\bm u} \\
(\bm{\sigma}-\alpha p\mathbb{I})\bm{n} &= \bm{g}^{\bm u} \quad \text{on } \Gamma_p, 
\end{aligned}
\label{eq:displacement_bc}
\end{empheq}
and
\begin{empheq}[left=\empheqlbrace]{equation}
\begin{aligned}
    p &= 0 \quad \text{on } \Gamma_p \\
\dfrac{\kappa}{\mu_f}\partial_{\bm n} p &= g^p \quad \text{on } \Gamma_{\bm u}.
\end{aligned}
\label{eq:pressure_bc}
\end{empheq}
In \eqref{eq:displacement_bc}--\eqref{eq:pressure_bc}, the terms $\bm{g}^{\bm u}$ and $g^p$ denote the prescribed traction and normal hydraulic flux amplitudes on $\Gamma_p$ and $\Gamma_{\bm{u}}$, respectively, with the same harmonic time dependence $e^{i\omega t}$.

\subsection{Weak formulation}
Let $L^2(\Omega)$ and $H^1(\Omega)$ denote the standard Sobolev spaces of complex-valued functions, equipped with the inner products
\begin{align*}
    (u,v)_{\Omega}
    &:=(u,v)_{L^2(\Omega)}
    :=\int_{\Omega}u\,\overline{v}\,dx,\\
    (u,v)_{1,\Omega}
    &:=(u,v)_{H^1(\Omega)}
    =(u,v)_{\Omega}
    +\ell^2\sum_{i=1}^{d}
    \left(
    \frac{\partial u}{\partial x_i},
    \frac{\partial v}{\partial x_i}
    \right)_{\Omega},
\end{align*}
where $\overline{v}$ denotes the complex conjugate of $v$, and $\ell$ is a characteristic length of $\Omega$ introduced to ensure dimensional consistency. The corresponding norms and seminorm are denoted by $\|\cdot\|_0$, $\|\cdot\|_1$, and $|\cdot|_1$, respectively. For any measurable subset $\Gamma\subset\partial\Omega$, we denote by $L^2(\Gamma)$ the space of square-integrable complex-valued functions on $\Gamma$, with inner product $\langle\cdot,\cdot\rangle_{\Gamma}$.

With the above notation, let us introduce the functional spaces
\begin{align}
    &\mathbf{H}:=\{\bm{v}:\Omega \to \mathbb{C}^d, \bm{v} \in H^1(\Omega)^d: \bm{v} = \textbf{0} \text{ on } \Gamma_{\bm{u}} \}, \nonumber \\
    &P := \{ q: \Omega \to \mathbb{C}, q \in H^1(\Omega):q=0 \text{ on } \Gamma_p\}, 
    \label{space} \\
    &S:=L^2(\Omega),\nonumber
\end{align}
as well as the product space $\bm{U}:=\mathbf{H}\times P\times S$, equipped with the norm
\begin{equation}
    \|(\bm{v},q,\xi)\|^2_{\bm{U}} := 2\mu_e\|\varepsilon(\bm{v})\|^2_0+\frac{\kappa}{\mu_f\omega\alpha}\|q\|^2_1+\lambda^{-1}\|\xi\|^2_0.
    \label{norm_U}
\end{equation}

For deriving the weak form of \eqref{eq:frequency_domain},  we define the (dimensionless) parameter
\begin{equation}
\theta := \frac{S_\epsilon\lambda}{\alpha} + \alpha = \frac{3\nu(1-\alpha B)}{B(1+\nu)} + \alpha.
\label{theta}
\end{equation}
Moreover, we introduce the sesquilinear forms
\begin{align}
     a_1:\mathbf{H} \times \mathbf{H} \to \mathbb{C},\quad 
    a_1(\bm{u},\bm{v}) &:= 2\mu_e(\varepsilon(\bm{u}),\varepsilon(\bm{v}))_\Omega,
    \nonumber \\
    a_2:P\times P \to \mathbb{C}, \quad a_2(p,q) &:= i\theta\lambda^{-1}(p,q)_{\Omega} 
    + \frac{\kappa}{\mu_f\omega\alpha}(\nabla p,\nabla q)_{\Omega},
    \nonumber \\
    \tilde{b}:\mathbf{H} \times S \to \mathbb{C}, \quad 
    \tilde{b}(\bm{u},\xi) &:= -(\mathrm{div}\,\bm{u},\xi)_{\Omega},
    \nonumber \\
    \tilde{b}^*:S\times \mathbf{H} \to \mathbb{C},\quad 
    \tilde{b}^*(\phi,\bm{v}) &:= -(\phi,\mathrm{div}\,\bm{v})_{\Omega},
    \nonumber \\
     c_1:\mathbf{H} \times \mathbf{H} \to \mathbb{C},\quad 
    c_1(\bm{u},\bm{v}) &:= -\omega^2\rho(\bm{u},\bm{v})_\Omega,
    \nonumber \\
    c_2:P \times S \to \mathbb{C},\quad 
    c_2(p,\xi)&:=-\alpha\lambda^{-1}(p,\xi)_{\Omega},
    \nonumber \\
    c_3:S \times P \to \mathbb{C},\quad 
    c_3(\phi,q)&:= -\lambda^{-1}(\phi,q)_{\Omega},
    \nonumber\\
    d:S \times S \to \mathbb{C}, \quad 
    d(\phi,\xi) &:= \lambda^{-1}(\phi,\xi)_{\Omega}.\nonumber
\end{align}
Multiplying the equations of the system \eqref{eq:frequency_domain} by $\bm{v} \in \mathbf{H},\, q \in P,\, \xi \in S$, respectively, integrating by parts, and imposing the boundary conditions \eqref{eq:displacement_bc}--\eqref{eq:pressure_bc} lead to the following problem: find $\vec{\bm{u}}=(\bm{u},p,\phi)\in \bm{U}$ such that
\begin{equation}
    \langle \mathcal{A}(\vec{\bm{u}}),\vec{\bm{v}} \rangle = \langle \mathcal{F}, \vec{\bm{v}} \rangle
    \label{problem}
\end{equation}
for all $\vec{\bm{v}}=(\bm{v},q,\xi)\in \bm{U}$, where $\mathcal{A}:\bm{U} \to \bm{U}^*$ is defined by
\begin{equation}
    \langle \mathcal{A}(\vec{\bm{u}}),\vec{\bm{v}} \rangle := a_1(\bm{u},\bm{v})
    + a_2(p,q)
    + \tilde{b}^*(\phi,\bm{v})
    - \tilde{b}(\bm{u},\xi)
    + d(\phi,\xi)
    + c_1(\bm{u},\bm{v})
    + c_2(p,\xi)
    + ic_3(\phi,q),
    \label{A}
\end{equation}
and $\mathcal{F} \in \bm{U}^*$ is defined by
\begin{equation}
    \langle \mathcal{F}, \vec{\bm{v}} \rangle :=
    (\bm{f}, \bm{v})_{\Omega} + 
    \langle \bm{g}^{\bm{u}}, \bm{v} \rangle_{\Gamma_p} +
    \frac{1}{\omega\alpha}(g,q)_{\Omega} + 
    \frac{1}{\omega\alpha}\langle g^p,q\rangle _{\Gamma_{\bm{u}}}.
    \label{F}
\end{equation}

We also introduce the operators $\mathcal{B}:\bm{U} \to \bm{U}^*$ and $\mathcal{C}:\bm{U} \to \bm{U}^*$ defined by
\begin{align}
    \langle \mathcal{B}(\vec{\bm{u}}),\vec{\bm{v}} \rangle &:=a_1(\bm{u},\bm{v})
    + a_2(p,q)
    + \tilde{b}^*(\phi,\bm{v})
    - \tilde{b}(\bm{u},\xi)
    + d(\phi,\xi),
    \label{B}\\
    \langle \mathcal{C}(\vec{\bm{u}}),\vec{\bm{v}} \rangle &:=  c_1(\bm{u},\bm{v})
    + c_2(p,\xi)
    + ic_3(\phi,q),
    \label{C}
\end{align}
respectively. With the above operator notation, the variational problem can be expressed compactly as
\begin{equation}
    \mathcal{A}=\mathcal{B}+\mathcal{C}.
    \label{A=B+C}
\end{equation}
The decomposition \eqref{A=B+C} plays a central role in the well-posedness analysis of the weak formulation \eqref{A}--\eqref{F}. Specifically, we shall prove that $\mathcal{B}$ is an isomorphism, whereas $\mathcal{C}$ is a compact operator. Consequently, $\mathcal{A}$ is a Fredholm operator of index zero. By the Fredholm alternative \cite{evans2010partial}, the well-posedness of the variational problem is therefore reduced to establishing the injectivity of $\mathcal{A}$, since every injective Fredholm operator of index zero is necessarily bijective.

\section{Analysis of the continuous problem}

\subsection{Preliminaries}
Before proceeding with the analysis, we recall several fundamental results that will be used throughout the paper. These results can be found, for example, in \cite{brenner2007mathematical,ciarlet2013linear}.

\begin{lemma} [Poincaré inequality]
\label{lem:Poincare}
Assume that \(\mathcal{H}^{d-1}(\Gamma_p)>0\). Then there exists a positive constant $C_{P}$, depending on $\Omega$ and $\Gamma_p$, such that
\begin{equation}
    \|q\|_{0} \leq C_{P}\ell |q|_{1},\quad 
    \|q\|_1^2 \leq (C_P^2+1)\ell^2|q|_1^2,
    \label{Poincare}
\end{equation}
for all $q \in P$.
\end{lemma}

\begin{lemma}[Trace inequality]\label{lem:Trace}
Assume that \(\Omega\) is a bounded Lipschitz domain and let
\(\Gamma\subset \partial\Omega\) be a measurable boundary portion. Then there exists a constant \(C_{\rm tr}>0\), depending only on \(\Omega\) and \(\Gamma\), such that
\begin{equation}
    \|\bm{v}\|_{0,\Gamma}
\le
C_{\rm tr}\|\bm{v}\|_{1},
    \label{Trace}
\end{equation}
for all $\bm{v} \in H^{1}(\Omega)^{d}$.
\end{lemma}

\begin{lemma} [Korn's inequality] 
\label{lem:Korn}
Assume that \(\mathcal{H}^{d-1}(\Gamma_{\bm{u}})>0\). There exists a positive constant $C_{K}$, dependent on $\Omega$, such that
\begin{equation}
    \|\bm{v}\|_{1} \leq C_{K} \|\varepsilon(\bm{v})\|_{0},
    \label{Korn}
\end{equation}
for all $\bm{v} \in \mathbf{H}$.
\end{lemma}
\begin{lemma}[Estimate of divergence by strain tensor] \label{lem:div}
Let $\Omega \subset \mathbb{R}^d$ be a bounded domain and $\bm{v} \in H^1(\Omega)^d$. Then the following inequality holds
\begin{equation}
    \|\nabla \cdot \bm{v}\|_0 \leq \sqrt{d} \|\varepsilon(\bm{v})\|_0.
    \label{div}
\end{equation}
\end{lemma}

\subsection{Well-posedness}
Our well-posedness analysis relies on the spectral properties of the elasticity operator associated with the displacement space $\mathbf H$. We therefore begin by introducing several weighted inner products on $\mathbf H$:
\begin{align*}
(\bm v,\bm w)_{0,\rho}
&:=\rho(\bm v,\bm w)_\Omega,\\
(\bm v,\bm w)_{1,\mu_e}
&:=2\mu_e(\varepsilon(\bm v),\varepsilon(\bm w))_\Omega,\\
(\bm v,\bm w)_{1,\mu_e,\rho}
&:=(\bm v,\bm w)_{1,\mu_e}
+(\bm v,\bm w)_{0,\rho},
\end{align*}
with the corresponding norms denoted by
$\|\cdot\|_{0,\rho}$,
$\|\cdot\|_{1,\mu_e}$,
and
$\|\cdot\|_{1,\mu_e,\rho}$,
respectively.

Since $\Omega$ is a bounded Lipschitz domain and
$\mathcal{H}^{d-1}(\Gamma_{\bm u}) >0$,
Korn's inequality holds on $\mathbf H$, and the embedding
$\mathbf H\hookrightarrow L^2(\Omega)^d$
is compact (see, e.g., \cite{evans2010partial}).
Consequently, the generalized elasticity eigenvalue problem
\[
(\bm w_n,\bm v)_{1,\mu_e}
=
\lambda_n
(\bm w_n,\bm v)_{0,\rho},
\qquad
\forall\,\bm v\in\mathbf H,
\]
admits a countable sequence of positive eigenvalues
\[
0<\lambda_1\le\lambda_2\le\cdots,
\qquad
\lambda_n\rightarrow\infty,
\]
whose associated eigenfunctions
$\{\bm w_n\}_{n\ge1}$
form an orthonormal basis of
$L^2(\Omega)^d$
with respect to
$(\cdot,\cdot)_{0,\rho}$, namely,
\[
(\bm w_n,\bm w_m)_{0,\rho}
=
\delta_{nm}.
\]
Therefore, every
$\bm v\in\mathbf H$
admits the expansion
\[
\bm v
=
\sum_{n=1}^{\infty}
\alpha_n\bm w_n,
\qquad
\alpha_n
=
(\bm v,\bm w_n)_{0,\rho},
\]
and Parseval's identities imply
\[
\|\bm v\|_{0,\rho}^2
=
\sum_{n=1}^{\infty}|\alpha_n|^2,
\qquad
\|\bm v\|_{1,\mu_e}^2
=
\sum_{n=1}^{\infty}\lambda_n|\alpha_n|^2,
\]
and hence
\[
\|\bm v\|_{1,\mu_e,\rho}^2
=
\sum_{n=1}^{\infty}(1+\lambda_n)|\alpha_n|^2.
\]
Similar to \cite{ciarlet2012t}, we make the following non-resonance assumption. 
\begin{assumption}[Non-resonance]\label{ass:nonres}
Throughout this paper, we assume that
\begin{equation}
\omega^2\notin\{\lambda_n:\,n\ge1\},
\label{omega}
\end{equation}
where $\{\lambda_n\}_{n\ge1}$ are the generalized elasticity eigenvalues introduced above.
\end{assumption}
Such an assumption guarantees the well-posedness of the following underlying elasticity problem
\begin{equation}
    -\omega^2\rho(\bm{v},\bm{w}) + 2\mu_e(\varepsilon(\bm{v}),\varepsilon(\bm{w})) = F(\bm{w}), \quad \forall \bm{w}\in \mathbf{H}.
    \label{elasticity_problem}
\end{equation}
In particular, when $F=0$, the unique solution is $\bm{v}=\textbf{0}$.

\begin{lemma}[Continuity]\label{lem:Continuity_A}
There exist positive constants $\eta_1$ and $\eta_2$, depending on the physical parameters and the spatial dimension $d$, such that
\begin{equation*}
\|\mathcal{A}(\vec{\bm{u}})\|_{\bm{U}^*}
\le
\eta_1\|\vec{\bm{u}}\|_{\bm{U}},
\label{continuity_A}
\end{equation*}
and
\begin{equation*}
\|\mathcal{F}\|_{\bm{U}^*}
\le
\eta_2\Big(
\|\bm{f}\|_0
+\|\bm{g}^{\bm{u}}\|_{0,\Gamma_p}
+\|g\|_0
+\|g^p\|_{0,\Gamma_{\bm{u}}}
\Big).
\label{continuity_F}
\end{equation*}
\end{lemma}

\begin{proof}
The estimates follow directly from Korn's inequality \eqref{Korn}, the divergence estimate \eqref{div}, and the Cauchy--Schwarz inequality. Indeed,
\begin{align}
    |\langle\mathcal{A}(\vec{\bm{u}}),\vec{\bm{v}} \rangle |
    &\leq 
    \omega^2\rho\|\bm{u}\|_0\|\bm{v}\|_0 
    +
    2\mu_e\|\varepsilon(\bm{u})\|_0\|\varepsilon(\bm{v})\|_0
    +
    \frac{\theta}{\lambda}\|p\|_0\|q\|_0
    +
    \frac{\kappa}{\mu_f\omega\alpha}\|\nabla p\|_0\|\nabla q\|_0
    \nonumber \\
    &+
    \|\phi\|_0\|\nabla\cdot\bm{v}\|_0
    +
    \|\nabla\cdot\bm{u}\|_0\|\xi\|_0
    +
    \alpha\lambda^{-1}\|p\|_0\|\xi\|_0
    +
    \lambda^{-1}\|\phi\|_0\|q\|_0
    +
    \lambda^{-1}\|\phi\|_0\|\xi\|_0
    \nonumber \\
    &\leq
    \Big( \frac{\omega^2\rho C_K^2}{2\mu_e} + 1 \Big)
    2\mu_e\|\varepsilon(\bm{u})\|_0\|\varepsilon(\bm{v})\|_0
    +
    \Big( \frac{\theta}{\lambda} \Big( \frac{\kappa}{\mu_f\omega\alpha} \Big)^{-1} + 1 \Big)
    \frac{\kappa}{\mu_f\omega\alpha}\|p\|_1\|q\|_1 
    \nonumber \\
    &+
    \frac{\lambda^{1/2}\sqrt{d}}{\sqrt{2\mu_e}}
    \lambda^{-1/2}\|\phi\|_0\sqrt{2\mu_e}\|\varepsilon(\bm{v})\|_0
    +
    \frac{\lambda^{1/2}\sqrt{d}}{\sqrt{2\mu_e}}
    \sqrt{2\mu_e}\|\varepsilon(\bm{u})\|_0\lambda^{-1/2}\|\xi\|_0
    \nonumber \\
    &+
    \alpha\lambda^{-1/2}\Big(\frac{\kappa}{\mu_f\omega\alpha}\Big)^{-1/2}
    \Big(\frac{\kappa}{\mu_f\omega\alpha}\Big)^{1/2}
    \|p\|_1\lambda^{-1/2} \|\xi\|_0
    +
    \lambda^{-1}\|\phi\|_0\|\xi\|_0
    \nonumber \\
    &+
    \lambda^{-1/2}\Big(\frac{\kappa}{\mu_f\omega\alpha}\Big)^{-1/2}
    \lambda^{-1/2}\|\phi\|_0
    \Big(\frac{\kappa}{\mu_f\omega\alpha}\Big)^{1/2} \|q\|_1
    \nonumber \\
    &\le
    \eta_1 \|\vec{\bm{u}}\|_{\bm{U}}\|\vec{\bm{v}}\|_{\bm{U}},
    \nonumber
\end{align}
where
\begin{equation*}
    \eta_1 = 3\max\Bigg\{
    \frac{\omega^2\rho C_K^2}{2\mu_e} + 1,
    \frac{\theta}{\lambda} \Big( \frac{\kappa}{\mu_f\omega\alpha} \Big)^{-1} + 1,
    \frac{\lambda^{1/2}\sqrt{d}}{\sqrt{2\mu_e}},
    \lambda^{-1/2}\Big(\frac{\kappa}{\mu_f\omega\alpha}\Big)^{-1/2}
    \Bigg\}.
    \label{eta_1}
\end{equation*}
For the right-hand side, using the inequalities \eqref{Trace} and \eqref{Korn}, we have
\begin{align*}
    |\langle \mathcal{F},\vec{\bm{v}} \rangle|
    &\le
    \|\bm{f}\|_0\|\bm{v}\|_0
    +
    \|\bm{g}^{\bm{u}}\|_{0,\Gamma_p}\|\bm{v}\|_{0,\Gamma_p}
    +
    \frac{1}{\omega\alpha}\|g\|_0\|q\|_0
    +
    \frac{1}{\omega\alpha}\|g^p\|_{0,\Gamma_{\bm{u}}}\|q\|_{0,\Gamma_{\bm{u}}} \\
    &\le
    C_K\|\bm{f}\|_0\|\varepsilon(\bm{v})\|_0
    +
    C_{tr}C_K\|\bm{g}^{\bm{u}}\|_{0,\Gamma_p} \|\varepsilon(\bm{v})\|_0
    +
    \frac{1}{\omega\alpha}\|g\|_0\|q\|_1
    +
    \frac{C_{tr}}{\omega\alpha} \|g^p\|_{0,\Gamma_{\bm{u}}}\|q\|_1 \\
    &\le
    \frac{C_K\|\bm{f}\|_0 + C_{tr}C_K\|\bm{g}^{\bm{u}}\|_{0,\Gamma_p}}{\sqrt{2\mu_e}} 
    \sqrt{2\mu_e}\|\varepsilon(\bm{v})\|_0
    +
    \frac{\|g\|_0 + C_{tr}\|g^p\|_{0,\Gamma_{\bm{u}}}}{\omega\alpha}\Big(\frac{\kappa}{\mu_f\omega\alpha}\Big)^{-1/2}
    \Big(\frac{\kappa}{\mu_f\omega\alpha}\Big)^{1/2}\|q\|_1 \\
    &\le
    \eta_2\Big(\|\bm{f}\|_0 + \|\bm{g}^{\bm{u}}\|_{0,\Gamma_p} + \|g\|_0 + \|g^p\|_{0,\Gamma_{\bm{u}}}\Big)
    \|\vec{\bm{v}}\|_{\bm{U}}
\end{align*}
where 
\begin{equation*}
    \eta_2 = \max\Bigg\{
    \frac{C_K}{\sqrt{2\mu_e}},
    \frac{C_{tr}C_K}{\sqrt{2\mu_e}},
    \frac{1}{\omega\alpha}\Big(\frac{\kappa}{\mu_f\omega\alpha}\Big)^{-1/2},
    \frac{C_{tr}}{\omega\alpha}\Big(\frac{\kappa}{\mu_f\omega\alpha}\Big)^{-1/2}
    \Bigg\}.
\end{equation*}
\end{proof}

\begin{lemma}
\label{lem:B_bijective}
The operator $\mathcal{B} : \bm{U} \to \bm{U}^*$ defined in \eqref{B} is real-part coercive: there exists a constant $C_B > 0$ such that
\[
\mathrm{Re}\,\langle \mathcal{B}(\vec{\bm{u}}), \vec{\bm{u}}\rangle \ge C_B \|\vec{\bm{u}}\|_{\bm{U}}^2
\quad \forall\vec{\bm{u}}=(\bm{u},p,\phi)\in \bm{U}.
\]
Consequently, $\mathcal{B}$ is bijective.
\end{lemma}

\begin{proof}
Let $\vec{\bm{u}}=(\bm{u},p,\phi)\in \bm{U}$. From \eqref{B} we have
\[
\langle \mathcal{B}(\vec{\bm{u}}), \vec{\bm{u}}\rangle
= a_1(\bm{u},\bm{u}) + a_2(p,p) + \tilde{b}^*(\phi,\bm{u}) - \tilde{b}(\bm{u},\phi) + d(\phi,\phi).
\]
Recall that
\begin{align*}
a_1(\bm{u},\bm{u}) &= 2\mu_e \|\varepsilon(\bm{u})\|_0^2,\\
a_2(p,p) &= i\theta\lambda^{-1}\|p\|_0^2 + \frac{\kappa}{\mu_f\omega\alpha}\|\nabla p\|_0^2,\\
\tilde{b}^*(\phi,\bm{u}) - \tilde{b}(\bm{u},\phi) &= -(\phi,\operatorname{div}\bm{u})_\Omega + (\operatorname{div}\bm{u},\phi)_\Omega
= 2i\,\mathrm{Im}\,(\operatorname{div}\bm{u},\phi)_\Omega,\\
d(\phi,\phi) &= \lambda^{-1}\|\phi\|_0^2.
\end{align*}
The coupling term $\tilde{b}^*-\tilde{b}$ is purely imaginary, so its real part vanishes. Taking the real part gives
\[
\mathrm{Re}\,\langle \mathcal{B}(\vec{\bm{u}}), \vec{\bm{u}}\rangle
= 2\mu_e \|\varepsilon(\bm{u})\|_0^2 + \frac{\kappa}{\mu_f\omega\alpha}\|\nabla p\|_0^2 + \lambda^{-1}\|\phi\|_0^2.
\]
Using the Poincar\'e inequality (\cref{lem:Poincare}) on $p\in P$ yields
\[
\mathrm{Re}\,\langle \mathcal{B}(\vec{\bm{u}}), \vec{\bm{u}}\rangle
\ge \min\!\left\{1,\; \frac{1}{(C_P^2+1)\ell^2}\right\}
\left( 2\mu_e \|\varepsilon(\bm{u})\|_0^2 + \frac{\kappa}{\mu_f\omega\alpha}\|p\|_1^2 + \lambda^{-1}\|\phi\|_0^2 \right)
= C_B \|\vec{\bm{u}}\|_{\bm{U}}^2,
\]
with $C_B := \min\{1,\, (C_P^2+1)^{-1}\ell^{-2}\} > 0$. 

By the Lax--Milgram theorem (see \cite[Proposition 11.29]{brezis2011functional}) applied to the bounded coercive sesquilinear form
associated with $\mathcal{B}$, the operator \(\mathcal{B}:\bm{U}\to \bm{U}^*\) is bijective.
\end{proof}

\begin{lemma}\label{lem:C_compact}
    The operator $\mathcal{C}$, defined in \eqref{C}, is compact.
\end{lemma}

\begin{proof}
We divide $\mathcal C$ into three operators,
\[
\mathcal C=\mathcal K_1+\mathcal K_2+\mathcal K_3
\]
where, for $\vec{\bm{u}},\vec{\bm{v}} \in \bm{U}$, we define
\begin{equation*}
    \langle \mathcal K_1(\bm{u}),\bm{v}\rangle:= -\omega^2\rho\,(\bm{u},\bm{v})_\Omega, \quad
    \langle \mathcal K_2(\phi),q\rangle:= -i\lambda^{-1}(\phi,q)_\Omega,\quad
    \langle \mathcal K_3(p),\xi\rangle:=
    -\alpha\lambda^{-1}(p,\xi)_\Omega.
\end{equation*}
\smallskip
\noindent
\textbf{Compactness of $\mathcal K_1$.}
Define $\mathcal{J}_{\mathbf H}:\mathbf H\hookrightarrow L^2(\Omega)^d$, which is a compact embedding, since $\Omega$ is bounded. We define
\[
\mathcal{M}_1:L^2(\Omega)^d\to \mathbf H^*,
\quad
\langle \mathcal{M}_1(\bm{w}),\bm v\rangle:=
-\omega^2\rho\,(\bm w,\bm v)_\Omega.
\]
By Cauchy--Schwarz and the continuous embedding $\mathbf H\hookrightarrow L^2(\Omega)^d$,
\[
|\langle \mathcal{M}_1(\bm w),\bm v\rangle|
\le
\omega^2\rho\,\|\bm w\|_0\|\bm v\|_0
\le
\frac{\omega^2\rho C_K}{\sqrt{2\mu_e}}\|\bm w\|_0\sqrt{2\mu_e}\|\varepsilon(\bm v)\|_0,
\]
so $\mathcal{M}_1$ is bounded.
Then the displacement contribution in $\mathcal K_1$ is the composition
\[
\mathbf H \xrightarrow{\,\mathcal{J}_{\mathbf H}\,} L^2(\Omega)^d
\xrightarrow{\,\mathcal{M}_1\,} \mathbf H^*
\]
Thus, $\mathcal{K}_1:\bm{U} \to \bm{U}^*$ is compact.

\smallskip
\noindent
\textbf{Compactness of $\mathcal K_2$.} The operator $\mathcal K_2$ maps the $S=L^2(\Omega)$-component to $P^*$. Consider again the compact embedding $\mathcal{J}_P:P\hookrightarrow L^2(\Omega)$. Its dual operator $\mathcal{J}_P^*:(L^2(\Omega))^*\to P^*$ is also compact. Define
\[
\mathcal{M}_2:S\to (L^2(\Omega))^*,\quad
\langle \mathcal{M}_2(\phi),w\rangle:=-i\lambda^{-1}(\phi,w)_\Omega.
\]
Since $\mathcal{M}_2$ is bounded,  $\mathcal K_2$ factors as
\[
S\xrightarrow{\,\mathcal{M}_2\,}(L^2(\Omega))^*
\xrightarrow{\,\mathcal{J}_P^*\,}P^*.
\]
Thus, $\mathcal{K}_2:\bm{U} \to \bm{U}^*$ is compact.

\smallskip
\noindent
\textbf{Compactness of \(\mathcal K_3\).}
The operator \(\mathcal K_3\) maps the \(P\)-component into \(S^*\). Since
\(S=L^2(\Omega)\), we identify \(S^*\) with \(L^2(\Omega)\) through the Riesz
isomorphism
\[
\mathcal R_S:L^2(\Omega)\to S^*,
\qquad
\langle \mathcal R_S(w),\xi\rangle
:=(w,\xi)_\Omega .
\]
Define
\[
\mathcal M_3:L^2(\Omega)\to S^*,
\qquad
\langle \mathcal M_3(w),\xi\rangle
:=
-\alpha\lambda^{-1}(w,\xi)_\Omega .
\]
Equivalently, \(\mathcal M_3=-\alpha\lambda^{-1}\mathcal R_S\), and hence
\(\mathcal M_3\) is bounded. Using the compact embedding
\(\mathcal J_P:P\hookrightarrow L^2(\Omega)\), we obtain the factorization
\[
P
\xrightarrow{\ \mathcal J_P\ }
L^2(\Omega)
\xrightarrow{\ \mathcal M_3\ }
S^* .
\]
Thus, $\mathcal{K}_3:\bm{U} \to \bm{U}^*$ is compact.

Therefore, the sum $\mathcal C=\mathcal K_1+\mathcal K_2+\mathcal K_3$ is compact.
\end{proof}

\begin{lemma}\label{lem:A}
 Under \textbf{Assumption \ref{ass:nonres}}, the operator $\mathcal{A}$ is injective.
\end{lemma}

\begin{proof}
    Let $\vec{\bm{u}} \in \bm{U}$ satisfy $\langle\mathcal{A}(\vec{\bm{u}}),\vec{\bm{v}}\rangle = 0$ for all $\vec{\bm{v}} \in \bm{U}$.
Taking $\vec{\bm{v}} = (\textbf{0},q,0)$, we see that
    \begin{equation}
        a_2(p,q) + ic_3(\phi,q) = 
        i\theta\lambda^{-1}(p,q)_{\Omega} 
        + \frac{\kappa}{\mu_f\omega\alpha}(\nabla p,\nabla q)_{\Omega}
        - i\lambda^{-1}(\phi,q)_{\Omega} = 0.
        \label{v_0_q_0}
    \end{equation}
Choosing $q = p$, then using \eqref{phi} and \eqref{theta}, we obtain
    \begin{equation}
        i\frac{S_{\epsilon}}{\alpha}\|p\|_0^2 + \frac{\kappa}{\mu_f\omega\alpha}\|\nabla p\|_0^2 + i(\mathrm{div}\,\bm{u},p)_{\Omega} = 0.
        \label{v_0_p_0}
    \end{equation}

    Taking $\vec{\bm{v}} = (\bm{u},0,0)$, we have
    \begin{equation*}
        c_1(\bm{u},\bm{u}) + a_1(\bm{u},\bm{u}) + \tilde{b}^*(\phi,\bm{u}) = -\omega^2\rho\|\bm{u}\|_0^2 + 2\mu_e\|\varepsilon(\bm{u})\|_0^2 - (\phi,\mathrm{div}\,\bm{u})_{\Omega} = 0.
    \end{equation*}
By the third variational equation in \eqref{problem}, together with \eqref{phi}, we obtain $\phi=\alpha p-\lambda\,\mathrm{div}\,\bm{u}$ in $L^2(\Omega)$.
Substituting this expression into the preceding identity yields
 
    \begin{equation}
        -\omega^2\rho\|\bm{u}\|_0^2 + 2\mu_e\|\varepsilon(\bm{u})\|_0^2 + \lambda\|\mathrm{div}\,\bm{u}\|_0^2 - \alpha(p, \mathrm{div}\,\bm{u})_{\Omega} = 0.
        \label{v_u_0_0}
    \end{equation}

    Taking the complex conjugate of \eqref{v_0_p_0} and multiplying it by $\alpha i$ gives
    \begin{equation}
        S_{\epsilon}\|p\|_0^2 + i\frac{\kappa}{\mu_f\omega}\|\nabla p\|_0^2 + \alpha(p,\mathrm{div}\,\bm{u})_{\Omega} = 0
        \label{v_0_p_0_i}.
    \end{equation}
    The sum of \eqref{v_u_0_0} and \eqref{v_0_p_0_i} is
    \begin{equation}
        (-\omega^2\rho\|\bm{u}\|_0^2 + 2\mu_e\|\varepsilon(\bm{u})\|_0^2 + \lambda\|\mathrm{div}\,\bm{u}\|_0^2) + S_{\epsilon}\|p\|^2_0 + i\frac{\kappa}{\mu_f\omega}\|\nabla p\|_0^2 = 0.
        \label{v_sum}
    \end{equation}
    Taking the imaginary part of \eqref{v_sum} yields $\frac{\kappa}{\mu_f\omega}\|\nabla p\|_0^2=0$, and hence $\nabla p = 0\,\text{ a.e. in } \Omega$. Since $p\in P$ by definition (see \eqref{space}), the Poincaré inequality (\cref{lem:Poincare}) gives
    \begin{equation*}
        \|p\|_0 \le C_P\ell|p|_1 = C_P\ell\|\nabla p \|_0 = 0,
    \end{equation*}
    so $\|p\|_0 = 0$ and consequently $p = 0$ in $\Omega$.

    Substituting $p = 0$ into \eqref{v_0_q_0}, we get
    \begin{equation*}
        (\phi,q)_{\Omega} = 0 \quad \forall q \in P.
    \end{equation*}
    Since \(C_c^\infty(\Omega)\subset P\) and \(C_c^\infty(\Omega)\) is dense in \(L^2(\Omega)\), the space \(P\) is dense in \(L^2(\Omega)\). We conclude that $\phi = 0$ in $L^2(\Omega)$. Using $p = \phi = 0$ in \eqref{phi}, we obtain $\mathrm{div}\,\bm{u} = 0$.
  Then $\langle\mathcal{A}(\vec{\bm{u}}),\vec{\bm{v}}\rangle = 0$ reduces to
    \begin{equation*}
        -\omega^2\rho(\bm{u},\bm{v})
        +2\mu_e(\varepsilon(\bm{u}),\varepsilon(\bm{v}))
        =0
        \quad\forall \bm{v}\in \mathbf{H}.
    \end{equation*}
    By the spectral assumption \eqref{omega}, problem \eqref{elasticity_problem} is well-posed; the homogeneous elasticity problem has only a trivial solution. Hence, $\bm{u} = \bm{0}$, concluding the proof.
\end{proof}

 
\begin{remark}
Assumption~\ref{ass:nonres} is a sufficient non-resonance condition based on the full elasticity spectrum. As shown in the proof of \cref{lem:A}, only elasticity eigenfunctions satisfying $\nabla\cdot\bm{u}=0$ can contribute to the homogeneous kernel of the coupled problem. A sharper characterization of the resonance set is not pursued here.
\end{remark}


Combining \cref{lem:B_bijective,lem:C_compact,lem:A} with the Fredholm alternative \cite[Theorem 4.2.9]{sauter2011boundary}, we obtain the following main stability result.

\begin{theorem} [Well-posedness] 
\label{thm:Well_Posedness_continuous}
 Under \textbf{Assumption \ref{ass:nonres}}, the problem \eqref{problem} has a unique solution $\vec{\bm{u}}^*\in \bm{U}$ and there exists a positive constant $C$ such that 
\begin{equation*}
    \|\vec{\bm{u}}^*\|_{\bm{U}} \leq C\|\mathcal{F}\|_{\bm{U}^*} \leq C\Big{[}\|\bm{f}\|_0 + \|\bm{g}^{\bm{u}}\|_{0,\Gamma_p} + \|g\|_0 + \|g^p\|_{0,\Gamma_{\bm{u}}} \Big{]}.
    \label{Well_Posedness_continuous}
\end{equation*}
\end{theorem}

\section{Approximation by stabilized finite elements and analysis of the discrete problem}
\label{section4}
\subsection{Stabilized finite element methods}
Let $\{\mathcal T_h\}_{h>0}$ be a shape-regular family of triangulations of $\Omega$. The mesh family is assumed to be boundary-fitted to the prescribed partition of $\partial\Omega$, so that $\Gamma_{\bm u}$ and $\Gamma_p$
are represented by unions of boundary faces. For an element $T\in\mathcal T_h$, we denote its diameter by $h_T$ and introduce
\[
h := \max\{h_T:T\in\mathcal T_h\}
\]
as the characteristic mesh size. 
In a regular mesh family of shapes, any nominal mesh size that is uniformly equivalent to $h$ can be used in the scaling arguments below; the equivalence constants are absorbed into the mesh-independent constants of the analysis.
We consider a mesh family with a fixed upper bound \(h_0>0\) on the mesh size, and we restrict to triangulations satisfying \(0<h\le h_0\). For $k\ge1$, let $\mathcal P_k(T)$ denote the space of polynomials of degree at most $k$ on $T\in\mathcal T_h$. We consider continuous equal-order finite element spaces
\begin{align}
\mathbf H_h
&:=
\Bigl\{
\bm v^h\in C^0(\overline\Omega)^d:
\bm v^h|_T\in\mathcal P_k(T)^d,\;
\forall T\in\mathcal T_h
\Bigr\}
\cap\mathbf H,
\nonumber\\
P_h
&:=
\Bigl\{
q^h\in C^0(\overline\Omega):
q^h|_T\in\mathcal P_k(T),\;
\forall T\in\mathcal T_h
\Bigr\}
\cap P,
\label{finite_element_space}\\
S_h
&:=
\Bigl\{
\xi^h\in C^0(\overline\Omega):
\xi^h|_T\in\mathcal P_k(T),\;
\forall T\in\mathcal T_h
\Bigr\},
\nonumber
\end{align}
and define
\[
\bm U_h:=\mathbf H_h\times P_h\times S_h.
\]
Equal-order discretizations are attractive because of their ease of implementation and computational efficiency, particularly for large-scale three-dimensional problems.  However, neither the pair $\mathbf H_h\times P_h$ nor the pair $\mathbf H_h\times S_h$ satisfies the corresponding discrete inf--sup condition, and the standard Galerkin formulation is not stable. To recover discrete stability, we introduce two complementary stabilization mechanisms.

For the momentum equation in \eqref{eq:frequency_domain}, we define the
elementwise residual
\begin{equation*}
\mathbf R(\bm v^h,\xi^h)
:=
\omega^2\rho\,\bm v^h
+
2\mu_e\operatorname{div}\varepsilon(\bm v^h)
-
\nabla\xi^h,
\label{R}
\end{equation*}
which satisfies
\[
\mathbf R(\bm u,\phi)+\bm f=\mathbf0
\]
for the exact solution. This residual stabilization compensates for the lack of
discrete inf--sup stability of equal-order finite element spaces \cite{oyarzua2016locking,phillips2009overcoming,vuong2016computational}. 
Here, $\mathbf R$ is understood elementwise; namely, for each
$T\in\mathcal T_h$,
\[
\mathbf R(\bm v^h,\xi^h)|_T
=
\omega^2\rho\,\bm v^h|_T
+
2\mu_e\operatorname{div}\varepsilon(\bm v^h|_T)
-
\nabla(\xi^h|_T).
\]

To improve pressure stability without penalizing resolved polynomial modes, we further introduce a local projection stabilization \cite{ganesan2008local}. For each $T\in\mathcal T_h$, let
\[
D_h(T):=\mathbb P_{k-2}(T),
\qquad
\mathbb P_{-1}(T):=\{0\},
\]
and define the discontinuous projection space
\[
D_h
:=
\left\{
\chi^h\in L^2(\Omega):
\chi^h|_T\in D_h(T),
\ \forall T\in\mathcal T_h
\right\}.
\]
Let
$\pi_T:L^2(T)\rightarrow D_h(T)$
be the local $L^2$-orthogonal projection, and define the global projection
$\pi_h:L^2(\Omega)\rightarrow D_h$ elementwise by
\[
(\pi_h w)|_T
=
\pi_T(w|_T),
\qquad
\forall T\in\mathcal T_h.
\]
The corresponding fluctuation operator is
\[
\mathcal K_h:=I-\pi_h,
\]
where $I$ denotes the identity operator on $L^2(\Omega)$. Both operators are
applied componentwise to vector-valued functions.

Then, the stabilization operator
$\mathcal S_h:\bm U_h\rightarrow\bm U_h^*$
is defined by
\begin{equation}
\label{S_h}
\begin{aligned}
\langle\mathcal S_h(\vec{\bm u}^h),\vec{\bm v}^h\rangle
=&\
\delta_1
\sum_{T\in\mathcal T_h}
h_T^2
\bigl(
\mathbf R(\bm u^h,\phi^h),
\mathbf R(\bm v^h,\xi^h)
\bigr)_T
\\
&
+\;
\delta_2
\sum_{T\in\mathcal T_h}
\frac{h_T^2}{\mu_f\omega\alpha}
\bigl(
\mathcal K_h\nabla p^h,
\mathcal K_h\nabla q^h
\bigr)_T,
\end{aligned}
\end{equation}
where
$\delta_1>0$
and
$\delta_2\ge0$
are stabilization parameters.

The second stabilization term is particularly effective in the
low-permeability regime. Since the pressure-gradient contribution in the norm
\eqref{norm_U} is weighted by
$\kappa/(\mu_f\omega\alpha)$,
its control deteriorates as $\kappa$ decreases. The local projection
stabilization restores this control while penalizing only the unresolved
component
$\mathcal K_h\nabla p^h$,
thereby preserving the resolved polynomial modes. For $k=1$,
$D_h(T)=\{0\}$ and hence $\mathcal K_h=I$, so that
\eqref{S_h} reduces to the classical Brezzi--Pitk\"aranta stabilization (BP) \cite{brezzi1984stabilization},
\[
s^{\rm BP}(p^h,q^h)
=
\delta_2
\sum_{T\in\mathcal T_h}
\frac{h_T^2}{\mu_f\omega\alpha}
(\nabla p^h,\nabla q^h)_T.
\]

In sum, our stabilized finite element formulation is: find
$
\vec{\bm u}^h=(\bm u^h,p^h,\phi^h)\in\bm U_h
$
such that
\begin{equation}
\label{discrete_problem}
\langle
\mathcal A_h(\vec{\bm u}^h),
\vec{\bm v}^h
\rangle
=
\langle
\mathcal F_h,
\vec{\bm v}^h
\rangle,
\qquad
\forall\,
\vec{\bm v}^h\in\bm U_h,
\end{equation}
where
\begin{equation}
\label{A_h}
\mathcal A_h
=
\mathcal B+\mathcal S_h+\mathcal C,
\qquad
\langle
\mathcal F_h,
\vec{\bm v}^h
\rangle
=
\langle
\mathcal F,
\vec{\bm v}^h
\rangle
-
\delta_1
\sum_{T\in\mathcal T_h}
h_T^2
(\bm f,\mathbf R(\bm v^h,\xi^h))_T.
\end{equation}


\subsection{Selection and scaling of the stabilization parameters}\label{sec:parsel-parameter-selection}

Next, we specify the scaling of the two stabilization parameters in
\eqref{S_h}.  The two terms have distinct purposes and, consequently, require
different selection principles.  The residual parameter $\delta_1$ is scaled
with the characteristic magnitude of the dynamic momentum operator, while
$\delta_2$ is selected to
provide the smallest prescribed control of the LPS-unresolved pressure modes.

\textbf{Momentum residual stabilization parameter.}
The elementwise weight multiplying the momentum-residual inner product in
\eqref{S_h} is $\delta_1h_T^2$.  Since a single global value of $\delta_1$ is
used throughout this work, its selection is based on the representative mesh
size $h$ introduced above; accordingly, we set
\begin{equation*}
 \tau_1:=\delta_1h^2.
 \end{equation*}
Here \(h\) is understood as any representative mesh size uniformly equivalent to \(\max_{T}h_T\); the constant relating the two is absorbed into \(\gamma_1\).
In the numerical experiments, \(h=1/N\) is used.

We first recall the two standard inverse estimates needed below. For a shape-regular mesh family and a fixed polynomial degree, there exist positive
constants $C_I$ and $\tilde C_I$, independent of $h$ and $T$, such that
\begin{equation}
 \begin{aligned}
 h_T\|\operatorname{div}\varepsilon(\bm v^h)\|_{0,T}
 &\leq C_I\|\varepsilon(\bm v^h)\|_{0,T},
 && \bm v^h|_T\in\mathbb P_k(T)^d,
 \\
 h_T\|\nabla z^h\|_{0,T}
 &\leq\tilde C_I\|z^h\|_{0,T},
 && z^h\in\mathbb P_k(T)^m,\quad m\in\{1,d\}.
 \end{aligned}
 \label{inverse}
\end{equation}
In particular, the second estimate, together with
$\|\varepsilon(\bm v^h)\|_{0,T}\leq
\|\nabla\bm v^h\|_{0,T}$, shows that
\begin{equation*}
 \|\operatorname{div}\varepsilon(\bm v^h)\|_{0,T}
 \leq
 C_I\tilde C_I\,h_T^{-2}\|\bm v^h\|_{0,T}.
\end{equation*}
Hence, up to constants depending only on the fixed finite element family and
the mesh shape-regularity class, the elastic and inertial parts of the strong
momentum operator have the respective scales
\begin{equation*}
 \frac{2\mu_e}{h^2}
 \quad\text{and}\quad
 \rho\omega^2.
 \end{equation*}
The inertial contribution is a zeroth-order term and therefore does not carry an inverse power of $h$. We choose $\tau_1$ proportional to the inverse of the
sum of these operator magnitudes,
\begin{equation*}
 \tau_1
 =
 \frac{\gamma_1}
 {2\mu_e/h^2+\rho\omega^2},
\end{equation*}
and hence
\begin{equation}
 \delta_1
 =
 \frac{\gamma_1}
 {2\mu_e+\rho\omega^2h^2}.
 \label{eq:parsel-delta1-choice}
\end{equation}
Here $\gamma_1>0$ is a dimensionless calibration factor, fixed for the chosen
finite element family and mesh shape-regularity class.  The mesh-independent
constants hidden in the two operator-scale estimates are incorporated into
$\gamma_1$; its value used in the computations is reported together with the
other numerical parameters.  In dimensional variables, both terms in the
denominator of~\eqref{eq:parsel-delta1-choice} have units of stress, and hence
$[\delta_1]=[\mathrm{stress}]^{-1}$.

The rule~\eqref{eq:parsel-delta1-choice} exhibits the following asymptotic behavior:
\begin{equation*}
\begin{gathered}
 \begin{aligned}
 \delta_1
 &\simeq
 \frac{\gamma_1}{2\mu_e},
 &&\text{if }\rho\omega^2h^2\ll2\mu_e,
 \\
 \delta_1
 &\simeq
 \frac{\gamma_1}{\rho\omega^2h^2},
 &&\text{if }\rho\omega^2h^2\gg2\mu_e.
 \end{aligned}
\end{gathered}
 \end{equation*}
Moreover, it immediately yields the uniform bound
\begin{equation*}
    \delta_1(h) \leq \frac{\gamma_1}{2\mu_e}.
\end{equation*}
The elastic and inertial contributions are added in the denominator because the scaling is based on their magnitudes. Using a signed difference would make the resulting stabilization weight singular near an elastic resonance. The momentum residual stabilization is not intended to remove a genuine resonance of the underlying frequency-domain problem; the non-resonance assumption therefore remains necessary.

The complete momentum residual also contains $-\nabla\phi^h$. This term does not enter the leading displacement-operator scaling in \eqref{eq:parsel-delta1-choice}, but the scalar case of \eqref{inverse} gives
\begin{equation*}
 h_T\|\nabla\phi^h\|_{0,T}
 \leq\tilde C_I\|\phi^h\|_{0,T}
 =\tilde C_I\lambda^{1/2}
  \bigl(\lambda^{-1/2}\|\phi^h\|_{0,T}\bigr).
\end{equation*}
It is therefore controlled in the continuity analysis, although the resulting continuity constant may depend on $\tilde C_I$ and $\lambda$. Consequently, \eqref{eq:parsel-delta1-choice} is used as an explicit leading-order scaling
rule, rather than as a parameter-robust optimization formula for the complete residual. Throughout this work, $\delta_1$ remains a single global parameter and $h$ is the representative mesh size used in its evaluation.


\textbf{Pressure stabilization parameter.}
The physical pressure diffusion and the LPS term are, respectively,
\begin{align*}
 d_{\kappa,h}(p^h,q^h)
 &:={}
 \frac{\kappa}{\mu_f\omega\alpha}
 (\nabla p^h,\nabla q^h)_\Omega,
 \\
 \delta_2\ell_h(p^h,q^h)
 &:={}
 \delta_2
 \sum_{T\in\mathcal T_h}
 \frac{h_T^2}{\mu_f\omega\alpha}
 (\mathcal K_h\nabla p^h,
  \mathcal K_h\nabla q^h)_T.
\end{align*}
Thus, $\delta_2h_T^2$ has units of permeability and can be interpreted as an artificial permeability acting on the fluctuation component of the discrete pressure gradient. In particular, $\delta_2$ is dimensionless under the scaling in~\eqref{S_h}.

This interpretation is exact only after separating the parts of the pressure gradient that are resolved and not resolved by the local projection. Since $\pi_T$ is the orthogonal projection in $L^2(T)$ and $\mathcal K_h=I-\pi_h$, one has
\begin{align*}
 &\kappa(\nabla p^h,\nabla q^h)_T
 +\delta_2h_T^2
 (\mathcal K_h\nabla p^h,\mathcal K_h\nabla q^h)_T
 \\
 &\qquad={}
 \kappa(\pi_T\nabla p^h,\pi_T\nabla q^h)_T
 +(\kappa+\delta_2h_T^2)
 (\mathcal K_h\nabla p^h,\mathcal K_h\nabla q^h)_T.
\end{align*}
Here, the projected gradient $\pi_T\nabla q^h$
represents the resolved gradient component, while
$(I-\pi_T)\nabla q^h$ represents its fluctuation component. These components may be interpreted heuristically as the low and high spatial frequency parts, respectively, but they are defined by the projection rather than by a Fourier decomposition.

Consequently, the resolved component retains the physical permeability $\kappa$, while only the fluctuation component experiences the augmented coefficient $\kappa+\delta_2h_T^2$.  For $k=1$, $D_h(T)=\{0\}$ and $\mathcal K_h=I$, so the LPS term acts on the complete pressure gradient. However, for $k\geq2$, the scalar identity $\kappa_{\mathrm{eff}}=\kappa+\delta_2h^2$ is not valid on the complete pressure space.
We therefore select the parameter $\delta_2$ from the pressure Schur complement restricted to the fluctuation space defined in~\ref{app:parsel-pressure-decomposition}. Provided that the complementary block $A_c(\kappa)$ defined in \ref{app:parsel-schur-complement} is invertible, let $\mathcal S_f^0(\kappa)$ denote the reduced pressure operator obtained with $\delta_2=0$, where the superscript $0$ refers only to the absence of the second stabilization term, and let $P_h^f$ be the fluctuation pressure space. We set
\begin{equation*}
 \lambda_f^0(\kappa)
 :=
 \min_{0\neq q_f^h\in P_h^f}
 \frac{
  \operatorname{Re}
  \langle\mathcal S_f^0(\kappa)q_f^h,q_f^h\rangle
 }{
  \ell_h(q_f^h,q_f^h)
 }.
 \label{eq:parsel-lambda-f-main}
\end{equation*}
For a prescribed dimensionless Hermitian stability margin $\beta_*>0$, the
pressure stabilization parameter is chosen as
\begin{equation}
 \delta_{2}(\kappa)
 =
 \left[\beta_*-\lambda_f^0(\kappa)\right]_+,
 \qquad
 [a]_+:=\max\{a,0\}.
 \label{eq:parsel-delta2-main}
\end{equation}
As shown in~\ref{app:parsel-delta2-selection}, this is the smallest nonnegative global value that satisfies the prescribed Hermitian stability margin on the fluctuation pressure space. The criterion applies to complex-valued high-order LPS discretizations in two and three dimensions and provides a load-independent bound for the fluctuation component. It is an
energy stability criterion rather than a nodal monotonicity condition.

\begin{remark} \label{rem:delta2-bounded}
The selection rule \eqref{eq:parsel-delta2-main} defines $\delta_2$ for each fixed mesh. For the discrete compactness and convergence analysis, we assume, in addition, that the resulting mesh-dependent values remain uniformly bounded, namely,
\[
    \sup_h \delta_2(h) < \infty.
\]
Equivalently, there exists a constant $C_{\delta_2}>0$, independent of $h$, such that
\[
    0\le \delta_2(h)\le C_{\delta_2}.
\]
\end{remark}

\subsection{Well-posedness of the discrete problem}
This section establishes the well-posedness and stability of the discrete problem \eqref{discrete_problem} based on the operator decomposition \eqref{A_h}. The bound on $\delta_1$ and the assumed mesh-uniform bound on $\delta_2$ are used throughout the following analysis.

First, we introduce the mesh-dependent norm on $\bm{U}_h$:
\begin{equation}
\|\vec{\bm{v}}^h\|_{\bm{U}_h}^2
:=
\|\vec{\bm{v}}^h\|_{\bm{U}}^2
+
\delta_1
\sum_{T\in\mathcal T_h}
h_T^2
\|\mathbf{R}(\bm{v}^h,\xi^h)\|_{0,T}^2
+
\delta_2
\sum_{T\in\mathcal T_h}
\frac{h_T^2}{\mu_f\omega\alpha}
\|\mathcal{K}_h\nabla q^h\|_{0,T}^2.
\label{norm_U_h}
\end{equation}

The corresponding dual norm on $\bm{U}_h^*$ is defined by
\begin{equation*}
\|\mathcal{G}_h\|_{\bm{U}_h^*}
:=
\sup_{0\ne\vec{\bm{v}}^h\in\bm{U}_h}
\frac{|\langle\mathcal{G}_h,\vec{\bm{v}}^h\rangle|}
{\|\vec{\bm{v}}^h\|_{\bm{U}_h}},
\qquad \mathcal{G}_h\in\bm{U}_h^*.
\end{equation*}

To establish the stability and convergence of
\eqref{discrete_problem}--\eqref{A_h}, we assume that the fluctuation operator satisfies the following local approximation property: there exists a constant $C_\pi>0$, independent of $h$, such that
\begin{equation}
\|\mathcal{K}_h w\|_{0,T}
\le
C_\pi h_T^l|w|_{l,T},
\qquad
\forall\,w\in H^l(T),\;
\forall\,T\in\mathcal T_h,\;
0\le l\le k-1.
\label{Assumption2}
\end{equation}
This assumption holds whenever $P_{k-2}(T)\subset D_h(T)$. For $k=1$, the convention $P_{-1}(T)=\{0\}$ yields $\mathcal{K}_h=I$, while for $k>1$ it follows from the Bramble--Hilbert lemma \cite{brenner2007mathematical}.

The following lemma establishes the continuity of the discrete operator $\mathcal{A}_h$.
\begin{lemma}[Continuity of $\mathcal{A}_h$]
\label{lem:continuity_Ah}
There exist two positive constants $\eta_3$ and $\eta_4$ such that, for any $\vec{\bm{u}}^h \in \bm{U}_h$, we have
\begin{equation*}
\|\mathcal{S}_h(\vec{\bm{u}}^h)\|_{\bm{U}_h^*} \leq \eta_3 \|\vec{\bm{u}}^h\|_{\bm{U}} \quad \text{and} \quad \|\mathcal{A}_h(\vec{\bm{u}}^h)\|_{\bm{U}_h^*} \leq \eta_4 \|\vec{\bm{u}}^h\|_{\bm{U}}.
\end{equation*}
Here, $\eta_3$ and $\eta_4$ are independent of $h$.
\end{lemma}

\begin{proof}
Let $\vec{\bm{u}}^h, \vec{\bm{v}}^h \in \bm{U}_h$. Using the inverse inequalities \eqref{inverse}, we have
\begin{align}
&\Bigl|\delta_1 \sum_{T \in \mathcal{T}_h} h_T^2 \bigl(
\omega^2 \rho \bm{u}^h + 2\mu_e \operatorname{div} \varepsilon(\bm{u}^h) - \nabla \phi^h,
\omega^2 \rho \bm{v}^h + 2\mu_e \operatorname{div} \varepsilon(\bm{v}^h) - \nabla \xi^h 
\bigr)_T\Bigr| \nonumber\\
&\leq \delta_1 \sum_{T \in \mathcal{T}_h} h_T^2 \Big(
\omega^2 \rho \|\bm{u}^h\|_{0,T} + 2\mu_e \|\operatorname{div} \varepsilon(\bm{u}^h)\|_{0,T} + \|\nabla \phi^h\|_{0,T} \Big) 
\Big( \omega^2 \rho \|\bm{v}^h\|_{0,T} + 2\mu_e \|\operatorname{div} \varepsilon(\bm{v}^h)\|_{0,T} + \|\nabla \xi^h\|_{0,T} \Big) \nonumber\\
&\leq \delta_1 \sum_{T \in \mathcal{T}_h} \Big( \omega^2 \rho h_T \|\bm{u}^h\|_{0,T} + 2\mu_e C_I \|\varepsilon(\bm{u}^h)\|_{0,T} + \tilde{C}_I \|\phi^h\|_{0,T} \Big) 
\Big( \omega^2 \rho h_T \|\bm{v}^h\|_{0,T} + 2\mu_e C_I \|\varepsilon(\bm{v}^h)\|_{0,T} + \tilde{C}_I \|\xi^h\|_{0,T} \Big) \nonumber\\
&\leq \delta_1 \max \bigl\{ \omega^2 \rho h_0 (2\mu_e)^{-1/2}, (2\mu_e)^{1/2} C_I, \tilde{C}_I \lambda^{1/2} \bigr\}^2 \left[ \sum_{T \in \mathcal{T}_h} \Big( (2\mu_e)^{1/2} \|\bm{u}^h\|_{0,T} + (2\mu_e)^{1/2} \|\varepsilon(\bm{u}^h)\|_{0,T} + \lambda^{-1/2} \|\phi^h\|_{0,T} \Big)^2 \right]^{1/2} \nonumber\\
&\quad \times \left[ \sum_{T \in \mathcal{T}_h} \Big( (2\mu_e)^{1/2} \|\bm{v}^h\|_{0,T} + (2\mu_e)^{1/2} \|\varepsilon(\bm{v}^h)\|_{0,T} + \lambda^{-1/2} \|\xi^h\|_{0,T} \Big)^2 \right]^{1/2} 
\nonumber\\
&\leq 3\delta_1 \max \bigl\{ \omega^2 \rho h_0 (2\mu_e)^{-1/2}, (2\mu_e)^{1/2} C_I, \tilde{C}_I \lambda^{1/2} \bigr\}^2 \left[ \sum_{T \in \mathcal{T}_h} \Big( 2\mu_e \|\bm{u}^h\|_{0,T}^2 + 2\mu_e \|\varepsilon(\bm{u}^h)\|_{0,T}^2 + \lambda^{-1} \|\phi^h\|_{0,T}^2 \Big) \right]^{1/2}\nonumber \\
&\quad \times \left[ \sum_{T \in \mathcal{T}_h} \Big( 2\mu_e \|\bm{v}^h\|_{0,T}^2 + 2\mu_e \|\varepsilon(\bm{v}^h)\|_{0,T}^2 + \lambda^{-1} \|\xi^h\|_{0,T}^2 \Big) \right]^{1/2} \nonumber\\
&= 3\delta_1 \max \bigl\{ \omega^2 \rho h_0 (2\mu_e)^{-1/2}, (2\mu_e)^{1/2} C_I, \tilde{C}_I \lambda^{1/2} \bigr\}^2\Big( 2\mu_e \|\bm{u}^h\|_{0}^2 + 2\mu_e \|\varepsilon(\bm{u}^h)\|_{0}^2 + \lambda^{-1} \|\phi^h\|_{0}^2 \Big)^{1/2}\nonumber \\
&\quad\times \Big( 2\mu_e \|\bm{v}^h\|_{0}^2 + 2\mu_e \|\varepsilon(\bm{v}^h)\|_{0}^2 + \lambda^{-1} \|\xi^h\|_{0}^2 \Big)^{1/2} \nonumber\\
&\leq 3\delta_1 \max \bigl\{ \omega^2 \rho h_0 (2\mu_e)^{-1/2}, (2\mu_e)^{1/2} C_I, \tilde{C}_I \lambda^{1/2} \bigr\}^2  \Big( 2\mu_e C_K^2 \|\varepsilon(\bm{u}^h)\|_{0}^2 + 2\mu_e \|\varepsilon(\bm{u}^h)\|_{0}^2 + \lambda^{-1} \|\phi^h\|_{0}^2 \Big)^{1/2}\nonumber\\
&\quad \times \Big( 2\mu_e C_K^2 \|\varepsilon(\bm{v}^h)\|_{0}^2 + 2\mu_e \|\varepsilon(\bm{v}^h)\|_{0}^2 + \lambda^{-1} \|\xi^h\|_{0}^2 \Big)^{1/2} \nonumber\\
&\leq 3\delta_1 \max \bigl\{ \omega^2 \rho h_0 (2\mu_e)^{-1/2}, (2\mu_e)^{1/2} C_I, \tilde{C}_I \lambda^{1/2} \bigr\}^2 (1+C_K^2) \|\vec{\bm{u}}^h\|_{\bm{U}} \|\vec{\bm{v}}^h\|_{\bm{U}}.\nonumber
\end{align}
On the other hand,
\begin{align}
\Bigl|\delta_2 \sum_{T \in \mathcal{T}_h} \frac{h_T^2}{\mu_f\omega\alpha} (\mathcal{K}_h\nabla p^h, \mathcal{K}_h\nabla q^h)_T \Bigr|
&\leq
\delta_2 \sum_{T \in \mathcal{T}_h} \frac{h_T^2}{\mu_f\omega\alpha} \|\mathcal{K}_h\nabla p^h\|_{0,T} \|\mathcal{K}_h\nabla q^h\|_{0,T} \nonumber\\
&\leq
\delta_2 \sum_{T \in \mathcal{T}_h} \frac{h_T^2}{\mu_f\omega\alpha} \|\nabla p^h\|_{0,T} \|\nabla q^h\|_{0,T} \nonumber\\
&\leq 
\delta_2 \left(\frac{h_0^2}{\kappa}\right)  \frac{\kappa}{\mu_f\omega\alpha} \left( \sum_{T \in \mathcal{T}_h} \|\nabla p^h\|_{0,T}^2 \right)^{1/2} \left( \sum_{T \in \mathcal{T}_h} \|\nabla q^h\|_{0,T}^2 \right)^{1/2}
\nonumber\\
&\leq \delta_2 \left(\frac{h_0^2}{\kappa}\right) \|\vec{\bm{u}}^h\|_{\bm{U}} \|\vec{\bm{v}}^h\|_{\bm{U}}.\nonumber
\end{align}
So, taking $\eta_3 = 3\delta_1 \max \bigl\{ \omega^2 \rho h_0 (2\mu_e)^{-1/2}, (2\mu_e)^{1/2} C_I, \tilde{C}_I \lambda^{1/2} \bigr\}^2 (1+C_K^2) + \delta_2 \left(\frac{h_0^2}{\kappa}\right)$, we arrive at
\begin{equation}
    \|S_h(\vec{\bm{u}}^h)\|_{\bm{U}_h^*} \leq \eta_3 \|\vec{\bm{u}}^h\|_{\bm{U}}.
    \label{S_h_continuity}
\end{equation}
Inequality \eqref{S_h_continuity} proves the continuity of $\mathcal{S}_h$. Combined with the continuity of $\mathcal{A}$ shown in \cref{lem:Continuity_A}, we can conclude that $\mathcal{A}_h$ is also continuous, i.e.,
\begin{equation}
    |\langle \mathcal{A}_h(\vec{\bm{u}}^h), \vec{\bm{v}}^h \rangle|  \leq \eta_4 \|\vec{\bm{u}}^h\|_{\bm{U}} \|\vec{\bm{v}}^h\|_{\bm{U}},
    \label{A_h_continuity}
\end{equation}
with $\eta_4 = \eta_1 + \eta_3$.
\end{proof}

\begin{lemma}
\label{lem:B_h_coercive}
Let $\mathcal{B}_h := \mathcal{B} + \mathcal{S}_h$. For any fixed angular frequency $\omega>0$, there exists a constant $c_0 > 0$, independent of $h$, such that
\begin{equation}
    \mathrm{Re}\langle\mathcal{B}_h(\vec{\bm{u}}^h),\vec{\bm{u}}^h\rangle \ge c_0\|\vec{\bm{u}}^h\|_{\bm{U}_h}^2,
    \quad \forall \vec{\bm{u}}^h\in \bm{U}_h.
\label{eq:Bhhat-coercive}
\end{equation}
\end{lemma}

\begin{proof}
Since
\begin{align*}
    \mathrm{Re}(\widetilde{b}^*(\phi^h,\bm{u}^h)-\widetilde b(\bm{u}^h,\phi^h)) &= \mathrm{Re}(-(\phi^h,\operatorname{div}\bm{u}^h)_\Omega + (\operatorname{div}\bm{u}^h,\phi^h)_\Omega) = 0,\\
    \mathrm{Re}\,a_2(p^h,p^h) &= \frac{\kappa}{\mu_f\omega\alpha}\|\nabla p^h\|_0^2 \ge (C_P^2+1)^{-1}\ell^{-2}\frac{\kappa}{\mu_f\omega\alpha}\|p^h\|_1^2.
\end{align*}
Therefore,
\begin{align*}
\mathrm{Re}\langle\mathcal{B}_h(\vec{\bm{u}}^h),\vec{\bm{u}}^h\rangle \ge & 2\mu_e\|\varepsilon(\bm{u}^h)\|_0^2 + \lambda^{-1}\|\phi^h\|_0^2
+(C_P^2+1)^{-1}\ell^{-2}\frac{\kappa}{\mu_f\omega\alpha}\|p^h\|_1^2 \\
&+\delta_1\sum_{T\in\mathcal T_h} h_T^2\|\mathbf{R}(\bm{u}^h,\phi^h)\|_{0,T}^2
+\delta_2\sum_{T\in\mathcal T_h} \frac{h_T^2}{\mu_f\omega\alpha}\|\mathcal{K}_h\nabla p^h\|_{0,T}^2.
\end{align*}

Combining the above estimates with the definition of $\|\cdot\|_{\bm{U}_h}$ in \eqref{norm_U_h} yields \eqref{eq:Bhhat-coercive}. 
\end{proof}

The injectivity proof of $\mathcal{A}_h$ will not rely on a direct algebraic argument, as in the continuous case, because the presence of the stabilization terms alters the exact cancellation structure of the original variational formulation. Instead, we prove the discrete well-posedness by a compactness argument of Schatz type.

The next lemma derives an a priori estimate for the discrete problem in a mesh-dependent norm, in which the norm of a discrete function is controlled by the residual $\mathcal{A}_h(\vec{\bm{u}}^h)$ together with lower-order $L^2$-type terms.

\begin{lemma} \label{lem:priori_estimate}
    For every $\vec{\bm{u}}^h\in \bm{U}_h$, there exists a constant $C_{\mathrm{pr}}>0$, independent of $h$, such that
    \begin{equation*}
        \|\vec{\bm{u}}^h\|_{\bm{U}_h} \le C_{\mathrm{pr}}\left(
        \sup_{\substack{\vec{\bm{v}}^h \in \bm{U}_h \\
        \vec{\bm{v}}^h \neq \vec{\bm{0}}}}
        \frac{|\langle \mathcal{A}_h(\vec{\bm{u}}^h), \vec{\bm{v}}^h \rangle|}{\|\vec{\bm{v}}^h\|_{\bm{U}_h}}
        + \|\bm{u}^h\|_0 + \|p^h\|_0
        \right).
    \end{equation*}
    In particular, if $\langle\mathcal{A}_h(\vec{\bm{u}}^h), \vec{\bm{v}}^h\rangle = 0$, then
    \begin{equation*}
        \|\vec{\bm{u}}^h\|_{\bm{U}_h} \le C_{\mathrm{pr}}\Big(
        \|\bm{u}^h\|_0+\|p^h\|_0\Big).
    \end{equation*}
\end{lemma}

\begin{proof}
Since
\[
\mathcal{A}_h=\mathcal{B}_h+\mathcal{C},
\]
we have
\[
\langle\mathcal{B}_h (\vec{\bm{u}}^h),\vec{\bm{v}}^h\rangle
=
\langle\mathcal{A}_h (\vec{\bm{u}}^h),\vec{\bm{v}}^h\rangle
-
\langle\mathcal{C} (\vec{\bm{u}}^h),\vec{\bm{v}}^h\rangle.
\]
Taking $\vec{\bm{u}}^h$ as the test function and applying the coercivity of $\mathcal{B}_h$ from \cref{lem:B_h_coercive}, we obtain
\[
c_0\|\vec{\bm{u}}^h\|_{\bm{U}_h}^2
\le
|\langle\mathcal{A}_h(\vec{\bm{u}}^h),\vec{\bm{u}}^h\rangle|
+
|\langle\mathcal{C}(\vec{\bm{u}}^h),\vec{\bm{u}}^h\rangle|.
\]
Moreover,
\[
\begin{aligned}
\langle\mathcal{C}(\vec{\bm{u}}^h),\vec{\bm{u}}^h\rangle
&=
c_1(\bm{u}^h,\bm{u}^h)
+i\,c_3(\phi^h,p^h)
+c_2(p^h,\phi^h) \\
&=
-\omega^2\rho\|\bm{u}^h\|_0^2
-i\lambda^{-1}(\phi^h,p^h)_\Omega
-\alpha\lambda^{-1}(p^h,\phi^h)_\Omega,
\end{aligned}
\]
and therefore
\[
|\langle\mathcal{C}(\vec{\bm{u}}^h),\vec{\bm{u}}^h\rangle|
\le
\omega^2\rho\|\bm{u}^h\|_0^2
+
2\lambda^{-1}\|\phi^h\|_0\,\|p^h\|_0.
\]
To bound $|\langle\mathcal{A}_h (\vec{\bm{u}}^h),\vec{\bm{u}}^h\rangle|$, by using the dual norm, we have
\begin{equation*}
    |\langle\mathcal{A}_h (\vec{\bm{u}}^h),\vec{\bm{u}}^h\rangle| \le
    \|\mathcal{A}_h \vec{\bm{u}}^h\|_{\bm{U}_h^*}
    \|\vec{\bm{u}}^h\|_{\bm{U}_h}
\end{equation*}
By Young's inequality, we have
\begin{align*}
    2\lambda^{-1}\|\phi^h\|_0\|p^h\|_0
    &\le
    \frac{c_0}{4}\lambda^{-1}\|\phi^h\|_0^2
    +
    \frac{4}{c_0\lambda} \|p^h\|_0^2
    \le
    \frac{c_0}{4}\|\vec{\bm{u}}^h\|_{\bm{U}_h}^2 + \frac{4}{c_0\lambda} \|p^h\|_0^2,\\
    \|\mathcal{A}_h \vec{\bm{u}}^h\|_{\bm{U}_h^*}
    \|\vec{\bm{u}}^h\|_{\bm{U}_h} &\le \frac{1}{c_0}\|\mathcal{A}_h \vec{\bm{u}}^h\|_{\bm{U}_h^*}^2 + 
    \frac{c_0}{4}\|\vec{\bm{u}}^h\|_{\bm{U}_h}^2
\end{align*}
Therefore,
\[
\|\vec{\bm{u}}^h\|_{\bm{U}_h}^2
\le
\frac{2}{c_0^2}\|\mathcal{A}_h \vec{\bm{u}}^h\|_{\bm{U}_h^*}^2
+
\frac{2\omega^2\rho}{c_0} \|\bm{u}^h\|_0^2
+
\frac{8}{c_0^2\lambda}\|p^h\|_0^2.
\]
Hence, 
\[
\|\vec{\bm{u}}^h\|_{\bm{U}_h}
\le
C_{\mathrm{pr}}\left(
\|\mathcal{A}_h \vec{\bm{u}}^h\|_{\bm{U}_h^*}
+ \|\bm{u}^h\|_0 + \|p^h\|_0
\right),
\]
where
\[
C_{\mathrm{pr}}
=
\max\left\{
\frac{\sqrt{2}}{c_0},
\sqrt{\frac{2\omega^2\rho}{c_0}},
\sqrt{\frac{8}{c_0^2\lambda}}
\right\}.
\]
This completes the proof.
\end{proof}

For the Schatz-type compactness argument, we impose the following density assumption:
\[
H^2(\Omega)^d\cap \mathbf H
\quad\text{is dense in } \mathbf H.
\]
To pass to the limit in the stabilization term, we then introduce the subspace
\[
\mathbf D :=
\bigl(H^2(\Omega)^d\cap \mathbf H\bigr)\times P\times
\bigl(H^1(\Omega)\cap S\bigr)
\subset \bm{U}.
\]
Since the second component of \(\mathbf D\) is \(P\) itself and \(H^1(\Omega)\cap S=H^1(\Omega)\) is dense in \(S=L^2(\Omega)\), it follows that \(\mathbf D\) is dense in $\bm{U}$.

We denote the Scott--Zhang quasi-interpolation operators by
\begin{equation}
    \mathcal{I}_h^{\mathbf{H}}:\mathbf{H}\to \mathbf{H}_h,
    \quad
    \mathcal{I}_h^P:P\to P_h,
    \quad
    \mathcal{I}_h^S:H^1(\Omega)\to S_h.
    \label{Scott–Zhang1}
\end{equation}
The Scott--Zhang quasi-interpolation operators are chosen to preserve
the corresponding homogeneous Dirichlet boundary conditions and satisfy
the standard local stability and approximation estimates. We define
\begin{equation}
    \mathcal{I}_h:\mathbf{D}\to \bm{U}_h,
    \quad
    \mathcal{I}_h(\bm{v},q,\xi)
    :=
    \bigl(\mathcal{I}_h^\mathbf{H} \bm{v},\mathcal{I}_h^P q,\mathcal{I}_h^S \xi\bigr).
    \label{Scott–Zhang2}
\end{equation}

The following lemma is an estimate of asymptotic consistency for the stabilization term in the dense test space $\mathbf{D}$, and the corresponding constant depends on the additional regularity of the test function chosen. 
\begin{lemma}\label{lem:S_h=0} For every fixed $\vec{\bm{v}}=(\bm{v},q,\xi)\in \mathbf{D}$, there exists a constant $C(\vec{\bm{v}})>0$, independent of $h$, such that \[ \bigl|\langle \mathcal{S}_h(\vec{\bm{u}}^h),\mathcal{I}_h\vec{\bm{v}}\rangle\bigr| \le C(\vec{\bm{v}})\, h\, \|\vec{\bm{u}}^h\|_{\bm{U}_h}, \quad \forall \vec{\bm{u}}^h\in \bm{U}_h. \] In particular, \[ \bigl|\langle \mathcal{S}_h(\vec{\bm{u}}^h),\mathcal{I}_h\vec{\bm{v}}\rangle  \bigr| \to 0 \quad\text{as } h\to 0, \] uniformly for bounded sequences in $\bm{U}_h$. 
\end{lemma}

\begin{proof}
By the definition of $\mathcal{S}_h$ and the Cauchy--Schwarz inequality,
\[
\begin{aligned}
\bigl|\langle \mathcal{S}_h(\vec{\bm{u}}^h),\mathcal{I}_h\vec{\bm{v}}\rangle\bigr|
\le{}&
\delta_1
\Bigl(
\sum_{T\in\mathcal{T}_h}
h_T^2
\|\mathbf{R}(\bm{u}^h,\phi^h)\|_{0,T}^2
\Bigr)^{1/2}
\Bigl(
\sum_{T\in\mathcal{T}_h}
h_T^2
\|\mathbf{R}(\mathcal{I}_h^{\mathbf H}\bm v,\mathcal{I}_h^S\xi)\|_{0,T}^2
\Bigr)^{1/2}
\\
&
+\delta_2
\Bigl(
\sum_{T\in\mathcal{T}_h}
\frac{h_T^2}{\mu_f\omega\alpha}
\|\mathcal K_h\nabla p^h\|_{0,T}^2
\Bigr)^{1/2}
\Bigl(
\sum_{T\in\mathcal{T}_h}
\frac{h_T^2}{\mu_f\omega\alpha}
\|\mathcal K_h\nabla\mathcal I_h^Pq\|_{0,T}^2
\Bigr)^{1/2}.
\end{aligned}
\]
The stabilization-weighted seminorms of $\vec{\bm u}^h$, including the factors $\sqrt{\delta_1}$ and $\sqrt{\delta_2}$, are controlled by $\|\vec{\bm u}^h\|_{\bm U_h}$. Since $\delta_1$ and $\delta_2$ are uniformly bounded with respect to $h$, the factors $\sqrt{\delta_1}$ and $\sqrt{\delta_2}$ can be absorbed into constants independent of $h$. Moreover, by the $L^2$-stability of $\mathcal K_h$ and the $H^1$-stability of the Scott--Zhang interpolation \cite{scott1990finite},
\begin{equation}
\Bigl(
\sum_{T\in\mathcal T_h}
\frac{h_T^2}{\mu_f\omega\alpha}
\|\mathcal K_h\nabla\mathcal I_h^Pq\|_{0,T}^2
\Bigr)^{1/2}
\le
Ch\|q\|_1 .
\label{fourth_factor}
\end{equation}

Next, since
\[
\mathbf R(\mathcal I_h^{\mathbf H}\bm v,\mathcal I_h^S\xi)
=
\omega^2\rho\,\mathcal I_h^{\mathbf H}\bm v
+
2\mu_e\operatorname{div}
\varepsilon(\mathcal I_h^{\mathbf H}\bm v)
-
\nabla\mathcal I_h^S\xi,
\]
it follows that
\[
\sum_{T\in\mathcal T_h}
h_T^2
\|\mathbf R(\mathcal I_h^{\mathbf H}\bm v,\mathcal I_h^S\xi)\|_{0,T}^2
\le
C(I_1+I_2+I_3),
\]
where
\[
I_1
=
\sum_{T}
h_T^2
\|\mathcal I_h^{\mathbf H}\bm v\|_{0,T}^2,
\quad
I_2
=
\sum_{T}
h_T^2
\|\operatorname{div}\varepsilon(\mathcal I_h^{\mathbf H}\bm v)\|_{0,T}^2,
\quad
I_3
=
\sum_{T}
h_T^2
\|\nabla\mathcal I_h^S\xi\|_{0,T}^2.
\]

By the stability of the Scott--Zhang interpolation \cite{scott1990finite},
\begin{equation}
I_1
\le
Ch^2\|\bm v\|_0^2,
\qquad
I_3
\le
Ch^2\|\xi\|_1^2.
\label{I13}
\end{equation}
For $I_2$, let $p_T\in[\mathbb P_1(\omega_T)]^d$ be the Bramble--Hilbert polynomial.
Since $\operatorname{div}\varepsilon(p_T)=0$,
\[
\begin{aligned}
h_T
\|
\operatorname{div}
\varepsilon(\mathcal I_h^{\mathbf H}\bm v)
\|_{0,T}
&=
h_T
\|
\operatorname{div}
\varepsilon(\mathcal I_h^{\mathbf H}\bm v-p_T)
\|_{0,T}
\\
&\le
C
|
\mathcal I_h^{\mathbf H}\bm v-p_T
|_{1,T}
\\
&\le
C
\Bigl(
|
\mathcal I_h^{\mathbf H}\bm v-\bm v
|_{1,T}
+
|
\bm v-p_T
|_{1,T}
\Bigr)
\\
&\le
Ch_T
|\bm v|_{2,\omega_T},
\end{aligned}
\]
where we have used the inverse inequality, the Scott--Zhang approximation property, and the Bramble--Hilbert estimate.
Using the finite overlap of the patches,
\begin{equation}
I_2
\le
Ch^2|\bm v|_2^2.
\label{I2}
\end{equation}
Combining
\eqref{I13}
and
\eqref{I2},
\begin{equation}
\Bigl(
\sum_{T\in\mathcal T_h}
h_T^2
\|
\mathbf R(\mathcal I_h^{\mathbf H}\bm v,\mathcal I_h^S\xi)
\|_{0,T}^2
\Bigr)^{1/2}
\le
Ch
\left(
\|\bm v\|_0
+
|\bm v|_2
+
\|\xi\|_1
\right).
\label{second_factor}
\end{equation}

Substituting
\eqref{fourth_factor}
and
\eqref{second_factor}
into the above estimate gives
\[
\bigl|
\langle
\mathcal S_h(\vec{\bm u}^h),
\mathcal I_h\vec{\bm v}
\rangle
\bigr|
\le
Ch
\left(
\|\bm v\|_0
+
|\bm v|_2
+
\|\xi\|_1
+
\|q\|_1
\right)
\|
\vec{\bm u}^h
\|_{\bm U_h}
=
C(\vec{\bm v})
h
\|
\vec{\bm u}^h
\|_{\bm U_h}.
\]
The desired convergence follows immediately.
\end{proof}

Finally, the injectivity of $\mathcal{A}_h$ is proven in the following lemma.

\begin{lemma}[Injectivity of $\mathcal{A}_h$]
\label{lem:A_h_injectivity}
There exists $h_{\mathcal{A}}>0$ such that, for every $0<h<h_{\mathcal{A}}$, the discrete operator $\mathcal{A}_h$ is injective.
\end{lemma}

\begin{proof}
We argue by contradiction. Suppose the statement is false. Then there exists a sequence $h_n\to 0$ and a sequence of nonzero vectors 
\[
\vec{\bm{u}}^{h_n}=(\bm{u}^{h_n},p^{h_n},\phi^{h_n})\in \bm{U}_{h_n}
\]
such that
\begin{equation}
    \langle\mathcal{A}_{h_n}(\vec{\bm{u}}^{h_n}), \vec{\bm{v}}^{h_n} \rangle=0,
    \quad
    \|\vec{\bm{u}}^{h_n}\|_{\bm{U}_{h_n}}=1,
    \quad
    \forall n\ge 1.
    \label{bounded}
\end{equation}
For simplicity, we write
\[
\vec{\bm{u}}_n :=\vec{\bm{u}}^{h_n}, \quad \bm{u}_n:=\bm{u}^{h_n}, \quad 
p_n:=p^{h_n}, \quad
\phi_n:=\phi^{h_n}
\]

From \eqref{bounded} the sequence $(\vec{\bm{u}}_{n})$ is bounded in $\bm{U}_{h_n}$ and therefore bounded in $\bm{U}$, since $\|\cdot\|_{\bm{U}} \le \|\cdot\|_{\bm{U}_h}$.
Consequently, there exists $\vec{\bm{u}} = (\bm{u},p,\phi) \in \bm{U}$ such that, up to a subsequence,
\[
\bm{u}_{n}\rightharpoonup \bm{u} \quad\text{in } \mathbf{H},
\quad
p_{n}\rightharpoonup p \quad\text{in } P,
\quad
\phi_{n}\rightharpoonup \phi \quad\text{in } S.
\]
Since the embeddings $\mathbf{H} \hookrightarrow L^2(\Omega)^d$ and $P \hookrightarrow L^2(\Omega)$ are compact, we also have strong convergences
\begin{equation}
\bm{u}_n\to \bm{u} \quad \text{in }L^2(\Omega)^d,
\qquad
p_n\to p \quad \text{in }L^2(\Omega).
\label{strong_convergence}
\end{equation}

Now fix an arbitrary test function $\vec{\bm{v}}=(\bm{v},q,\xi)\in \mathbf{D}$. Set $\vec{\bm{v}}_{n}:=\mathcal{I}_{h_n}\vec{\bm{v}}\in \bm{U}_{h_n}$. By the definition of $\mathcal{A}_{h_n}$ \eqref{A_h}, we have
\begin{equation}
    0
    =
    \langle\mathcal{A}_{h_n}(\vec{\bm{u}}_{n}),\vec{\bm{v}}_n\rangle
    =
    \langle\mathcal{A}(\vec{\bm{u}}_{n}),\vec{\bm{v}}_n\rangle
    +
    \langle\mathcal{S}_{h_n}(\vec{\bm{u}}_{n}),\vec{\bm{v}}_n\rangle.
    \label{A_h_n}
\end{equation}

For the stabilization term, \cref{lem:S_h=0}  together with $\|\vec{\bm{u}}_n\|_{\bm{U}_{h_n}} = 1$ yields the uniform estimate 
\begin{equation*}
    |\langle\mathcal{S}_{h_n}(\vec{\bm{u}}_{n}),\vec{\bm{v}}_n\rangle| \le C(\vec{\bm{v}})h_n \to 0 \quad \text{as } n \to \infty.
    \label{S_h_n=0}
\end{equation*}

It remains to study $\langle \mathcal{A}(\vec{\bm{u}}_n), \vec{\bm{v}}_n\rangle$. Since $\vec{\bm{v}} \in \mathbf{D}$, the Scott--Zhang interpolation provides the local approximation estimates
\begin{equation*}
    \|\bm{v}-\mathcal{I}_{h_n}^{\mathbf{H}}\bm{v}\|_{H_1} \le Ch_n\|\bm{v}\|_{H^2}, \quad \|\xi-\mathcal{I}_{h_n}^S\xi\|_0 \le Ch_n\|\xi\|_{H^1}.
\end{equation*}
For the pressure component, we simply use the strong convergence $\|q-\mathcal{I}_{h_n}^Pq\|_{H^1} \to 0$ (without rate). Consequently,
\[
\|\vec{\bm{v}}_n - \vec{\bm{v}}\|_{\bm{U}} \to 0,
\]
hence $\vec{\bm{v}}_n \to \vec{\bm{v}}$ strongly in $\bm{U}$. Moreover, we have
\[
\langle \mathcal{A}(\vec{\bm{u}}_n), \vec{\bm{v}}_n\rangle
= \langle \mathcal{A}(\vec{\bm{u}}_n), \vec{\bm{v}}_n - \vec{\bm{v}} \rangle
+ \langle \mathcal{A}(\vec{\bm{u}}_n), \vec{\bm{v}} \rangle.
\]
The first term vanishes by the weak--strong pairing: $\vec{\bm{u}}_n \rightharpoonup \vec{\bm{u}}$ with the continuity of $\mathcal{A}$ (\cref{lem:Continuity_A}) and $\vec{\bm{v}}_n - \vec{\bm{v}} \to 0$ strongly. The second term converges to $\langle \mathcal{A}(\vec{\bm{u}}), \vec{\bm{v}}\rangle$ because $\vec{\bm{u}}_n \rightharpoonup \vec{\bm{u}}$ in $\bm{U}$. Consequently,
\begin{equation}
    \lim_{n\to\infty} \langle \mathcal{A}(\vec{\bm{u}}_n), \vec{\bm{v}}_n\rangle
    = \langle \mathcal{A}(\vec{\bm{u}}), \vec{\bm{v}}\rangle.
    \label{Auv}
\end{equation}

Combining \eqref{A_h_n}, \eqref{Auv} and \cref{lem:S_h=0}, we obtain
\[
0 = \lim_{n\to\infty} \langle \mathcal{A}_{h_n}(\vec{\bm{u}}_n), \vec{\bm{v}}_n\rangle
= \langle \mathcal{A}(\vec{\bm{u}}), \vec{\bm{v}}\rangle, \quad \forall \vec{\bm{v}} \in \mathbf{D}.
\]
Since $\mathbf{D}$ is dense in $\bm{U}$ and $\mathcal{A}$ is continuous, we conclude $\langle \mathcal{A}(\vec{\bm{u}}), \vec{\bm{v}}\rangle = 0$ for all $\vec{\bm{v}} \in \bm{U}$. Then the injectivity of the continuous operator $\mathcal{A}$ (\cref{lem:A}) forces
\[
\vec{\bm{u}} = (\bm{u}, p, \phi) = \vec{\bm{0}}.
\]

Returning to \eqref{strong_convergence}, we therefore have $\bm{u}_n \to \mathbf{0}$ and $p_n \to 0$ strongly in $L^2$. Finally, apply \cref{lem:priori_estimate} (the special case where $\langle \mathcal{A}_h(\vec{\bm{u}}^h), \vec{\bm{v}}^h \rangle=0$) to obtain
\[
1 = \|\vec{\bm{u}}_n\|_{\bm{U}_{h_n}}
\le C_{\mathrm{pr}} \bigl( \|\bm{u}_n\|_0 + \|p_n\|_0 \bigr) \to 0,
\]
which is a contradiction. Therefore, there exists $h_{\mathcal{A}} > 0$ such that $\mathcal{A}_h$ is injective for all $0 < h < h_{\mathcal{A}}$.
\end{proof}

\begin{lemma}[Uniform discrete inf-sup condition]
\label{lem:uniform_infsup}
There exists $h_{\mathcal{A}}>0$ and a constant $\beta_2>0$ independent of $h$ such that for every $0<h<h_{\mathcal{A}}$,
\begin{equation}
    \beta_2 \|\vec{\bm{u}}^h\|_{\bm{U}_h}
    \le \sup_{\substack{\vec{\bm{v}}^h \in \bm{U}_h \\ \vec{\bm{v}}^h \neq \vec{\mathbf{0}}}}
    \frac{|\langle \mathcal{A}_h(\vec{\bm{u}}^h), \vec{\bm{v}}^h \rangle|}
    {\|\vec{\bm{v}}^h\|_{\bm{U}_h}}
    \qquad \forall \vec{\bm{u}}^h \in \bm{U}_h.
    \label{uniform_infsup}
\end{equation}
\end{lemma}

\begin{proof}
We argue in a similar way to \cref{lem:A_h_injectivity}. Suppose that such $\beta_2$ does not exist. Then there exists a sequence $h_n\to 0$ and nonzero vectors $\vec{\bm{u}}_n \in \bm{U}_{h_n}$ with
\[
\|\vec{\bm{u}}_n\|_{\bm{U}_{h_n}} = 1,
\qquad
\sup_{\substack{\vec{\bm{v}}^h \in \bm{U}_h \\ \vec{\bm{v}}^h \neq \vec{\mathbf{0}}}}
\frac{|\langle \mathcal{A}_{h_n}(\vec{\bm{u}}_n), \vec{\bm{v}}^h\rangle|}
{\|\vec{\bm{v}}^h\|_{\bm{U}_{h_n}}} \to 0.
\]

Following exactly the same compactness argument as in the proof of \cref{lem:A_h_injectivity}, we extract a weak limit $\vec{\bm{u}} = (\bm{u},p,\phi) \in \bm{U}$ and obtain strong $L^2$ convergence $\bm{u}_n \to \bm{u}$, $p_n \to p$. For any $\vec{\bm{v}} \in \mathbf{D}$ and $\vec{\bm{v}}_n := \mathcal{I}_{h_n}\vec{\bm{v}}$, the vanishing of the dual norm gives $\langle \mathcal{A}_{h_n}(\vec{\bm{u}}_n), \vec{\bm{v}}_n\rangle \to 0$. The stabilization term also vanishes by \cref{lem:S_h=0}, hence $\langle \mathcal{A}(\vec{\bm{u}}_n), \vec{\bm{v}}_n\rangle \to 0$. The identical limit argument as in \cref{lem:A_h_injectivity} yields $\langle \mathcal{A}(\vec{\bm{u}}), \vec{\bm{v}}\rangle = 0$ for all $\vec{\bm{v}}\in\mathbf{D}$, and by density, for all $\vec{\bm{v}}\in\bm{U}$. Injectivity of $\mathcal{A}$ (\cref{lem:A}) forces $\vec{\bm{u}} = \vec{\mathbf{0}}$.

Now \cref{lem:priori_estimate} gives
\[
1 = \|\vec{\bm{u}}_n\|_{\bm{U}_{h_n}}
\le C_{\mathrm{pr}} \left(
\sup_{\substack{\vec{\bm{v}}^h \in \bm{U}_h \\ \vec{\bm{v}}^h \neq \vec{\mathbf{0}}}}
\frac{|\langle \mathcal{A}_{h_n}(\vec{\bm{u}}_n), \vec{\bm{v}}^h\rangle|}
{\|\vec{\bm{v}}^h\|_{\bm{U}_{h_n}}}
+ \|\bm{u}_n\|_0 + \|p_n\|_0
\right) \to 0,
\]
which is a contradiction. Hence, $\beta_2$ must exist uniformly in $h$.
\end{proof}

\begin{theorem}[Well-posedness of the discrete problem]
\label{thm:Well_Posedness_discrete}
Assume that $H^2(\Omega)^d\cap\mathbf H$ is dense in $\mathbf H$. There exists $h_{\mathcal{A}}>0$ such that for every $0<h<h_{\mathcal{A}}$,
the discrete problem \eqref{discrete_problem} admits a unique solution $\vec{\bm{u}}^{h*} \in \bm{U}_h$. Moreover, there exists a constant $C_4>0$ independent of $h$ such that
\begin{equation}
    \|\vec{\bm{u}}_h^{*}\|_{\bm{U}_h}
    \le
    C_4\|\mathcal{F}_h\|_{\bm{U}_h^*}.
    \label{stability_estimate}
\end{equation}
\end{theorem}

\begin{proof}
    \cref{lem:A_h_injectivity} guarantees that $\mathcal{A}_h$ is injective for $0 < h < h_{\mathcal{A}}$. In finite dimensions, injectivity implies bijectivity, hence \eqref{discrete_problem} has a unique solution. The stability estimate \eqref{stability_estimate}  follows from \cref{lem:uniform_infsup}.
\end{proof}

\subsection{Convergence analysis}

This section is devoted to the convergence analysis of the proposed stabilized finite element method. To this end, we assume an additional regularity of the exact solution and introduce the space
\[
\mathbf{W}
:=
H^{k+1}(\Omega)^d
\times
H^{k+1}(\Omega)
\times
H^{k}(\Omega),
\]
where $k$ denotes the polynomial degree of the finite element spaces
defined in \eqref{finite_element_space}. We continue to use the Scott--Zhang quasi-interpolation operator $\mathcal I_h:\mathbf D\to\bm U_h$ introduced in \eqref{Scott–Zhang1} and \eqref{Scott–Zhang2}. By the standard approximation properties of the Scott--Zhang operators and the uniformly bounded overlap of the element patches (see, e.g., \cite[Lemma 1.130]{ern2004theory}), there exists a constant $C_{\mathrm{SZ}}>0$, independent of $h$, such that

\begin{equation}
 \begin{aligned}
 \|\bm v-I_h^{\mathbf H}\bm v\|_{1}
&\le
C_{\mathrm{SZ}}h^k|\bm v|_{k+1}, &&
\forall\,\bm v\in\mathbf H\cap H^{k+1}(\Omega)^d,
\\
\|p-I_h^{P}p\|_{1}
&\le
C_{\mathrm{SZ}}h^k|p|_{k+1},&&
\forall\,p\in P\cap H^{k+1}(\Omega),
\\
\|\phi-I_h^{S}\phi\|_{0}
&\le
C_{\mathrm{SZ}}h^k|\phi|_{k},&&
\forall\,\phi\in H^{k}(\Omega).
 \end{aligned}
 \label{estimate}
\end{equation}

The next lemma establishes a perturbed Galerkin orthogonality relation between the continuous and discrete solutions.

\begin{lemma}[Perturbed Galerkin orthogonality]
\label{lem:Orthogonality}
    Let $\vec{\bm{u}}$ and $\vec{\bm{u}}^h$ be the solutions of \eqref{problem} and \eqref{discrete_problem}, respectively. Assume that $\bm{u} \in \mathbf{H} \cap H^2(\Omega)^d$, $p \in P \cap H^2(\Omega)$, and $\phi \in S \cap H^1(\Omega)$. Then
    \begin{equation*}
        \langle \mathcal{A}_h(\vec{\bm{u}} - \vec{\bm{u}}^h),\vec{\bm{v}}^h\rangle = 
        \delta_2 \sum_{T \in \mathcal{T}_h} \frac{h_T^2}{\mu_f\omega\alpha} \bigl( \mathcal{K}_h\nabla p, \mathcal{K}_h\nabla q^h \bigr)_T
    \end{equation*}
    for all $\vec{\bm{v}}^h \in \bm{U}_h$.
\end{lemma}

\begin{proof}
We have
    \begin{align*}
       \langle \mathcal{A}_h(\vec{\bm{u}} - \vec{\bm{u}}^h), \vec{\bm{v}}^h \rangle &= 
        \langle \mathcal{A}(\vec{\bm{u}}), \vec{\bm{v}}^h \rangle + \langle \mathcal{S}_h(\vec{\bm{u}}), \vec{\bm{v}}^h \rangle - \langle \mathcal{A}_h(\vec{\bm{u}}^h), \vec{\bm{v}}^h \rangle\\ &= \langle\mathcal{S}_h(\vec{\bm{u}}),\vec{\bm{v}}^h\rangle + 
        \delta_1 \sum_{T \in \mathcal{T}_h} h_T^2 \bigl( \bm{f}, \mathbf{R}(\bm{v}^h, \xi^h) \bigr)_T \\
        &=\delta_1 \sum_{T \in \mathcal{T}_h} h_T^2 \bigl( \mathbf{R}(\bm{u}, \phi) + \bm{f}, \mathbf{R}(\bm{v}^h, \xi^h) \bigr)_T 
    + \delta_2 \sum_{T \in \mathcal{T}_h} \frac{h_T^2}{\mu_f\omega\alpha} \bigl( \mathcal{K}_h\nabla p,\mathcal{K}_h \nabla q^h \bigr)_T.
    \end{align*}
     Since $\bm{u} \in \mathbf H\cap H^2(\Omega)^d$ and $\phi \in S\cap H^1(\Omega)$, the strong residual is well defined elementwise and satisfies $\mathbf{R}(\bm{u},\phi)+ \bm{f}=\bm{0}$. The stated identity follows.
     
\end{proof}
\begin{remark}
    Although the operator $\mathcal{S}_h$ was originally defined for discrete functions, it extends in an obvious elementwise manner to any functions for which the integrals in \eqref{S_h} are finite; in particular, it is well-defined for the interpolation error $\vec{\bm{u}} - I_h \vec{\bm{u}}$.
\end{remark}
Next, we derive an error estimate that serves as the key ingredient in the convergence analysis.
\begin{lemma}
\label{lem:error_estimation}
Let $\vec{\bm{u}}$ and $\vec{\bm{u}}^h$ be the solutions
of \eqref{problem} and \eqref{discrete_problem} respectively. In addition, assume that $\vec{\bm{u}} \in \bm{U} \cap \mathbf{W}$. There exists a constant $C>0$, independent of $h$, such that
\[
|\langle \mathcal{A}_h(\vec{\bm{u}}-I_h\vec{\bm{u}}), \vec{\bm{v}}^h\rangle|
\le C h^k \bigl( \|\bm{u}\|_{k+1} + \|p\|_{k+1} + \|\phi\|_k \bigr) \|\vec{\bm{v}}^h\|_{\bm{U}_h}
\quad \forall \vec{\bm{v}}^h\in\bm{U}_h.
\]
\end{lemma}

\begin{proof}
We first decompose
\begin{equation}
\langle \mathcal{A}_h(\vec{\bm{u}}-I_h\vec{\bm{u}}), \vec{\bm{v}}^h\rangle
=
\langle \mathcal{A}(\vec{\bm{u}}-I_h\vec{\bm{u}}), \vec{\bm{v}}^h\rangle
+
\langle \mathcal{S}_h(\vec{\bm{u}}-I_h\vec{\bm{u}}), \vec{\bm{v}}^h\rangle .
\label{4.29}
\end{equation}

For the first term,
\begin{align}
|\langle \mathcal{A}(\vec{\bm{u}}-I_h\vec{\bm{u}}), \vec{\bm{v}}^h\rangle|
&\le
\eta_1
\|\vec{\bm{u}}-I_h\vec{\bm{u}}\|_{\bm{U}}
\|\vec{\bm{v}}^h\|_{\bm{U}}
\nonumber\\
&\le
\eta_1
\|\vec{\bm{u}}-I_h\vec{\bm{u}}\|_{\bm{U}}
\|\vec{\bm{v}}^h\|_{\bm{U}_h}
\label{4.30}\\
&\le
\eta_1C_{\mathrm{SZ}}h^k
\bigl(
\|\bm{u}\|_{k+1}
+
\|p\|_{k+1}
+
\|\phi\|_{k}
\bigr)
\|\vec{\bm{v}}^h\|_{\bm{U}_h},
\nonumber
\end{align}
where we have used the continuity of $\mathcal{A}$ (\cref{lem:Continuity_A}), the inequality
$\|\cdot\|_{\bm U}\le\|\cdot\|_{\bm U_h}$, and the interpolation estimate \eqref{estimate}.

For the second term, using the definition of $\mathcal{S}_h$ in \eqref{S_h},
\[
\begin{aligned}
\langle \mathcal{S}_h(\vec{\bm{u}}-I_h\vec{\bm{u}}), \vec{\bm{v}}^h\rangle
={}&
\delta_1
\sum_{T\in\mathcal{T}_h}
h_T^2
\bigl(
\mathbf{R}(\bm{u}-I_h^{\mathbf{H}}\bm{u},\,\phi-I_h^S\phi),
\mathbf{R}(\bm{v}^h,\xi^h)
\bigr)_T
\\
&
+\delta_2
\sum_{T\in\mathcal{T}_h}
\frac{h_T^2}{\mu_f\omega\alpha}
\bigl(
\mathcal{K}_h\nabla(p-I_h^Pp),
\mathcal{K}_h\nabla q^h
\bigr)_T.
\end{aligned}
\]
We estimate the two stabilization terms separately.

\textbf{Momentum-residual stabilization.}
For any $T\in\mathcal{T}_h$, using the Cauchy--Schwarz inequality, we obtain
\begin{align}
    \Bigl|\delta_1 \sum_{T\in\mathcal{T}_h} &h_T^2
    \bigl( \mathbf{R}(\bm{u} - I_h^{\mathbf{H}}\bm{u},\, \phi - I_h^S\phi),\, \mathbf{R}(\bm{v}^h,\xi^h) \bigr)_T \Bigr|
    \le
    \delta_1 \sum_{T\in\mathcal{T}_h} h_T^2
    \|\mathbf{R}(\bm{u} - I_h^{\mathbf{H}}\bm{u},\, \phi - I_h^S\phi)\|_{0,T}
    \|\mathbf{R}(\bm{v}^h,\xi^h)\|_{0,T}\nonumber \\
    &\le
    \delta_1 \Bigl(\sum_{T\in\mathcal{T}_h} h_T^2 \|\mathbf{R}(\bm{u} - I_h^{\mathbf{H}}\bm{u},\, \phi - I_h^S\phi)\|_{0,T}^2\Bigr)^{1/2}
    \Bigl(\sum_{T\in\mathcal{T}_h} h_T^2 \|\mathbf{R}(\bm{v}^h,\xi^h)\|_{0,T}^2\Bigr)^{1/2} 
    \label{Residual_stabilization}\\
    &\le
    \sqrt{\delta_1} \Bigl(\sum_{T\in\mathcal{T}_h} h_T^2 \|\mathbf{R}(\bm{u} - I_h^{\mathbf{H}}\bm{u},\, \phi - I_h^S\phi)\|_{0,T}^2\Bigr)^{1/2}
    \|\vec{\bm{v}}^h\|_{\bm{U}_h}\nonumber
\end{align}
To bound the term containing the interpolation error, we employ the additional regularity $\bm{u}\in H^{k+1}(\Omega)^d$ and $\phi\in H^k(\Omega)$.
Note that $\mathbf{R}(\bm{v},\xi) = \omega^2\rho\bm{v} + 2\mu_e\operatorname{div}\varepsilon(\bm{v}) - \nabla\xi$.
Using the triangle inequality and standard interpolation estimates, we obtain, for each $T$,
\[
\|\mathbf{R}(\bm{u} - I_h^{\mathbf{H}}\bm{u},\, \phi - I_h^S\phi)\|_{0,T}
\le \omega^2\rho\|\bm{u} - I_h^{\mathbf{H}}\bm{u}\|_{0,T}
+ 2\mu_e\|\operatorname{div}\varepsilon(\bm{u} - I_h^{\mathbf{H}}\bm{u})\|_{0,T}
+ \|\nabla(\phi - I_h^S\phi)\|_{0,T}.
\]
By the local approximation properties of the Scott--Zhang quasi-interpolation operators,
\begin{align*}
\|\bm{u} - I_h^{\mathbf{H}}\bm{u}\|_{0,T} &\le C h_T^{k+1} |\bm{u}|_{k+1, \omega_T},\\
\|\operatorname{div}\varepsilon(\bm{u} - I_h^{\mathbf{H}}\bm{u})\|_{0,T} &\le C h_T^{k-1} |\bm{u}|_{k+1, \omega_T},\\
\|\nabla(\phi - I_h^S\phi)\|_{0,T} &\le C h_T^{k-1} |\phi|_{k, \omega_T},
\end{align*}
where $\omega_T$ denotes the patch of elements adjacent to $T$.
Because $h_T^{k-1}$ is the dominant power (for $k\ge 1$), the entire expression is bounded by $C h_T^{k-1} \bigl(|\bm{u}|_{k+1,\omega_T} + |\phi|_{k,\omega_T}\bigr)$.
Multiplying by $h_T$ and summing all elements gives
\[
\sum_{T\in\mathcal{T}_h} h_T^2 \|\mathbf{R}(\bm{u} - I_h^{\mathbf{H}}\bm{u},\, \phi - I_h^S\phi)\|_{0,T}^2
\le C h^{2k} \bigl( \|\bm{u}\|_{k+1}^2 + \|\phi\|_k^2 \bigr).
\]
Inserting this into \eqref{Residual_stabilization} yields
\begin{equation*}
    \delta_1 \sum_{T\in\mathcal{T}_h} h_T^2
    \bigl( \mathbf{R}(\bm{u} - I_h^{\mathbf{H}}\bm{u},\, \phi - I_h^S\phi),\, \mathbf{R}(\bm{v}^h,\xi^h) \bigr)_T
    \le
    C \sqrt{\delta_1} h^{k} \bigl( \|\bm{u}\|_{k+1} + \|\phi\|_k \bigr) \|\vec{\bm{v}}^h\|_{\bm{U}_h}.
\end{equation*}
 
\textbf{Pressure stabilization.}
Similarly, for the second part of $\mathcal{S}_h$ we have 
\begin{align}
    \Bigl|\delta_2 \sum_{T\in\mathcal{T}_h} \frac{h_T^2}{\mu_f\omega\alpha} \bigl( \mathcal{K}_h\nabla(p - I_h^P p), \mathcal{K}_h\nabla q^h \bigr)_T\Bigr|
    &\le 
    \delta_2 \Bigl(\sum_{T\in\mathcal{T}_h} \frac{h_T^2}{\mu_f\omega\alpha} \|\mathcal{K}_h\nabla(p - I_h^P p)\|_{0,T}^2\Bigr)^{1/2}
    \Bigl(\sum_{T\in\mathcal{T}_h} \frac{h_T^2}{\mu_f\omega\alpha} \|\mathcal{K}_h\nabla q^h\|_{0,T}^2\Bigr)^{1/2}
    \nonumber \\
    &\le
    \sqrt{\delta_2} \Bigl(\sum_{T\in\mathcal{T}_h} \frac{h_T^2}{\mu_f\omega\alpha} \|\nabla(p - I_h^P p)\|_{0,T}^2\Bigr)^{1/2}
    \|\vec{\bm{v}}^h\|_{\bm{U}_h}
    \label{Pressure_stabilization} \\
    &\le
    \sqrt{\delta_2} \Bigl(\sum_{T\in\mathcal{T}_h} \frac{h_T^2}{\mu_f\omega\alpha} (C h_T^{k} |p|_{k+1, \omega_T})^2\Bigr)^{1/2}
    \|\vec{\bm{v}}^h\|_{\bm{U}_h}
    \nonumber \\
    &\le
    C \sqrt{\delta_2} h^{k+1} \|p\|_{k+1} \|\vec{\bm{v}}^h\|_{\bm{U}_h}.\nonumber
\end{align}
Combining \eqref{4.29}, \eqref{4.30}, \eqref{Residual_stabilization} and \eqref{Pressure_stabilization} we obtain
\[
|\langle \mathcal{A}_h(\vec{\bm{u}}-I_h\vec{\bm{u}}), \vec{\bm{v}}^h\rangle|
\le \Bigl[ \eta_1 C_{\mathrm{SZ}} h^k\bigl(\|\bm{u}\|_{k+1} + \|p\|_{k+1} + \|\phi\|_k\bigr)
+ C\sqrt{\delta_1} h^{k} \bigl( \|\bm{u}\|_{k+1} + \|\phi\|_k \bigr)
+ C\sqrt{\delta_2} h^{k+1} \|p\|_{k+1} \Bigr] \|\vec{\bm{v}}^h\|_{\bm{U}_h}.
\]
All terms on the right-hand side are of order $h^k$ (the $\delta_2$ term is
$O(h^{k+1})$ and can be absorbed into the $h^k$ estimate by increasing the
constant). Hence, using the uniform bounds on $\delta_1$ and $\delta_2$ and absorbing them into the constant, there exists a constant $C>0$, independent of $h$, such that
\[
|\langle \mathcal{A}_h(\vec{\bm{u}}-I_h\vec{\bm{u}}), \vec{\bm{v}}^h\rangle|
\le C h^k\bigl(\|\bm{u}\|_{k+1} + \|p\|_{k+1} + \|\phi\|_k\bigr) \|\vec{\bm{v}}^h\|_{\bm{U}_h},
\]
with a constant $C$ independent of $h$.
\end{proof}

The next result concerns the theoretical convergence rate for the stabilized finite element scheme \eqref{discrete_problem}.

\begin{theorem}[Convergence]\label{thm:Convergence}
 Under the assumptions of $\vec{\bm{u}}$ and $\vec{\bm{u}}^h$ in \textup{\cref{lem:error_estimation}}, there exists a constant $C_1 > 0$ independent of $h$ such that
\begin{equation}\label{Convergence}
    \|\vec{\bm{u}}-\vec{\bm{u}}^h\|_{\bm{U}}
    \le
    C_1 h^k
    \Bigl(
    \|\bm{u}\|_{k+1}
    +
    \|p\|_{k+1}
    +
    \|\phi\|_k
    \Bigr).
\end{equation}
\end{theorem}

\begin{proof}
We split the error as
\[
\vec{\bm{u}}-\vec{\bm{u}}^h
=
(\vec{\bm{u}}-I_h\vec{\bm{u}}) + (I_h\vec{\bm{u}}-\vec{\bm{u}}^h)
=:
E(\vec{\bm{u}}) + E_h(\vec{\bm{u}},\vec{\bm{u}}^h),
\]
where
\begin{equation*}
    E(\vec{\bm{u}}) =
    (\bm{u}-I_h^\mathbf{H} \bm{u}, p-I_h^P p, \phi-I_h^S\phi),\quad
    E_h(\vec{\bm{u}},\vec{\bm{u}}^h) =
    (I_h^\mathbf{H} \bm{u}-\bm{u}^h, I_h^P p-p^h, I_h^S\phi-\phi^h).
\end{equation*}

Let $\vec{\bm{v}}^h \in \bm{U}_h$ be arbitrary. Using \cref{lem:uniform_infsup} and \cref{lem:Orthogonality} we obtain
\begin{equation*}
\beta_2 \|E_h(\vec{\bm{u}},\vec{\bm{u}}^h)\|_{\bm{U}_h}
\le
\sup_{\substack{\vec{\bm{v}}^h \in \bm{U}_h \\ 
\vec{\bm{v}}^h \neq \vec{\mathbf{0}}}} 
\frac{|\langle \mathcal{A}_h(E_h(\vec{\bm{u}},\vec{\bm{u}}^h)), \vec{\bm{v}}^h \rangle|}
{\|\vec{\bm{v}}^h\|_{\bm{U}_h}}
=
\sup_{\substack{\vec{\bm{v}}^h \in \bm{U}_h \\ 
\vec{\bm{v}}^h \neq \vec{\mathbf{0}}}} 
\frac{\left| 
\langle \mathcal{A}_h(E(\vec{\bm{u}})), \vec{\bm{v}}^h \rangle 
-
\delta_2 \sum\limits_{T \in \mathcal{T}_h} \frac{h_T^2}{\mu_f\omega\alpha} \bigl( \mathcal{K}_h\nabla p, \mathcal{K}_h\nabla q^h \bigr)_T
\right|}{
\|\vec{\bm{v}}^h\|_{\bm{U}_h}
}.
\end{equation*}
Using \cref{lem:error_estimation}, we obtain
\begin{equation*}
    |\langle \mathcal{A}_h(E(\vec{\bm{u}})), \vec{\bm{v}}^h \rangle| 
    \le
    C h^k \Big(
    \|\bm{u}\|_{k+1} + \|p\|_{k+1} + \|\phi\|_k
    \Big)
    \|\vec{\bm{v}}^h\|_{\bm{U}_h}.
\end{equation*}
In addition, by the Cauchy-Schwarz inequality and \eqref{Assumption2}, we obtain
\begin{align*}
\Bigg|\delta_2 \sum\limits_{T \in \mathcal{T}_h} \frac{h_T^2}{\mu_f\omega\alpha} \bigl( \mathcal{K}_h\nabla p, \mathcal{K}_h\nabla q^h \bigr)_T\Bigg|
&\le
\left( \delta_2 
\sum_{T\in\mathcal{T}_h}\frac{h^2_T}{\mu_f\omega\alpha}\|\mathcal{K}_h\nabla p\|_{0,T}^2
\right)^{1/2}
\left( \delta_2
\sum_{T\in\mathcal{T}_h}\frac{h^2_T}{\mu_f\omega\alpha} \|\mathcal{K}_h\nabla q^h\|_{0,T}^2
\right)^{1/2}
\\
&\le
\left( \delta_2 
\sum_{T\in\mathcal{T}_h}\frac{h^2_T}{\mu_f\omega\alpha}C_{\pi}^2h_T^{2k-2}|p|_{k,T}^2
\right)^{1/2}\|\vec{\bm{v}}^h\|_{\bm{U}_h}
\\
&\le
\tilde{C}_2 h^k |p|_k\|\vec{\bm{v}}^h\|_{\bm{U}_h}
\end{align*}
with $\tilde{C}_2 = C_{\pi}(\frac{\delta_2}{\mu_f\omega\alpha})^{\frac{1}{2}}$, which allows us to conclude
\begin{equation*}
    \|E_h(\vec{\bm{u}},\vec{\bm{u}}^h)\|_{\bm{U}_h}
    \le
    \beta_2^{-1}C h^k \Big(
    \|\bm{u}\|_{k+1} + \|p\|_{k+1} + \|\phi\|_k
    \Big)
    +
    \beta_2^{-1} \tilde{C}_2 h^k |p|_k.
\end{equation*}

Finally, by the triangle inequality and the properties \eqref{estimate} with $|\cdot|_k \le \|\cdot\|_k \le \|\cdot\|_{k+1}$, we obtain
\begin{align*}
\|\vec{\bm{u}}-\vec{\bm{u}}^h\|_{\bm{U}}
&\le
\|E(\vec{\bm{u}})\|_{\bm{U}} + \|E_h(\vec{\bm{u}},\vec{\bm{u}}^h)\|_{\bm{U}}
\le
\|E(\vec{\bm{u}})\|_{\bm{U}} + \|E_h(\vec{\bm{u}},\vec{\bm{u}}^h)\|_{\bm{U}_h}
\nonumber\\
&\le
(C_{\text{SZ}}+\beta_2^{-1}C) h^k \Big(
\|\bm{u}\|_{k+1} + \|p\|_{k+1} + \|\phi\|_k
\Big)
+
\beta_2^{-1}\tilde{C}_2 h^k | p|_k \\
& \le
C_1 h^k
    \Bigl(
    \|\bm{u}\|_{k+1}
    +
    \|p\|_{k+1}
    +
    \|\phi\|_k
    \Bigr).
\end{align*}
\end{proof}

\section{Numerical experiments}

This section presents a series of numerical experiments to validate the proposed stabilized equal-order finite element method. We first consider benchmark problems with exact solutions in two and three spatial dimensions to verify the implementation and confirm the theoretical convergence rates established in Section~\ref{section4}.
For a complex-valued field $v$, we write
$v=\operatorname{Re}v+i\operatorname{Im}v$
and display its real and imaginary components
as the ordered pair
$(\operatorname{Re}v,\operatorname{Im}v)$.

We then consider a heterogeneous layered medium and an anatomically realistic brain geometry reconstructed from magnetic resonance images. These examples assess the performance of the method for discontinuous permeability fields and its applicability to magnetic resonance elastography.

Unless otherwise stated, all numerical experiments use
\[
\gamma_1=1,\quad
\beta_*=1,\quad
\mathrm{tol}_{\mathrm{ker}}=10^{-10} ~(\text{defined in Appendix A.1}).
\]
For each discrete problem, the global residual-stabilization parameter $\delta_1$ is evaluated from~\eqref{eq:parsel-delta1-choice}, and the pressure-stabilization parameter $\delta_2$ is computed from~\eqref{eq:parsel-delta2-main}. Thus, neither parameter is tuned against the numerical solution. The first rule depends explicitly on the representative mesh size and the elastic and inertial operator scales, whereas the second also incorporates the permeability, discrete spaces, boundary conditions, and excitation frequency through the reduced pressure operators.

With the normalization of $\ell_h$ adopted here, $\beta_*=1$ specifies a unit upper bound for the reduced inverse operator from the $\mathcal L_f^{-1}$ dual norm to the $\mathcal L_f$ energy norm in~\eqref{eq:parsel-pressure-amplification-bound}.

\subsection{Benchmark problems with exact solutions}
\label{sec:analytical-validation}

We first study benchmark problems with exact solutions in two and three dimensions. These examples follow the setting of Example~1 in \cite{carcamo2024stabilized} and are used to verify the implementation and convergence of the stabilized equal-order discretization.

\subsubsection{Example 1: Two-dimensional benchmark}

We consider the unit square $\Omega^{2D}=(0,1)^2$,
with boundary decomposition
\begin{equation*}
\Gamma_p=\{(x,y)\in\partial\Omega^{2D}:x=0\ \text{or}\ y=1\},
\qquad
\Gamma_{\bm u}=\{(x,y)\in\partial\Omega^{2D}:x=1\ \text{or}\ y=0\}.
\end{equation*}
The exact solution is given by
\begin{align}
\bm{u}
&=
\left(
\mathrm{Re}\,\bm{u}
\quad
\mathrm{Im}\,\bm{u}
\right)
=
\left(
\begin{array}{cc}
(x-1)^2y^2 &
(x-1)(x+2)^2y(y+1)
\\
xy(x-1) &
2x^2y(x-1)
\end{array}
\right),
\nonumber
\\
p
&=
\left(
\mathrm{Re}\,p
\quad
\mathrm{Im}\,p
\right)
=
\left(
\sin\frac{\pi x}{2}\cos\frac{\pi y}{2}
\quad
(1-\cos\pi x)(1+\cos\pi y)
\right),
\nonumber
\\
\phi
&=
\left(
\mathrm{Re}\,\phi
\quad
\mathrm{Im}\,\phi
\right)
=
\left(
\alpha
\sin\frac{\pi x}{2}
\cos\frac{\pi y}{2}
-
\lambda
\nabla\!\cdot
(\mathrm{Re}\,\bm u)
\quad
\alpha
(1-\cos\pi x)
(1+\cos\pi y)
-
\lambda
\nabla\!\cdot
(\mathrm{Im}\,\bm u)
\right).
\nonumber
\end{align}

\subsubsection{Example 2: Three-dimensional benchmark}

We next consider the unit cube $\Omega^{3D}=(0,1)^3$,
with boundary decomposition
\begin{equation*}
\Gamma_p=\{(x,y,z)\in\partial\Omega^{3D}:x=0\ \text{or}\ x=1\ \text{or} \ y=0\},
\qquad
\Gamma_{\bm u}=\{(x,y,z)\in\partial\Omega^{3D}:y=1\ \text{or}\ z=0\ \text{or} \ z=1\}.
\end{equation*}
The exact solution is given by
\begin{align}
\bm{u}
&=
\left(
\mathrm{Re}\,\bm{u}
\quad
\mathrm{Im}\,\bm{u}
\right)
=
\left(
\begin{array}{cc}
(x-1)^2xy^2(z+2) &
(x-1)x^2y^2(z+1)
\\
x^3y(z+1)^2(x-1)^3 &
(x-1)^2x^3y^3(z-2)
\\
(x-1)xy(z+2)^2 &
x^2y(x-1)(z-2)^2
\end{array}
\right),
\nonumber
\\
p
&=
\left(
\mathrm{Re}\,p
\quad
\mathrm{Im}\,p
\right)
=
\left(
\cos\frac{\pi y}{2}\sin(\pi z)
\quad
\sin\frac{\pi z}{2}
(1+\cos\pi y)
(1+\cos\pi z)
\right),
\nonumber
\\
\phi
&=
\left(
\mathrm{Re}\,\phi
\quad
\mathrm{Im}\,\phi
\right)
=
\left(
\alpha
\cos\frac{\pi y}{2}
\sin(\pi z)
-
\lambda
\nabla\!\cdot
(\mathrm{Re}\,\bm u)
\quad
\alpha
\sin\frac{\pi z}{2}
(1+\cos\pi y)
(1+\cos\pi z)
-
\lambda
\nabla\!\cdot
(\mathrm{Im}\,\bm u)
\right).
\nonumber
\end{align}

The physical parameters used in both benchmark problems are listed in \cref{tab:mms-parameters}. The values of
$E$, $\nu$, $\mu_f$, $\kappa$, $\omega$, and $\rho$
are taken from \cite{carcamo2024stabilized}. The method is implemented in Firedrake with PETSc, and the resulting linear systems are solved using the direct solver MUMPS. Uniformly refined triangular and tetrahedral meshes are used in two and three dimensions, respectively, with $h=1/N$ and $N=2,4,\ldots$ as reported in Tables~\ref{tab:ex1_2d_p1}--\ref{tab:ex1_3d_p2}. The linear computations use the full pressure-gradient stabilization (BP), and the quadratic computations use LPS.

\begin{table}[!ht]
\centering
\caption{Physical parameters used in Example~1 and Example~2.}
\label{tab:mms-parameters}
\begingroup
\setlength{\tabcolsep}{10pt}
\begin{tabular}{lllllllll}
\toprule
Parameter &
$E\,[\mathrm{dyn}\,\mathrm{cm}^{-2}]$ &
$\nu$ &
$\mu_f\,[\mathrm{Poise}]$ &
$\kappa\,[\mathrm{cm}^2]$ &
$\omega\,[\mathrm{rad}\,\mathrm{s}^{-1}]$ &
$\rho\,[\mathrm{g}\,\mathrm{cm}^{-3}]$ &
$\alpha$ &
$B$
\\
\midrule
Value &
$10^2$ &
$0.4$ &
$1.0$ &
$0.1$ &
$1.0$ &
$1.0$ &
$0.5$ &
$0.8$
\\
\bottomrule
\end{tabular}
\endgroup
\end{table}

\begin{figure}[!ht]
\centering
\includegraphics[width=1.00\textwidth]{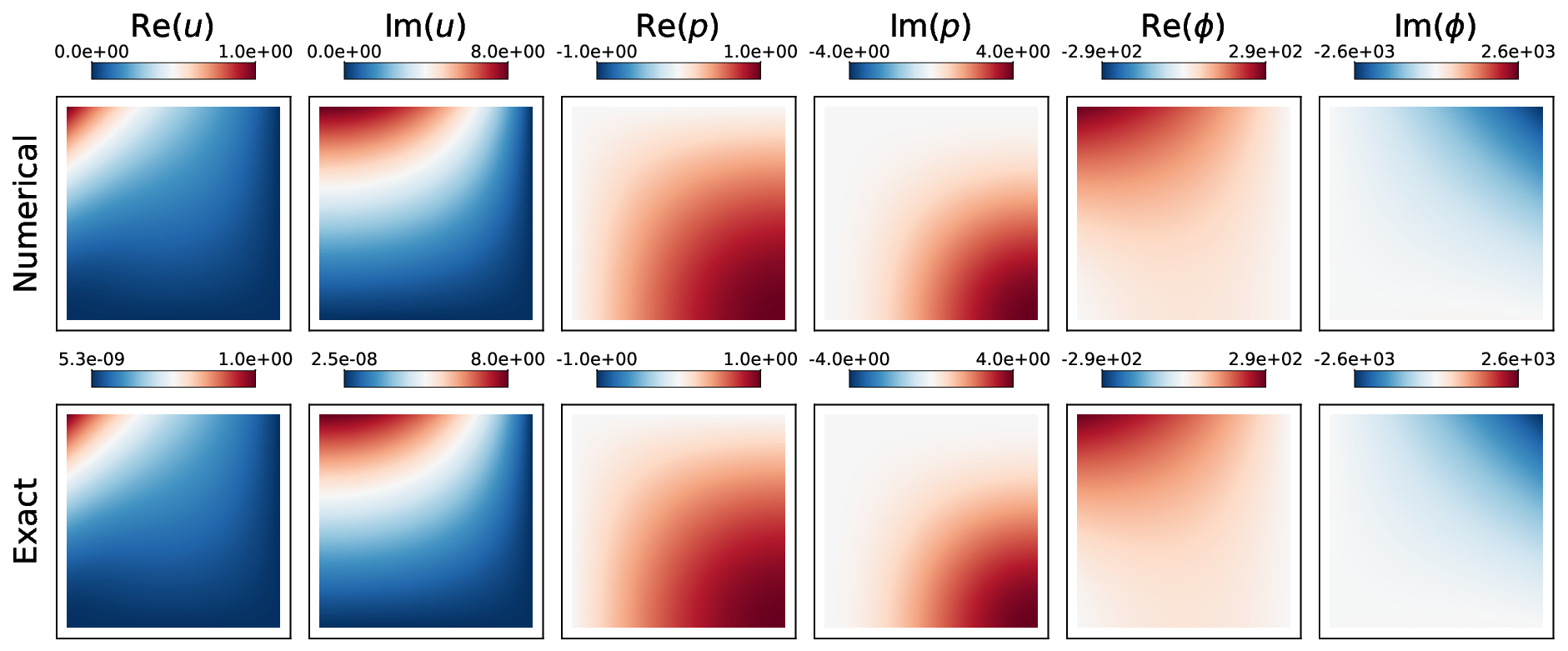}
\caption{Example~1 (2D): Numerical and exact solutions for the real and imaginary parts of $\bm u$, $p$, and $\phi$, obtained with quadratic LPS elements on the mesh $h=1/128$.}
\label{fig:mms-fields}
\end{figure}

For the benchmark data in \cref{tab:mms-parameters}, the computed residual parameter varies only slightly over the mesh sequences,
\[
1.3951\times10^{-2}\leq \delta_1\leq1.4000\times10^{-2}.
\]
The pressure parameter is nonzero only on the two coarsest meshes. In two dimensions, the BP and LPS calculations give $\delta_2=0.8000$ for $h=1/2$ and $\delta_2=0.2000$ for $h=1/4$. In three dimensions, the corresponding values are $0.8667$ and $0.4667$. For $h\leq1/8$, the computed value of $\lambda_f^0(\kappa)$ already exceeds the prescribed margin $\beta_*=1$, and hence \eqref{eq:parsel-delta2-main} gives $\delta_2=0$. 

\begin{figure}[!ht]
\centering
\includegraphics[width=0.835\textwidth]{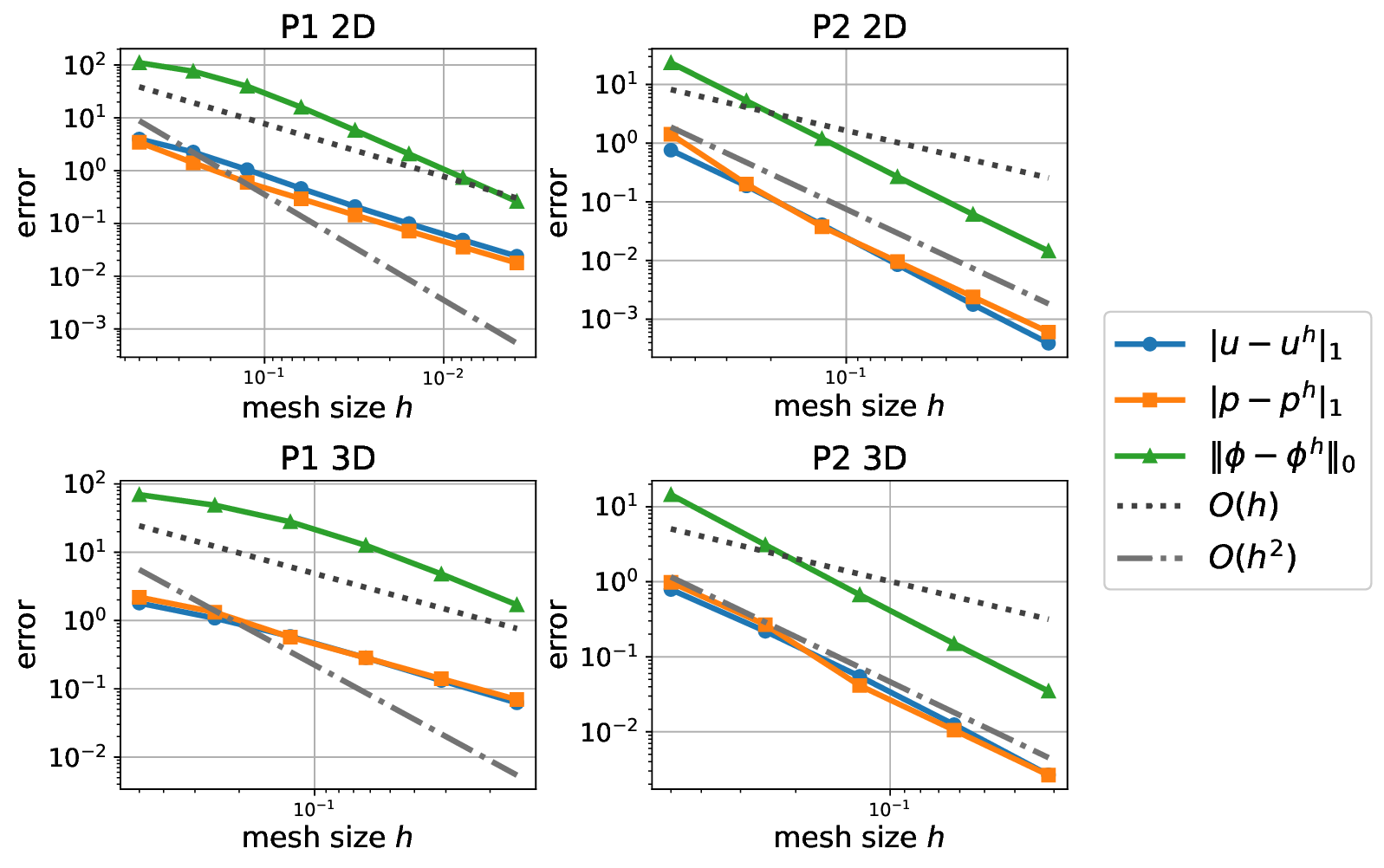}
\caption{Errors of $\bm u$, $p$, and $\phi$ as functions of the mesh size. Top row: Example~1; bottom row: Example~2. Left column: linear elements; right column: quadratic elements. The reference lines indicate the expected $O(h)$ and $O(h^2)$ rates.}
\label{fig:mms-convergence}
\end{figure}

The representative field calculation in \cref{fig:mms-fields} uses quadratic LPS elements with $h=1/128$, $\delta_1=1.4000\times10^{-2}$, and $\delta_2=0$. It compares the numerical and exact solutions for the real and imaginary parts of the displacement $\bm u$, pore pressure $p$, and total pressure $\phi$. The two sets of fields are visually indistinguishable at the plotting scale. \cref{fig:mms-convergence} shows the convergence histories obtained with the mesh-dependent stabilization parameters specified above, and Tables~\ref{tab:ex1_2d_p1}--\ref{tab:ex1_3d_p2} report the corresponding errors and observed rates. In particular, the $H^1$-seminorm errors of $\bm u$ and $p$, together with the $L^2$-error of $\phi$, exhibit the predicted behavior $O(h^k)$ for the discretizations $P_k$. These results verify the implementation and the convergence estimate established in Section~\ref{section4}.

\begin{table}[!ht]
\centering
\caption{Errors and convergence rates for Example 1 using \(P_1\) equal-order elements.}
\label{tab:ex1_2d_p1}
\scriptsize
\resizebox{\textwidth}{!}{%
\begin{tabular}{c cc cc cc cc cc}
\toprule
\(1/h\)
& \(\|\bm{u}-\bm{u}^h\|_0\) & rate
& \(|\bm{u}-\bm{u}^h|_1\) & rate
& \(\|p-p^h\|_0\) & rate
& \(|p-p^h|_1\) & rate
& \(\|\phi-\phi^h\|_0\) & rate \\
\midrule
2   & \(1.071\mathrm{e}{+00}\) & --   & \(3.962\mathrm{e}{+00}\) & --   & \(1.316\mathrm{e}{+00}\) & --   & \(3.436\mathrm{e}{+00}\) & --   & \(1.102\mathrm{e}{+02}\) & --   \\
4   & \(5.413\mathrm{e}{-01}\) & 0.98 & \(2.245\mathrm{e}{+00}\) & 0.82 & \(4.091\mathrm{e}{-01}\) & 1.69 & \(1.407\mathrm{e}{+00}\) & 1.29 & \(7.602\mathrm{e}{+01}\) & 0.54 \\
8   & \(1.913\mathrm{e}{-01}\) & 1.50 & \(1.043\mathrm{e}{+00}\) & 1.11 & \(7.121\mathrm{e}{-02}\) & 2.52 & \(5.961\mathrm{e}{-01}\) & 1.24 & \(3.956\mathrm{e}{+01}\) & 0.94 \\
16  & \(5.615\mathrm{e}{-02}\) & 1.77 & \(4.596\mathrm{e}{-01}\) & 1.18 & \(2.530\mathrm{e}{-02}\) & 1.49 & \(2.924\mathrm{e}{-01}\) & 1.03 & \(1.585\mathrm{e}{+01}\) & 1.32 \\
32  & \(1.523\mathrm{e}{-02}\) & 1.88 & \(2.093\mathrm{e}{-01}\) & 1.14 & \(7.135\mathrm{e}{-03}\) & 1.83 & \(1.439\mathrm{e}{-01}\) & 1.02 & \(5.787\mathrm{e}{+00}\) & 1.45 \\
64  & \(3.980\mathrm{e}{-03}\) & 1.94 & \(9.944\mathrm{e}{-02}\) & 1.07 & \(1.861\mathrm{e}{-03}\) & 1.94 & \(7.151\mathrm{e}{-02}\) & 1.01 & \(2.068\mathrm{e}{+00}\) & 1.49 \\
128 & \(1.020\mathrm{e}{-03}\) & 1.96 & \(4.848\mathrm{e}{-02}\) & 1.04 & \(4.733\mathrm{e}{-04}\) & 1.98 & \(3.570\mathrm{e}{-02}\) & 1.00 & \(7.347\mathrm{e}{-01}\) & 1.49 \\
256 & \(2.589\mathrm{e}{-04}\) & 1.98 & \(2.395\mathrm{e}{-02}\) & 1.02 & \(1.192\mathrm{e}{-04}\) & 1.99 & \(1.784\mathrm{e}{-02}\) & 1.00 & \(2.605\mathrm{e}{-01}\) & 1.50 \\
\bottomrule
\end{tabular}%
}
\vspace{1em}
\end{table}

\begin{table}[!htb]
\centering
\caption{Errors and convergence rates for Example 1 using \(P_2\) equal-order elements.}
\label{tab:ex1_2d_p2}
\scriptsize
\resizebox{\textwidth}{!}{%
\begin{tabular}{c cc cc cc cc cc}
\toprule
\(1/h\)
& \(\|\bm{u}-\bm{u}^h\|_0\) & rate
& \(|\bm{u}-\bm{u}^h|_1\) & rate
& \(\|p-p^h\|_0\) & rate
& \(|p-p^h|_1\) & rate
& \(\|\phi-\phi^h\|_0\) & rate \\
\midrule
2  & \(6.851\mathrm{e}{-02}\) & --   & \(7.597\mathrm{e}{-01}\) & --   & \(2.197\mathrm{e}{-01}\) & --   & \(1.419\mathrm{e}{+00}\) & --   & \(2.345\mathrm{e}{+01}\) & --   \\
4  & \(9.364\mathrm{e}{-03}\) & 2.87 & \(1.867\mathrm{e}{-01}\) & 2.03 & \(1.300\mathrm{e}{-02}\) & 4.08 & \(2.004\mathrm{e}{-01}\) & 2.82 & \(5.265\mathrm{e}{+00}\) & 2.16 \\
8  & \(1.104\mathrm{e}{-03}\) & 3.08 & \(4.074\mathrm{e}{-02}\) & 2.20 & \(6.747\mathrm{e}{-04}\) & 4.27 & \(3.772\mathrm{e}{-02}\) & 2.41 & \(1.193\mathrm{e}{+00}\) & 2.14 \\
16 & \(1.161\mathrm{e}{-04}\) & 3.25 & \(8.487\mathrm{e}{-03}\) & 2.26 & \(8.637\mathrm{e}{-05}\) & 2.97 & \(9.562\mathrm{e}{-03}\) & 1.98 & \(2.667\mathrm{e}{-01}\) & 2.16 \\
32 & \(1.150\mathrm{e}{-05}\) & 3.34 & \(1.779\mathrm{e}{-03}\) & 2.26 & \(1.062\mathrm{e}{-05}\) & 3.02 & \(2.403\mathrm{e}{-03}\) & 1.99 & \(6.114\mathrm{e}{-02}\) & 2.13 \\
64 & \(1.130\mathrm{e}{-06}\) & 3.35 & \(3.891\mathrm{e}{-04}\) & 2.19 & \(1.303\mathrm{e}{-06}\) & 3.03 & \(6.022\mathrm{e}{-04}\) & 2.00 & \(1.448\mathrm{e}{-02}\) & 2.08 \\
\bottomrule
\end{tabular}%
}
\vspace{1em}
\end{table}

\begin{table}[!ht]
\centering
\caption{Errors and convergence rates for Example 2 using \(P_1\) equal-order elements.}
\label{tab:ex1_3d_p1}
\scriptsize
\resizebox{\textwidth}{!}{%
\begin{tabular}{c cc cc cc cc cc}
\toprule
\(1/h\)
& \(\|\bm{u}-\bm{u}^h\|_0\) & rate
& \(|\bm{u}-\bm{u}^h|_1\) & rate
& \(\|p-p^h\|_0\) & rate
& \(|p-p^h|_1\) & rate
& \(\|\phi-\phi^h\|_0\) & rate \\
\midrule
2  & \(2.935\mathrm{e}{-01}\) & --   & \(1.799\mathrm{e}{+00}\) & --   & \(3.118\mathrm{e}{-01}\) & --   & \(2.191\mathrm{e}{+00}\) & --   & \(6.956\mathrm{e}{+01}\) & --   \\
4  & \(1.519\mathrm{e}{-01}\) & 0.95 & \(1.080\mathrm{e}{+00}\) & 0.74 & \(1.347\mathrm{e}{-01}\) & 1.21 & \(1.322\mathrm{e}{+00}\) & 0.73 & \(4.876\mathrm{e}{+01}\) & 0.51 \\
8  & \(6.563\mathrm{e}{-02}\) & 1.21 & \(5.850\mathrm{e}{-01}\) & 0.88 & \(3.675\mathrm{e}{-02}\) & 1.87 & \(5.700\mathrm{e}{-01}\) & 1.21 & \(2.785\mathrm{e}{+01}\) & 0.81 \\
16 & \(2.272\mathrm{e}{-02}\) & 1.53 & \(2.845\mathrm{e}{-01}\) & 1.04 & \(1.495\mathrm{e}{-02}\) & 1.30 & \(2.836\mathrm{e}{-01}\) & 1.01 & \(1.259\mathrm{e}{+01}\) & 1.15 \\
32 & \(6.761\mathrm{e}{-03}\) & 1.75 & \(1.321\mathrm{e}{-01}\) & 1.11 & \(4.824\mathrm{e}{-03}\) & 1.63 & \(1.401\mathrm{e}{-01}\) & 1.02 & \(4.786\mathrm{e}{+00}\) & 1.40 \\
64 & \(1.840\mathrm{e}{-03}\) & 1.88 & \(6.21\mathrm{e}{-02}\) & 1.09 & \(1.344\mathrm{e}{-03}\) & 1.84 & \(6.959\mathrm{e}{-02}\) & 1.01 & \(1.696\mathrm{e}{+00}\) & 1.50 \\
\bottomrule
\end{tabular}%
}
\vspace{1em}
\end{table}

\begin{table}[!htb]
\centering
\caption{Errors and convergence rates for Example 2 using \(P_2\) equal-order elements.}
\label{tab:ex1_3d_p2}
\scriptsize
\resizebox{\textwidth}{!}{%
\begin{tabular}{c cc cc cc cc cc}
\toprule
\(1/h\)
& \(\|\bm{u}-\bm{u}^h\|_0\) & rate
& \(|\bm{u}-\bm{u}^h|_1\) & rate
& \(\|p-p^h\|_0\) & rate
& \(|p-p^h|_1\) & rate
& \(\|\phi-\phi^h\|_0\) & rate \\
\midrule
2  & \(5.204\mathrm{e}{-02}\) & --   & \(7.877\mathrm{e}{-01}\) & --   & \(9.203\mathrm{e}{-02}\) & --   & \(9.812\mathrm{e}{-01}\) & --   & \(1.447\mathrm{e}{+01}\) & --   \\
4  & \(7.477\mathrm{e}{-03}\) & 2.80 & \(2.190\mathrm{e}{-01}\) & 1.85 & \(1.273\mathrm{e}{-02}\) & 2.85 & \(2.677\mathrm{e}{-01}\) & 1.87 & \(3.099\mathrm{e}{+00}\) & 2.22 \\
8  & \(9.326\mathrm{e}{-04}\) & 3.00 & \(5.476\mathrm{e}{-02}\) & 2.00 & \(7.358\mathrm{e}{-04}\) & 4.11 & \(4.133\mathrm{e}{-02}\) & 2.70 & \(6.692\mathrm{e}{-01}\) & 2.21 \\
16 & \(1.105\mathrm{e}{-04}\) & 3.08 & \(1.244\mathrm{e}{-02}\) & 2.14 & \(9.164\mathrm{e}{-05}\) & 3.01 & \(1.051\mathrm{e}{-02}\) & 1.98 & \(1.494\mathrm{e}{-01}\) & 2.16 \\
32 & \(1.247\mathrm{e}{-05}\) & 3.15 & \(2.667\mathrm{e}{-03}\) & 2.22 & \(1.144\mathrm{e}{-05}\) & 3.00 & \(2.644\mathrm{e}{-03}\) & 1.99 & \(3.477\mathrm{e}{-02}\) & 2.10 \\
\bottomrule
\end{tabular}%
}
\vspace{1em}
\end{table}

\subsubsection{Pressure stabilization at low permeability}
\label{sec:lowperm-p2}

We compare the formulation with $\delta_2=0$, the penalty for full gradient BP, and the penalty for fluctuation gradient LPS using the same spaces of equal-order $P_2$ for displacement, pore pressure, and total pressure in Example~1. All three retain the momentum residual term with a common $\delta_1$ on each mesh. Apart from $\kappa\in\{10^{-4},10^{-6},10^{-8}\}$, the physical parameters, exact solution, and boundary decomposition remain unchanged; the manufactured sources and boundary data are recomputed for each permeability. The mesh spacings reported are $h=1/N$ with $N\in\{16,32,64,128\}$. BP and LPS use values of $\delta_2$ computed  from~\eqref{eq:parsel-delta2-main} for their respective stabilization operators separately.

Table~\ref{tab:lowperm-p2-parameters} reports the computed parameters. For $\kappa=10^{-4}$, refinement reduces $\delta_2$ from $0.9872$ to $0.1808$. For the two smaller permeabilities, $\lambda_f^0(\kappa)$ remains small relative to the prescribed margin, so $\delta_2$ stays close to unity. The pressure stabilization then supplies most of the prescribed Hermitian margin. The coefficient $\delta_2h_T^2$ nevertheless decreases on mesh refinement.

\begin{table}[!htbp]
\centering
\caption{Computed stabilization parameters for the $P_2$ comparison in Example~1. BP and LPS values of $\delta_2$ coincide at the displayed precision. All computations use unrounded values.}
\label{tab:lowperm-p2-parameters}
\begin{tabular}{ccccc}
\toprule
$1/h$ & $\delta_1$ & \multicolumn{3}{c}{$\delta_2$ (BP and LPS)} \\
\cmidrule(lr){3-5}
& & $\kappa=10^{-4}$ & $\kappa=10^{-6}$ & $\kappa=10^{-8}$ \\
\midrule
16  & $1.3999\times10^{-2}$ & $0.9872$ & $0.9999$ & $1.0000$ \\
32  & $1.4000\times10^{-2}$ & $0.9488$ & $0.9995$ & $1.0000$ \\
64  & $1.4000\times10^{-2}$ & $0.7952$ & $0.9980$ & $1.0000$ \\
128 & $1.4000\times10^{-2}$ & $0.1808$ & $0.9918$ & $0.9999$ \\
\bottomrule
\end{tabular}
\end{table}

\begin{figure}[!htbp]
\centering
\includegraphics[width=0.9\textwidth]{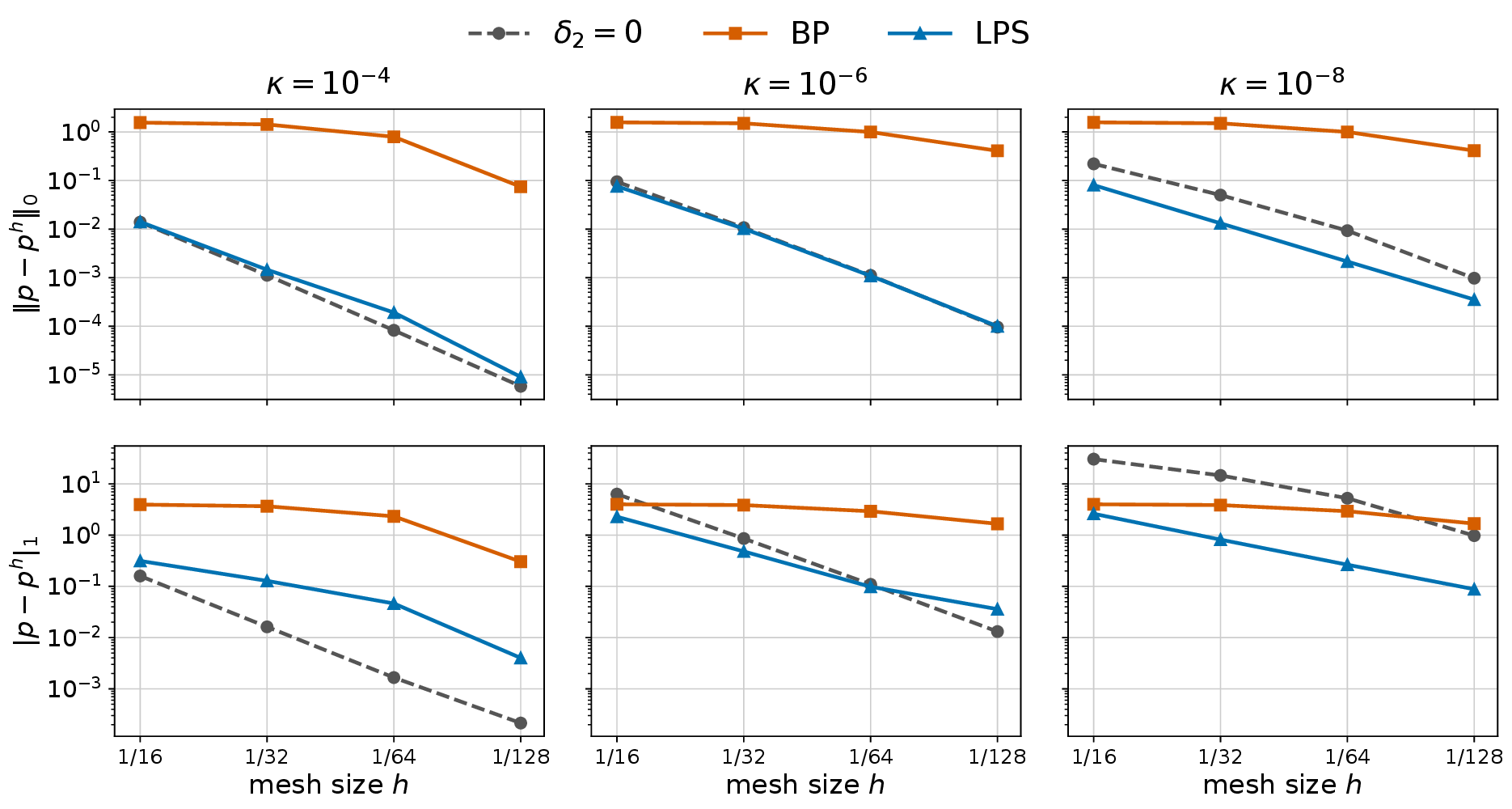}
\caption{Pressure errors for Example~1 with $P_2$ equal-order elements. Top row: $E_0=\|p-p^h\|_0$. Bottom row: $E_1=|p-p^h|_1$. All three formulations retain the momentum residual stabilization; BP and LPS use their respective computed values of $\delta_2$.}
\label{fig:lowperm-p2-errors}
\end{figure}

Figure~\ref{fig:lowperm-p2-errors} compares $E_0=\|p-p^h\|_0$ and $E_1=|p-p^h|_1$. Both norms include the real and imaginary components. The gradient error is unweighted, so its magnitude is not masked as $\kappa$ decreases. For $\kappa=10^{-8}$ and $h=1/128$, LPS reduces $E_1$ from $9.775\times10^{-1}$ to $8.801\times10^{-2}$ and reduces $E_0$ by a factor of $2.75$ relative to $\delta_2=0$. For the prescribed parameter choices, LPS gives substantially smaller pressure errors than BP on all tested meshes.

The benefit relative to $\delta_2=0$ depends on the regime. For $\kappa=10^{-4}$, LPS gives larger gradient errors on all four meshes. For $\kappa=10^{-6}$ at $h=1/128$, it also increases $E_1$ from $1.311\times10^{-2}$ to $3.592\times10^{-2}$. The prescribed Hermitian margin therefore does not guarantee smaller pressure errors.

\begin{figure}[!htbp]
\centering
\includegraphics[width=0.9\textwidth]{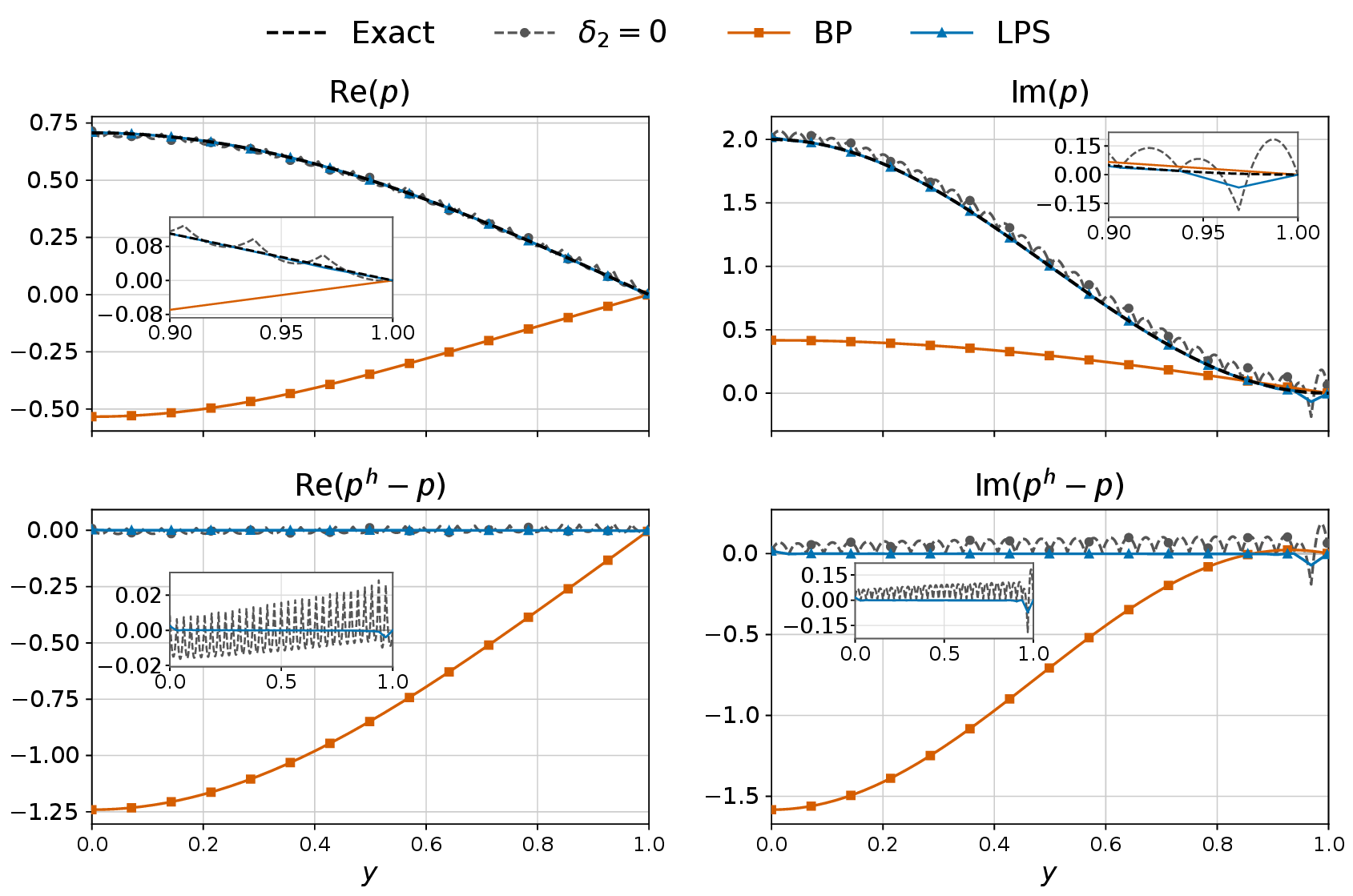}
\caption{Pressure profiles (top) and signed errors $p^h-p$ (bottom) along $x=0.5$ in Example~1 for $\kappa=10^{-8}$ and $h=1/32$. The black dashed curves denote the exact solution. The top insets show $0.9\leq y\leq1$; the bottom insets compare $\delta_2=0$ and LPS on an expanded vertical scale.}
\label{fig:lowperm-p2-profiles}
\end{figure}

Figure~\ref{fig:lowperm-p2-profiles} shows the pressure profiles and signed errors along $x=0.5$ for $\kappa=10^{-8}$ and $h=1/32$. The exact components decrease monotonically, whereas the $\delta_2=0$ approximation exhibits repeated spurious spatial oscillations. With the selected BP parameter, the real component is negative and increases toward zero, opposite to the nonnegative decreasing exact profile. BP also underestimates the imaginary component over most of the line. LPS closely follows the exact profiles in the interior and strongly attenuates the oscillations. A localized imaginary-pressure undershoot remains near $y=1$: at $y=31/32$, the LPS value is $-0.0671$, compared with the exact value $0.0048$. This residual illustrates the distinction between energy stability and nodal monotonicity.

\FloatBarrier

The generalized eigenproblems are solved by the Krylov--Schur method in SLEPc with a relative eigensolver tolerance of $10^{-8}$ and MUMPS for the subsidiary linear systems. These runs use Firedrake~2026.4.1 and PETSc/SLEPc~3.25.0, with eight MPI processes and one thread per process. Table~\ref{tab:lowperm-p2-costs} reports wall-clock times and outer eigenvalue iteration counts. The LPS Krylov subspace dimensions are $40$ and $80$ for the two mesh sizes shown. In the present implementation, LPS parameter computation at $h=1/128$ costs about $157$ times the forward calculation, whereas BP parameter computation is less costly than the forward calculation. The parameters can be reused for additional right-hand sides when the discrete operator remains fixed.

\begin{table}[!htbp]
\centering
\caption{Parameter computation and forward calculation for $\kappa=10^{-8}$ on eight MPI processes. Times are wall-clock seconds. Parameter time includes assembly, solver setup, and eigenvalue computation. Forward time includes assembly, solution, and error evaluation; profile sampling and plotting are excluded.}
\label{tab:lowperm-p2-costs}
\begin{tabular}{ccrrr}
\toprule
$1/h$ & Method & Eigenvalue iterations & Parameter time (s) & Forward time (s) \\
\midrule
32  & BP  & 2   & $0.51$    & $0.89$ \\
32  & LPS & 15  & $3.72$    & $0.88$ \\
128 & BP  & 1   & $8.90$    & $13.17$ \\
128 & LPS & 238 & $2166.93$ & $13.81$ \\
\bottomrule
\end{tabular}
\end{table}


\subsection{Example 3: Layered heterogeneous medium}
\label{subsec:ex2-layered}
 
This example investigates the performance of the proposed stabilized formulation for heterogeneous media with discontinuous and low permeability coefficients. The objective is to assess the robustness of the pressure stabilization in the presence of sharp material interfaces, where standard equal-order discretizations are prone to spurious pressure oscillations.
 
We consider the square domain $\Omega=(0,1)^2$,
which is partitioned into three horizontal layers,
\[
\Omega_1=(0,1)\times(0,1/3),\qquad
\Omega_2=(0,1)\times(1/3,2/3),\qquad
\Omega_3=(0,1)\times(2/3,1).
\]
The permeability is assumed to be piecewise constant,
\[
\kappa=
\begin{cases}
10^{-4}, & \text{in }\Omega_1,\\
10^{-5}, & \text{in }\Omega_2,\\
10^{-6}, & \text{in }\Omega_3,
\end{cases}
\]
representing a three-layer medium with progressively decreasing permeability. The remaining physical parameters are listed in Table~\ref{tab:layered-parameters}. Homogeneous Dirichlet boundary conditions are prescribed for the displacement on the top boundary. On the bottom boundary, we impose
\[
p=10^{-2},
\qquad
(\boldsymbol{\sigma}-\alpha p\mathbb I)\bm n=(0,-10^{-3}),
\]
whereas the remaining boundaries are subject to natural boundary conditions.

\begin{table}[!ht]
\centering
\caption{Key parameter values used in Example 3.}
\label{tab:layered-parameters}
\begingroup
\setlength{\tabcolsep}{8pt} 
\begin{tabular}{lllllllll}
\toprule
Parameter & $E\,[\mathrm{dyn}\,\mathrm{cm}^{-2}]$ & $\nu$ & $\mu_f\,[\mathrm{Poise}]$ & $\kappa\,[\mathrm{cm}^2]$ & $\omega\,[\mathrm{rad}\, \mathrm{s}^{-1}]$ & $\rho\,[\mathrm{g}\,\mathrm{cm}^{-3}]$ & $\alpha$ & $B$ \\
\midrule
Value & $10^2$ & $0.45$ & $10^{-2}$ & $10^{-4}\,|\,10^{-5}\,|\,10^{-6}$ & $1\,|\,10\,|\,25\,|\,50$ & $1.0$ & $0.5$ & $0.8$ \\
\bottomrule
\end{tabular}
\endgroup
\end{table}

\begin{figure}[!ht]
    \centering
    \includegraphics[width=0.45\textwidth]{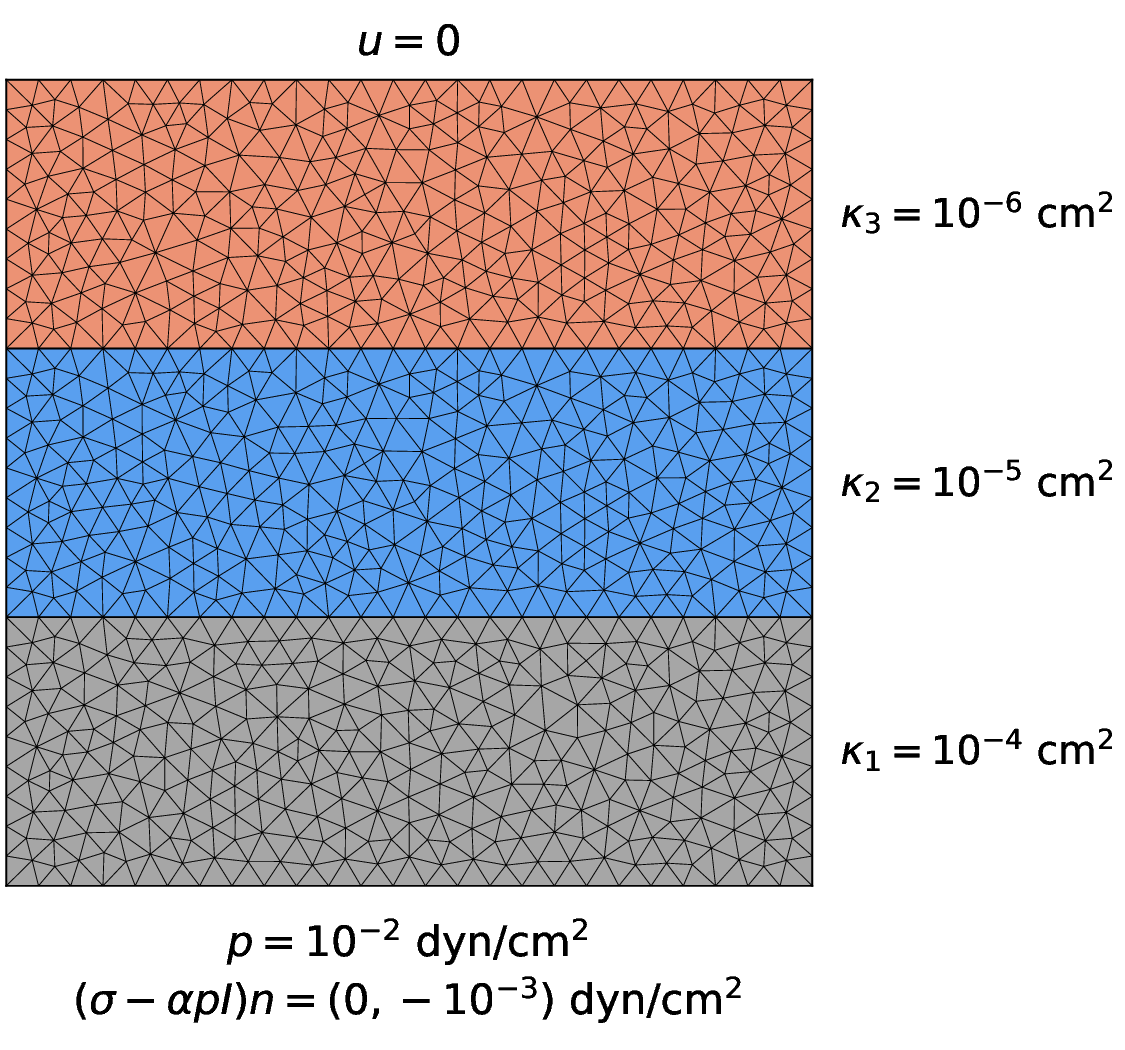}
    \caption{Example~3: Layered heterogeneous medium with discontinuous permeability coefficients (coarse mesh).}
    \label{fig:ex2-setup}
\end{figure}
 
The computational domain is discretized by an interface-fitted unstructured triangular mesh generated with Gmsh. The interfaces $y=1/3$ and $y=2/3$ coincide with the edges of the element so that each element belongs entirely to a single material layer; see \cref{fig:ex2-setup}. The characteristic mesh size is approximately $h=4\times10^{-3}$. The computations use linear equal-order finite elements with BP stabilization.
For each excitation frequency, the two
stabilization parameters are evaluated according to the selection rules in Section~\ref{sec:parsel-parameter-selection}, with the heterogeneous permeability field retained in the calculation of $\delta_2$. Their frequency
dependence is negligible at four decimal places, and the four cases can thus be summarized by
\[
\delta_1\approx0.0145,
\qquad
\delta_2\approx0.9665.
\]
The corresponding unrounded, frequency-specific values are used in the simulations. The implementation is carried out in Firedrake with PETSc, and the resulting linear systems are solved using the sparse direct solver MUMPS.

\begin{figure}[!ht]
    \centering
    \includegraphics[width=0.90\textwidth]{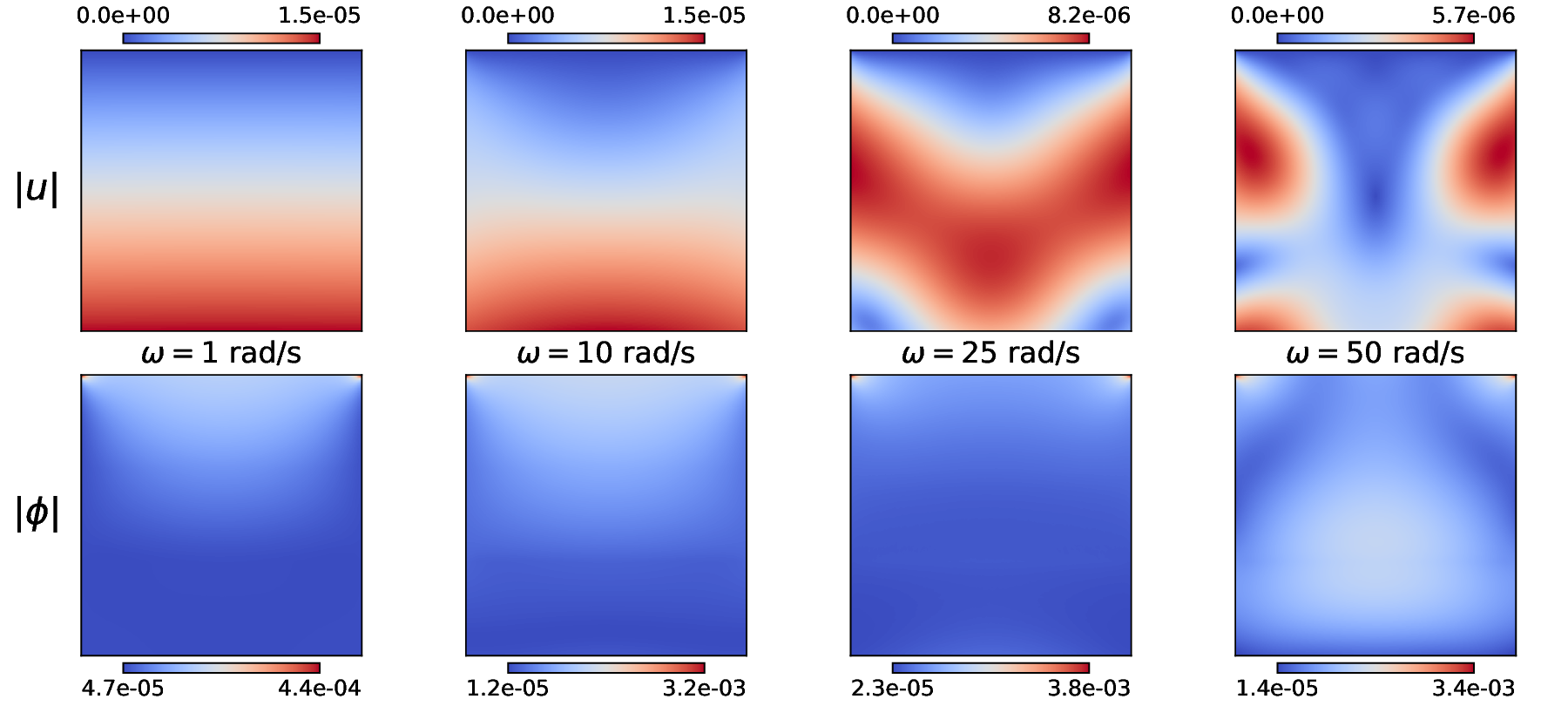}
    \caption{Magnitudes of the displacement \(\bm u\) and total pressure \(\phi\) for different excitation frequencies in the layered domain test.}
    \label{fig:ex2-uphi}
\end{figure}

\begin{figure}[!ht]
    \centering
    \includegraphics[width=0.90\textwidth]{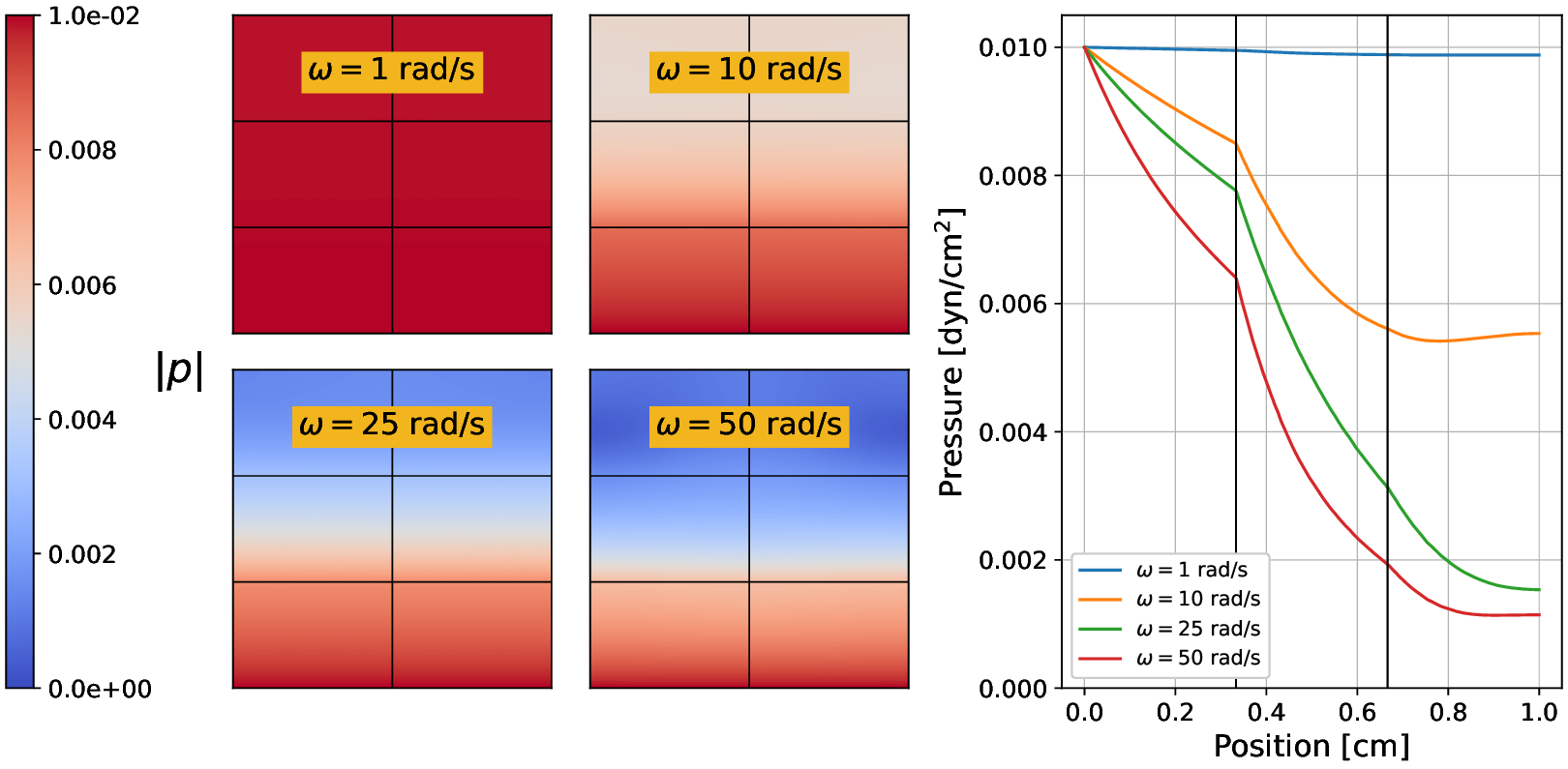}
    \caption{Pressure magnitude \(|p|\) and vertical pressure profiles along the line \(x=0.5\). The black lines indicate the permeability interfaces.}
    \label{fig:ex2-pressure}
\end{figure}

\cref{fig:ex2-uphi} displays the magnitudes of the displacement and total pressure for $\omega=1$, $10$, $25$, and $50~\mathrm{rad\,s^{-1}}$, while \cref{fig:ex2-pressure} presents the corresponding pore-pressure magnitudes and vertical profiles along $x=0.5$. As the excitation frequency increases, the displacement field exhibits increasingly pronounced oscillatory behavior, whereas the total pressure remains smooth across the permeability interfaces. The pore-pressure response also depends strongly on the excitation frequency. In particular, the pressure profiles remain continuous across the interfaces at $y=1/3$ and $y=2/3$, with changes in slope reflecting the discontinuous permeability distribution. Pressure continuity is enforced by the conforming finite element space. Despite the two orders of magnitude variation in permeability, no checkerboard patterns or spurious oscillations are observed near the material interfaces. These results demonstrate that the proposed pressure stabilization provides stable numerical approximations for heterogeneous low-permeability media with strongly discontinuous material coefficients.

\subsection{Example 4: Three-dimensional brain simulation}
This example assesses the stabilized equal-order discretization for the time-harmonic Biot model on an anatomically realistic three-dimensional brain geometry in the low-permeability regime. The calculation uses linear elements for displacement, pore pressure, and total pressure, together with the Brezzi--Pitk\"aranta (BP) stabilization.

Our mesh data are derived from \cite{mardal2022mathematical}, which converts segmented magnetic resonance images into anatomically meaningful volume meshes, including an exterior brain surface and an internal ventricle surface. Let $\Omega \subset \mathbb{R}^3$ denote the computational domain. The boundary decomposition and boundary conditions are based on the in silico brain setting of MRE \cite{galarce2023displacement}.  We decompose the boundary $\partial\Omega$ into disjoint sets:
\begin{equation*}
  \partial\Omega =
  \Gamma_{\mathrm{neck}}\cup
  \Gamma_{\mathrm{MRE}}\cup
  \Gamma_{\mathrm{ventricles}},
\end{equation*}
where \(\Gamma_{\mathrm{neck}}\) represents the inferior support region,
\(\Gamma_{\mathrm{MRE}}\) is the external surface on which the harmonic excitation is applied,
and \(\Gamma_{\mathrm{ventricles}}\) denotes the ventricular wall. The boundary conditions in the frequency domain are imposed as
\begin{alignat*}{2}
  \bm{u} &= \mathbf{0} 
      &\qquad &\text{on } \Gamma_{\mathrm{neck}}, \\
  \bm{\sigma}\bm{n}
      &= \bm{g}_{\mathrm{MRE}}
      &\qquad &\text{on } \Gamma_{\mathrm{MRE}}, \\
  \bm{\sigma}\boldsymbol{n}
      &= \mathbf{0}
      &\qquad &\text{on } \Gamma_{\mathrm{ventricles}}, \\
   p &= 0 
      &\quad &\text{on } \Gamma_{\mathrm{ventricles}}\cup\Gamma_{\mathrm{MRE}}, \\
  \partial_{\bm{n}} p &= 0 
      &\quad &\text{on } \Gamma_{\mathrm{neck}}. 
\end{alignat*}
Since $p=0$ on $\Gamma_{\mathrm{ventricles}}\cup\Gamma_{\mathrm{MRE}}$, the traction $\bm{\sigma}\bm{n}$ coincides with $(\bm{\sigma}-\alpha p\mathbb{I})\bm{n}$ in~\eqref{eq:displacement_bc}.
The pressure values prescribed on \(\Gamma_{\mathrm{ventricles}}\) and \(\Gamma_{\mathrm{MRE}}\) are
interpreted as harmonic perturbations around constant background  cerebrospinal fluid pressures; hence the
frequency-domain perturbation is set to zero on these boundaries.  The ventricular wall is
mechanically free, while the neck boundary is mechanically clamped and hydraulically impermeable.
Following the spatially distributed MRE forcing model of \cite{galarce2023displacement}, the
traction is chosen as
\begin{equation}
  \bm{g}_{\mathrm{MRE}}(\boldsymbol{x})
  = [0,A\,s(y),0],
  \qquad
  s(y)
  =
  1-\frac{y-y_{\min}}{y_{\max}-y_{\min}}.
  \label{brain_mre_load}
\end{equation}
Thus, the forcing amplitude is maximum on the posterior side and decays linearly along the frontal \(y\)-direction. The real-valued amplitude $A$ in \eqref{brain_mre_load} is chosen to obtain physiological magnitudes for the displacement fields.

\begin{table}[!ht]
\centering
\caption{Key parameter values used in the brain experiments.}
\label{tab:brain_parameters}
\begingroup
\setlength{\tabcolsep}{6pt} 
\begin{tabular}{llllllllll}
\toprule
Parameter & $E\,[\mathrm{dyn}\,\mathrm{cm}^{-2}]$ & $\nu$ & $\mu_f\,[\mathrm{Poise}]$ & $\kappa\,[\mathrm{cm}^2]$ & $\omega\,[\mathrm{rad}\, \mathrm{s}^{-1}]$ & $\rho\,[\mathrm{g}\,\mathrm{cm}^{-3}]$ & $\alpha$ & $B$ & $A\,[\mathrm{dyn}\,\mathrm{cm}^{-2}]$ \\
\midrule
Value & $10^4$ & $0.4$ & $0.01$ & $10^{-8}$ & $10.0$ & $1.0$ & $1.0$ & $0.9$ & $2.3\times10^{-3}$ \\
\bottomrule
\end{tabular}
\endgroup
\end{table}

\begin{figure}[!ht]
    \centering

    \captionsetup[subfigure]{
    format=hang,
    justification=raggedright,
    singlelinecheck=false,
    width=0.75\linewidth
    }

    \begin{subfigure}[t]{0.32\textwidth}
        \centering
        \includegraphics[
            width=\linewidth,  
            trim=300bp 10bp 300bp 10bp,
            clip
        ]{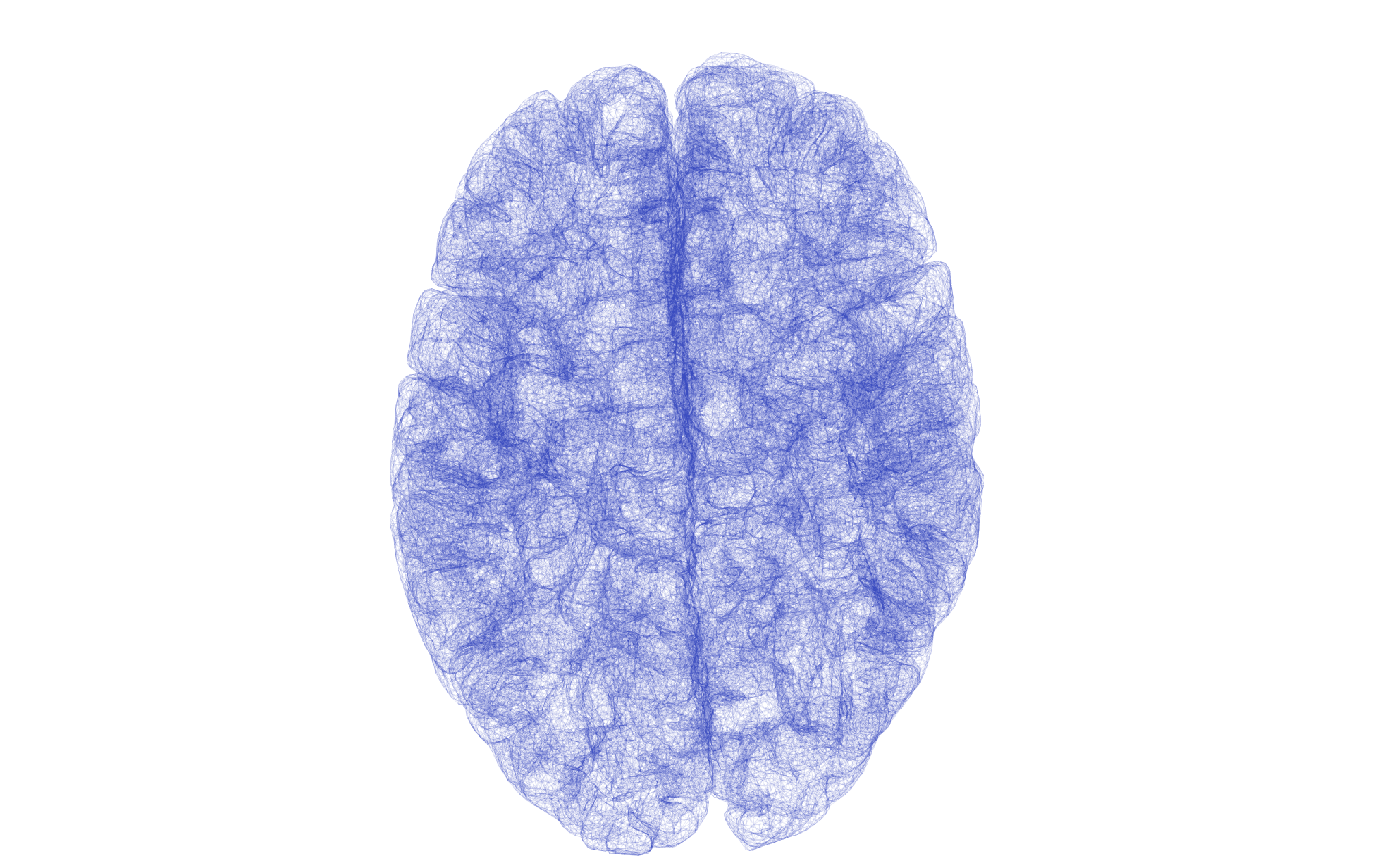}
        \caption{Brain mesh}
        \label{fig:mesh}
    \end{subfigure}%
    \begin{subfigure}[t]{0.339\textwidth}
        \centering
        \includegraphics[
            width=\linewidth,
            trim=400bp 120bp 200bp 10bp,
            clip
        ]{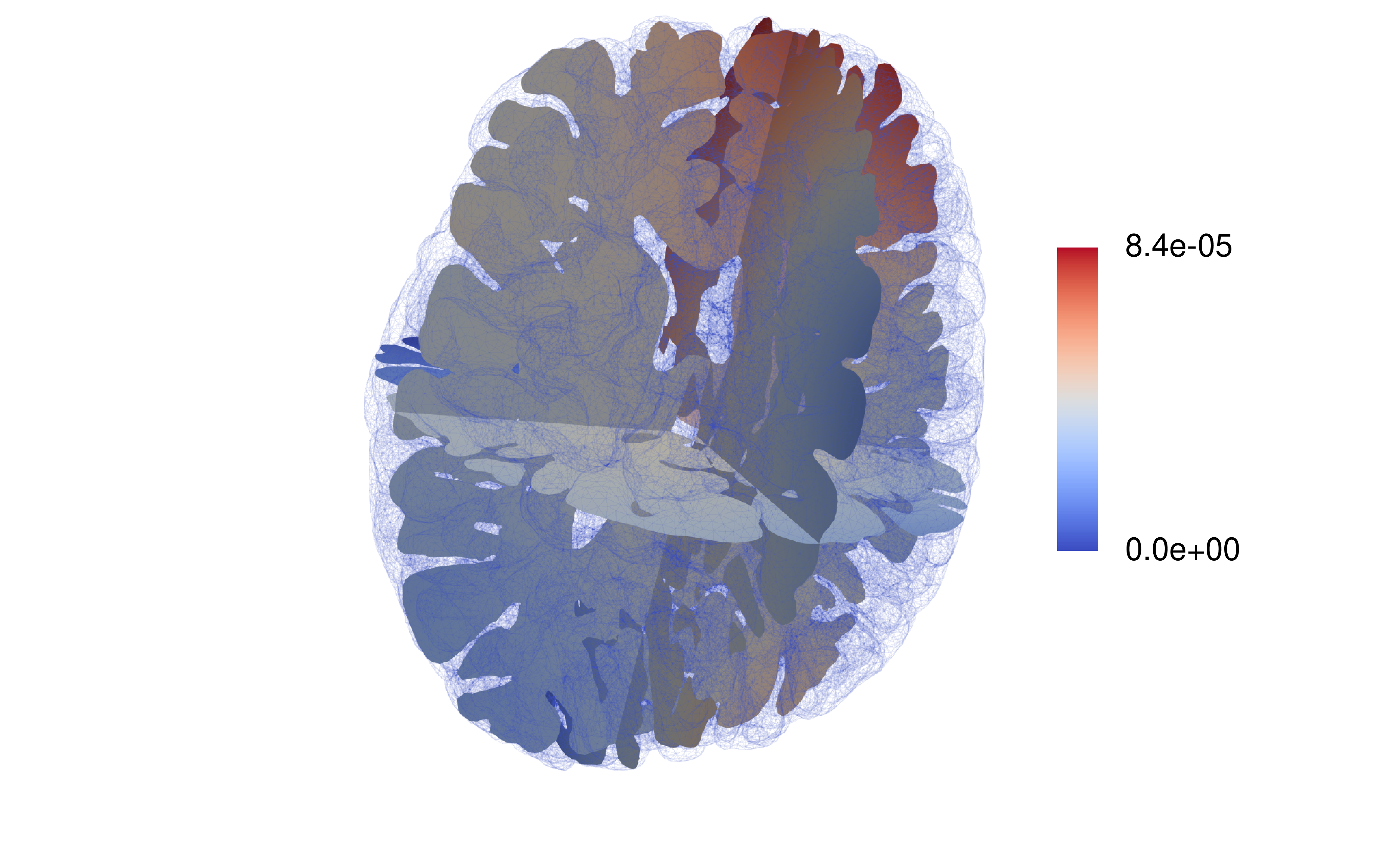}
        \caption{Displacement magnitude\protect\\slice view}
        \label{fig:u_slice}
    \end{subfigure}%
    \begin{subfigure}[t]{0.339\textwidth}
        \centering
        \includegraphics[
            width=\linewidth,
            trim=400bp 120bp 200bp 10bp,
            clip
        ]{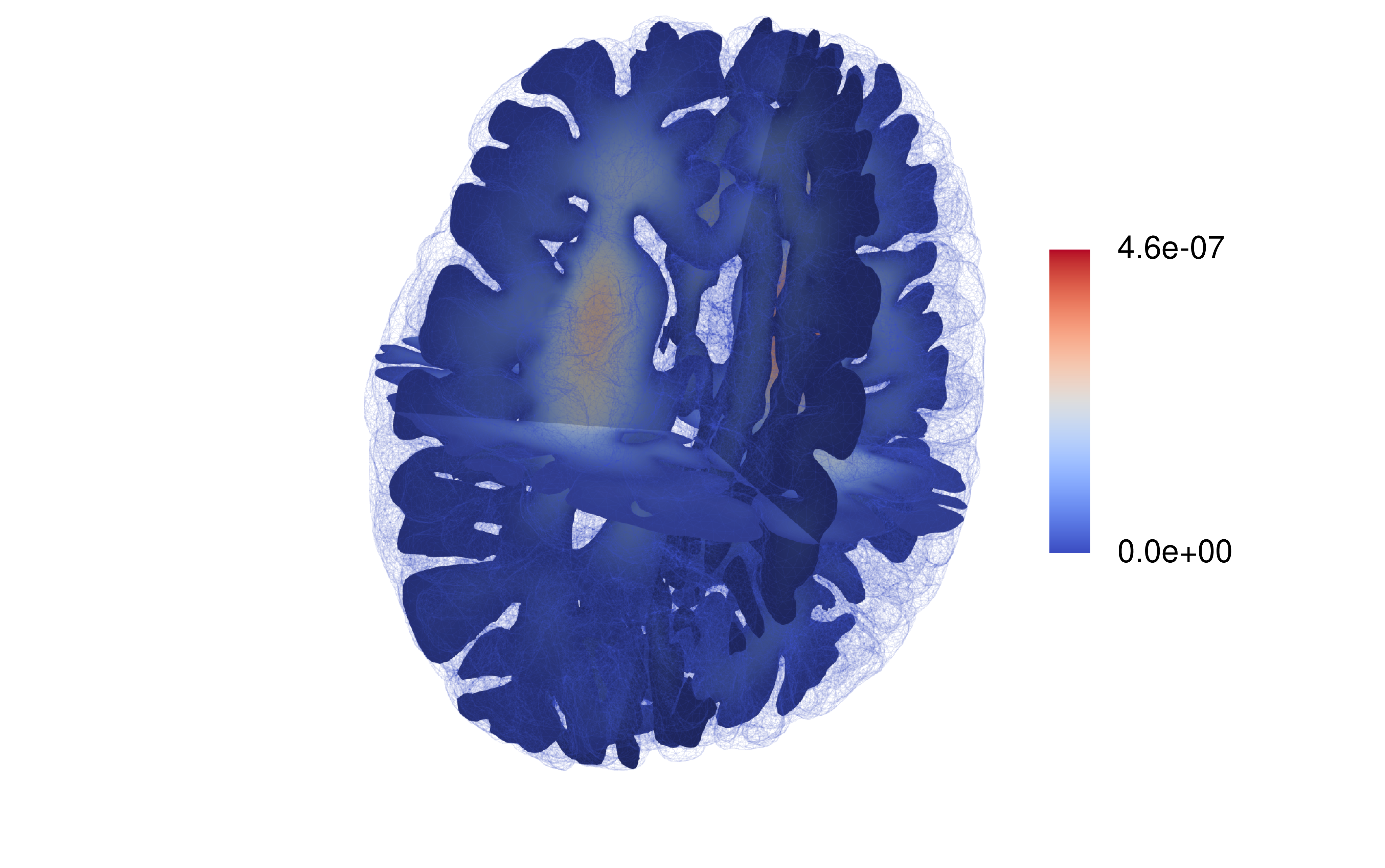}
        \caption{Pore-pressure magnitude
        \protect\\slice view}
        \label{fig:p_slice}
    \end{subfigure}%

    \caption{Example 4: Brain mesh and slice views of the displacement and pore-pressure magnitudes.}
    \label{fig:slice}
\end{figure}

The simulation is performed on a dual-socket AMD EPYC 9654 server using Firedrake, PETSc, and MPI, with 96 MPI ranks and $\texttt{OMP\_NUM\_THREADS}=1$. The mesh used for the calculations contains 178479 vertices and 917015 tetrahedra (see \cref{fig:mesh}). After the real and imaginary parts are separated, the linear system has 1784790 degrees of freedom. The coupled system is solved with flexible GMRES and a multiplicative two-field field-split preconditioner; the mechanics--total-pressure block and the pore-pressure block are both solved with MUMPS. The relative and absolute residual tolerances are $10^{-8}$ and $10^{-13}$, respectively, with a maximum of 1500 outer iterations.

All quantities are expressed in the CGS unit system. The physical parameters are listed in \cref{tab:brain_parameters}. The representative mesh size is $h=0.3123\,\mathrm{cm}$.
  Applying the parameter-selection rules in Section~\ref{sec:parsel-parameter-selection} gives
\[
\delta_1=1.3981\times10^{-4},
\qquad
\delta_2=1.0000.
\]
The SLEPc eigenvalue calculation uses Krylov--Schur with a relative tolerance of $10^{-6}$ and a Krylov subspace dimension of $40$. It converged in six outer iterations and required $143$ applications each of the Schur operator and its adjoint, including the residual check; the complementary-block solves reused the LU factors. Parameter selection took $119.07$~s, compared with $73.71$~s and one outer FGMRES iteration for the forward calculation, a time ratio of $1.62$. These timings are maxima over all MPI ranks after synchronized starts and include matrix assembly and solver setup. The parameter time also includes the eigenvalue calculation, verification, diagnostic output, and cleanup. Mesh loading and field post-processing are excluded. 

\begin{figure}[!ht]
    \centering

    \captionsetup[subfigure]{
    format=hang,
    justification=raggedright,
    singlelinecheck=false
    }
    
    \begin{subfigure}[t]{0.45\textwidth}
        \centering
        \includegraphics[
            width=\linewidth,
            trim=300bp 200bp 1bp 200bp, 
            clip
        ]{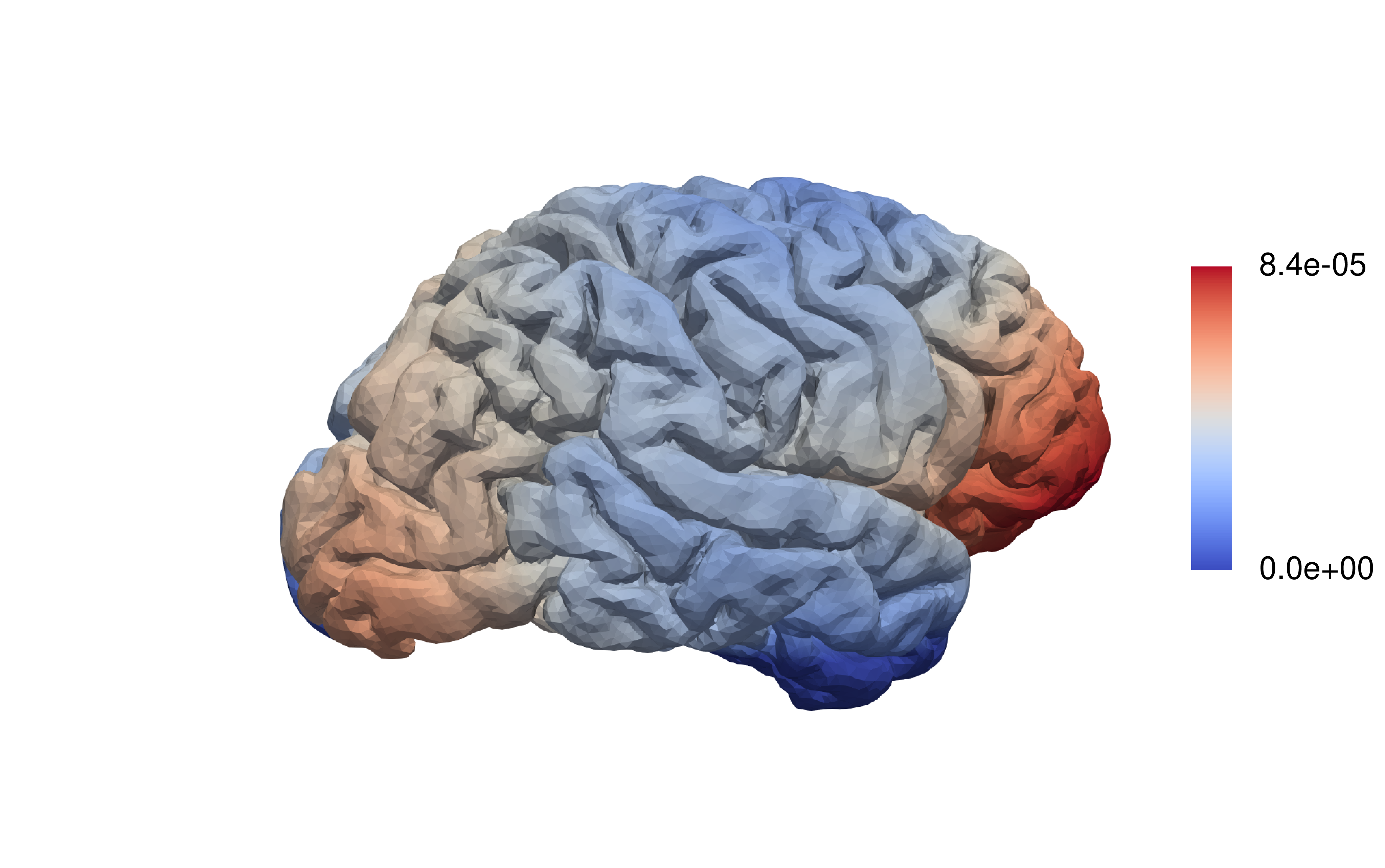}
        \caption{Displacement magnitude}
        \label{fig:u}
    \end{subfigure}%
    \begin{subfigure}[t]{0.45\textwidth}
        \centering
        \includegraphics[
            width=\linewidth,
            trim=300bp 200bp 1bp 200bp,
            clip
        ]{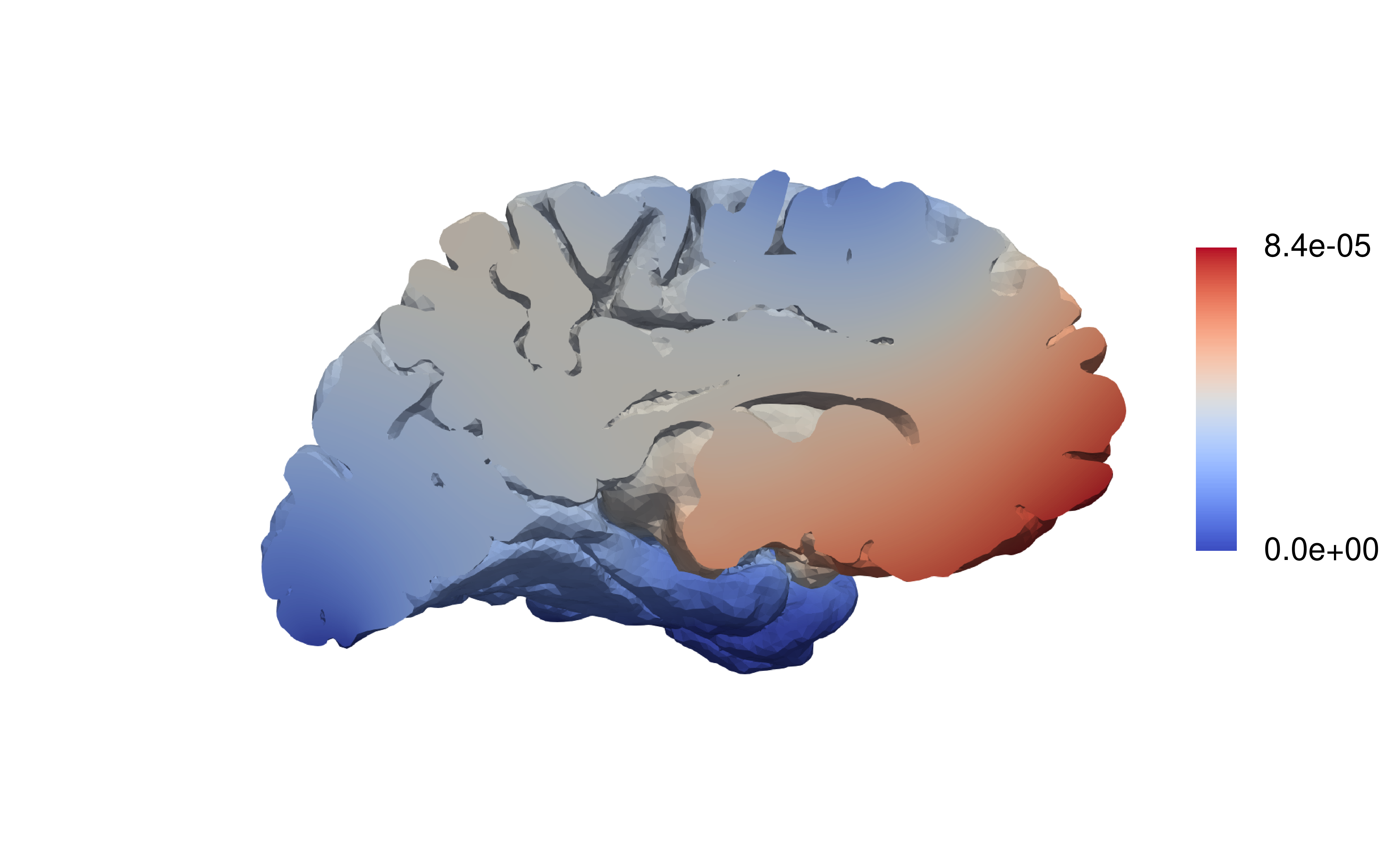}
        \caption{Displacement magnitude, clipped left half}
        \label{fig:u_cilp}
    \end{subfigure}%
    \\
    \begin{subfigure}[t]{0.43\textwidth}
        \centering
        \includegraphics[
            width=\linewidth,
            trim=300bp 200bp 1bp 100bp,
            clip
        ]{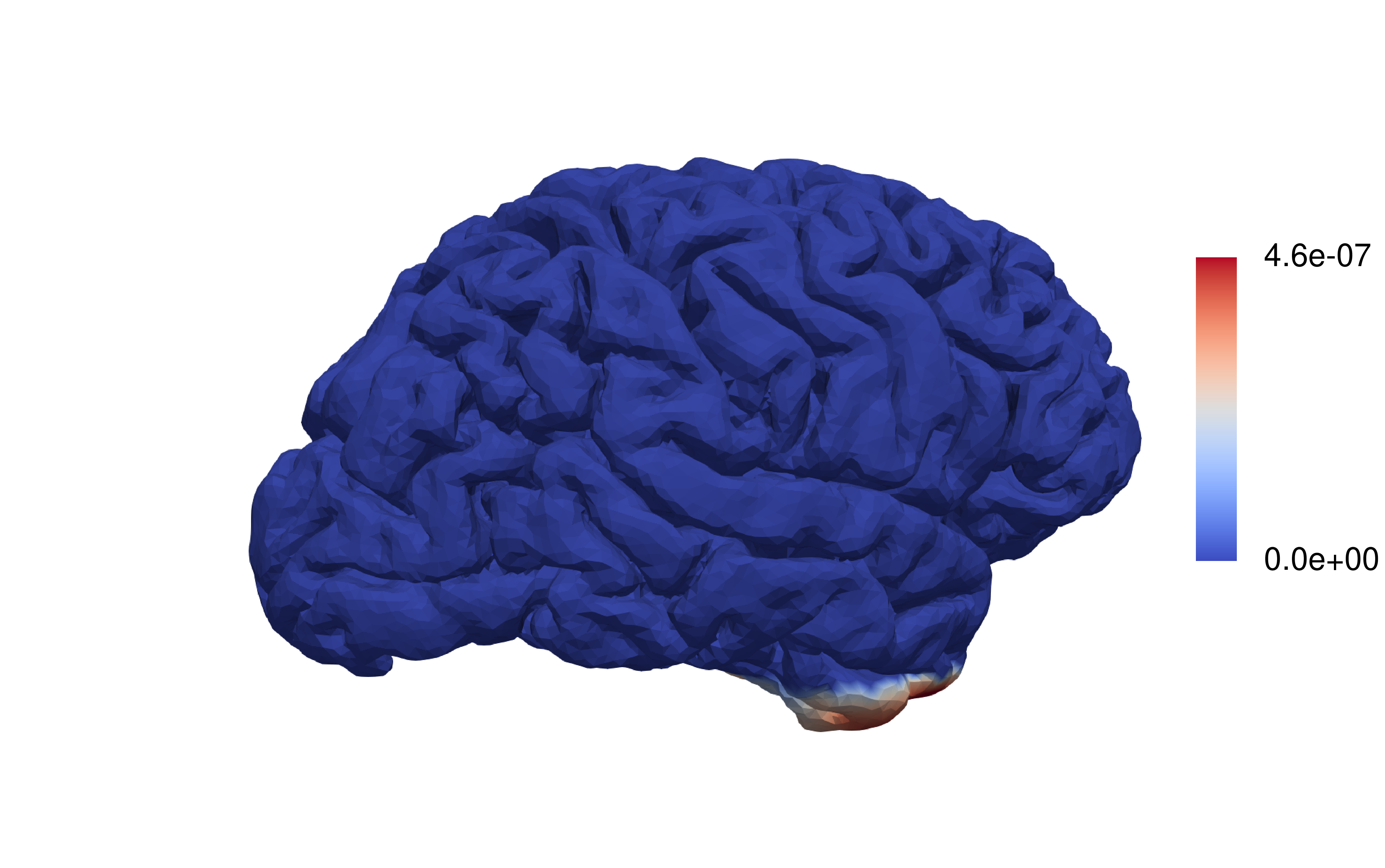}
        \caption{Pore-pressure magnitude}
        \label{fig:p}
    \end{subfigure}%
    \begin{subfigure}[t]{0.45\textwidth}
        \centering
        \includegraphics[
            width=\linewidth,
            trim=300bp 200bp 1bp 100bp,
            clip
        ]{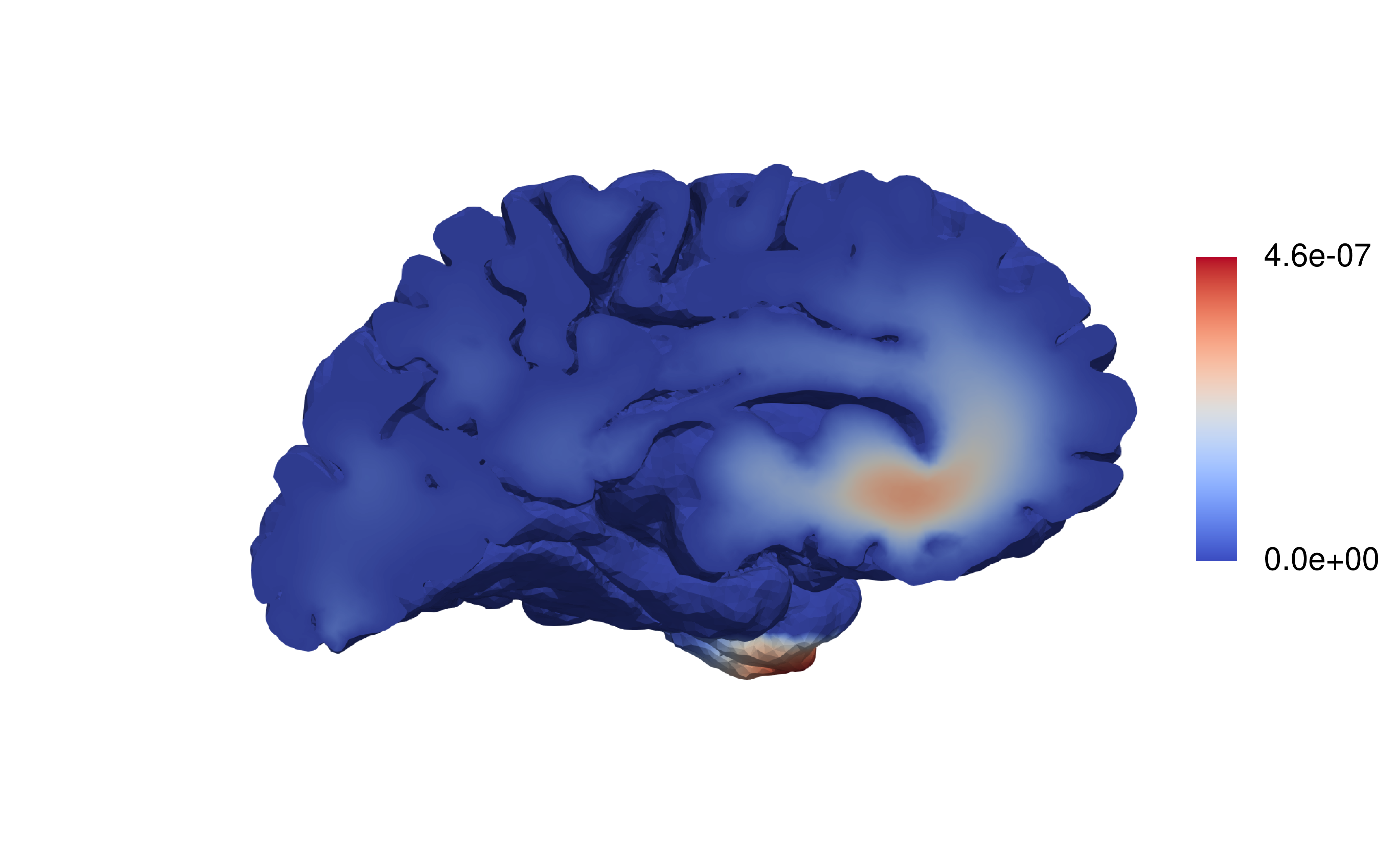}
        \caption{Pore-pressure magnitude, clipped left half}
        \label{fig:p_clip}
    \end{subfigure}%

    \caption{Example 4: Displacement and pore-pressure magnitudes, mesh clipped to the left half.}
    \label{fig:clip}
\end{figure}

The resulting displacement and pore-pressure magnitudes are shown in \cref{fig:clip}. The clipped views provide a clear visualization of the interior solution, revealing a smooth displacement field together with a spatially varying pore-pressure distribution that reflects the spatially varying pressure response induced by the geometry and applied loading. In addition,

\section{Conclusions and outlook}

In this work, we have developed and analyzed a stabilized equal-order finite element method for a three-field formulation of the time-harmonic Biot model in terms of displacement, pore pressure, and total pressure. At the continuous level, the well-posedness of the complex-valued mixed formulation was established under a non-resonance condition through a Fredholm framework based on the decomposition of the variational operator into an invertible principal part and a compact perturbation. At the discrete level, an elementwise momentum-residual stabilization was combined with a pressure local projection stabilization to overcome the lack of discrete inf--sup stability while avoiding excessive pressure-gradient penalization. The resulting finite element approximation was shown to be well-posed on sufficiently fine meshes, and optimal-order a priori error estimates were derived in a natural energy norm.

The proposed method was validated through a series of numerical experiments in both two and three spatial dimensions. Benchmark problems with exact solutions confirmed the predicted convergence rates for linear and quadratic equal-order discretizations. Simulations in layered heterogeneous media demonstrated the robustness of the proposed stabilization in the presence of strongly discontinuous and low-permeability coefficients, without exhibiting visible spurious pressure oscillations. Finally, simulations on anatomically realistic brain geometries illustrated the applicability of the proposed formulation to large-scale magnetic resonance elastography, producing smooth displacement and pore-pressure distributions on complex three-dimensional meshes.

Several directions merit further investigation. From a theoretical perspective, it is of interest to establish parameter-robust stability and error estimates with explicit dependence on the physical parameters and excitation frequency, and to develop adaptive strategies for selecting the stabilization parameters. On the computational side, scalable iterative solvers and efficient block preconditioners are needed for large-scale three-dimensional simulations. Future work will also extend the proposed forward solver to more realistic heterogeneous tissue models and integrate it into inverse problem frameworks for quantitative magnetic resonance elastography and other biomedical imaging applications.

\section*{CRediT authorship contribution statement}
\textbf{Mingchao Cai:} Conceptualization, Methodology, Formal analysis, Writing - Review \& Editing, Funding acquisition. \textbf{Yiming Fu:} Software, Validation, Formal analysis, Investigation, Writing - Review \& Editing. \textbf{Jingzhi Li:} Conceptualization, Methodology, Writing - Review \& Editing, Funding acquisition. 

\section*{Declaration of Generative AI and AI-assisted technologies in the writing process}
During the preparation of this work, the authors used ChatGPT to polish the writing. After using this tool/service, the authors reviewed and edited the content as needed and take full responsibility for the content of the published article.

\section*{Declaration of competing interest}
The authors declare that they have no known competing financial interests or personal relationships that could have appeared 
to influence the work reported in this paper.

\section*{Acknowledgements}
M. Cai is supported in part by the U. S. National Science Foundation under Grant No. DMS-2602787. Y. Fu and J. Li are partially supported by the Shenzhen Sci-Tech Fund No. RCJC20200714114556020, Guangdong Basic and Applied Research Fund No. 2023B1515250005,  National Center for Applied Mathematics Shenzhen, and SUSTech International Center for Mathematics.  

Any opinions, findings, conclusions, or recommendations expressed in this material are those of the author(s) and do not necessarily reflect the views of the National Science Foundation.


\appendix

\section{Algebraic selection of the pressure stabilization parameter}
\label{app:parsel-delta2-selection}

This appendix derives~\eqref{eq:parsel-delta2-main} for prescribed physical
parameters, mesh, polynomial degree, frequency, boundary conditions, and
$\delta_1$.  The derivation consists of three steps: decomposition of the
discrete pressure space by the local projection, elimination of the remaining
variables through a Schur complement, and imposition of a Hermitian stability
criterion on the fluctuation pressure space.

\subsection{Resolved and fluctuation pressure spaces}
\label{app:parsel-pressure-decomposition}

Define the discrete LPS operator
$\mathcal L_h:P_h\rightarrow P_h^*$ by
\begin{equation}
 \langle\mathcal L_hp^h,q^h\rangle
 :=
 \ell_h(p^h,q^h)
 =
 \sum_{T\in\mathcal T_h}
 \frac{h_T^2}{\mu_f\omega\alpha}
 (\mathcal K_h\nabla p^h,
  \mathcal K_h\nabla q^h)_T.
 \label{eq:parsel-Lh}
\end{equation}
We define the resolved pressure space and the fluctuation pressure space by
\begin{align}
 P_h^r
 &:={}
 \ker\mathcal L_h
 =\left\{q^h\in P_h:\ell_h(q^h,q^h)=0\right\},
 \label{eq:parsel-resolved-space}
 \\
 P_h^f
 &:={}
 \left(P_h^r\right)^{\perp_{L^2(\Omega)}}.
 \label{eq:parsel-fluctuation-space}
\end{align}
The superscripts $r$ and $f$ stand for resolved and fluctuation,
respectively. The latter is also referred to as the unresolved pressure
space.  By~\eqref{eq:parsel-Lh},
\begin{equation}
 q^h\in P_h^r
 \quad\Longleftrightarrow\quad
 \mathcal K_h\nabla q^h=0.
 \label{eq:parsel-resolved-characterization}
\end{equation}

Let $N_p:=\dim P_h$. Since the sesquilinear form
$\ell_h(\cdot,\cdot)$ in~\eqref{eq:parsel-Lh}
is Hermitian and positive semidefinite, the finite-dimensional spectral theorem gives an $L^2(\Omega)$ orthonormal basis $\{z_j^h\}_{j=1}^{N_p}$ of $P_h$ and
eigenvalues $0\leq\zeta_1\leq\cdots\leq\zeta_{N_p}$ such that
\begin{equation}
 \ell_h(z_j^h,q^h)
 =
 \zeta_j(z_j^h,q^h)_\Omega
 \quad\forall q^h\in P_h,
 \qquad
 (z_i^h,z_j^h)_\Omega=\delta_{ij}.
 \label{eq:parsel-projection-eigenproblem}
\end{equation}
Taking $q^h=z_j^h$ in~\eqref{eq:parsel-projection-eigenproblem} gives
\begin{equation*}
 \ell_h(z_j^h,z_j^h)
 =
 \zeta_j\|z_j^h\|_{0,\Omega}^2.
 \label{eq:parsel-eigenvalue-energy-identity}
\end{equation*}
Hence, $\zeta_j=0$ if and only if $z_j^h\in\ker\mathcal L_h$.  With
\begin{equation*}
 I_r:=\{j:\zeta_j=0\},
 \qquad
 I_f:=\{j:\zeta_j>0\},
 \label{eq:parsel-spectral-index-sets}
\end{equation*}
the spaces in~\eqref{eq:parsel-resolved-space} and
\eqref{eq:parsel-fluctuation-space} satisfy
\begin{equation*}
 P_h^r=\operatorname{span}\{z_j^h:j\in I_r\},
 \qquad
 P_h^f=\operatorname{span}\{z_j^h:j\in I_f\}.
 \label{eq:parsel-spectral-pressure-spaces}
\end{equation*}
Consequently, every $p^h\in P_h$ has the unique orthogonal decomposition
\begin{align}
 p^h&=p_r^h+p_f^h,
 \label{eq:parsel-pressure-splitting}\\
 p_r^h
 &:={}
 \sum_{j\in I_r}(p^h,z_j^h)_\Omega z_j^h,
 \qquad
 p_f^h
 :=
 \sum_{j\in I_f}(p^h,z_j^h)_\Omega z_j^h.
 \nonumber
\end{align}
In finite-precision arithmetic, the numerical set $I_r$ may be identified by
$\zeta_j\leq\mathrm{tol}_{\ker}\zeta_{\max}$, where
$\zeta_{\max}:=\max_j\zeta_j$ and the relative tolerance is reported with the computational parameters.

Let $\mathcal L_f:P_h^f\rightarrow(P_h^f)^*$ be the restriction of
$\mathcal L_h$ to $P_h^f$:
\begin{equation*}
 \langle\mathcal L_fp_f^h,q_f^h\rangle
 :=
 \ell_h(p_f^h,q_f^h)
 \qquad
 \forall p_f^h,q_f^h\in P_h^f.
 \label{eq:parsel-Lf}
\end{equation*}
The operator $\mathcal L_f$ is positive definite. Indeed, a function in $P_h^f$ with zero LPS energy also belongs to
$\ker\mathcal L_h=P_h^r$, and must therefore vanish.  Moreover, \eqref{eq:parsel-resolved-characterization} and the Hermitian property imply
\begin{equation}
 \ell_h(p^h,q^h)=\ell_h(p_f^h,q_f^h)
 \qquad
 \forall p^h=p_r^h+p_f^h,
 \quad
 q^h=q_r^h+q_f^h.
 \label{eq:parsel-LPS-fluctuation-only}
\end{equation}
Equation~\eqref{eq:parsel-LPS-fluctuation-only} explicitly shows why the second stabilization acts only on the fluctuation pressure block. For $k=1$, the resolved space contains only pressure functions with zero gradient,
subject to the imposed pressure boundary condition or gauge constraint. For $k\geq2$, it also contains pressure functions whose gradients belong to $D_h$.

\subsection{Complete block operator and fluctuation Schur complement}
\label{app:parsel-schur-complement}

All discrete operators below incorporate the same essential boundary conditions and pressure constraint as the forward problem. Introduce the displacement part of the strong momentum residual,
\begin{equation*}
 \mathbf R_u(\bm v^h)
 :=
 \omega^2\rho\bm v^h
 +2\mu_e\operatorname{div}\varepsilon(\bm v^h).
 \label{eq:parsel-Ru}
\end{equation*}
For $\bm u^h,\bm v^h\in\mathbf H_h$,
$p^h,q^h\in P_h$, and $\phi^h,\xi^h\in S_h$, define
\begin{align*}
 \langle\mathcal A_{uu}\bm u^h,\bm v^h\rangle
 &:={}
 2\mu_e(\varepsilon(\bm u^h),\varepsilon(\bm v^h))_\Omega
 -\omega^2\rho(\bm u^h,\bm v^h)_\Omega
 +\delta_1\sum_{T\in\mathcal T_h}h_T^2
 (\mathbf R_u(\bm u^h),\mathbf R_u(\bm v^h))_T,
 \\
 \langle\mathcal A_{u\phi}\phi^h,\bm v^h\rangle
 &:={}
 -(\phi^h,\operatorname{div}\bm v^h)_\Omega
 -\delta_1\sum_{T\in\mathcal T_h}h_T^2
 (\nabla\phi^h,\mathbf R_u(\bm v^h))_T,
 \\
 \langle\mathcal A_{\phi u}\bm u^h,\xi^h\rangle
 &:={}
 (\operatorname{div}\bm u^h,\xi^h)_\Omega
 -\delta_1\sum_{T\in\mathcal T_h}h_T^2
 (\mathbf R_u(\bm u^h),\nabla\xi^h)_T,
 \\
 \langle\mathcal A_{\phi\phi}\phi^h,\xi^h\rangle
 &:={}
 \lambda^{-1}(\phi^h,\xi^h)_\Omega
 +\delta_1\sum_{T\in\mathcal T_h}h_T^2
 (\nabla\phi^h,\nabla\xi^h)_T,
 \\
 \langle\mathcal A_{pp}(\kappa)p^h,q^h\rangle
 &:={}
 i\theta\lambda^{-1}(p^h,q^h)_\Omega
 +\frac{\kappa}{\mu_f\omega\alpha}
 (\nabla p^h,\nabla q^h)_\Omega,
 \\
 \langle\mathcal A_{p\phi}\phi^h,q^h\rangle
 &:={}
 -i\lambda^{-1}(\phi^h,q^h)_\Omega,
 \\
 \langle\mathcal A_{\phi p}p^h,\xi^h\rangle
 &:={}
 -\alpha\lambda^{-1}(p^h,\xi^h)_\Omega.
\end{align*}
These definitions contain the elasticity, inertia, storage, physical pressure diffusion, coupling, and residual stabilization contributions. The discrete system without the second stabilization term is
\begin{equation*}
 \begin{bmatrix}
  \mathcal A_{uu} & 0 & \mathcal A_{u\phi}\\
  0 & \mathcal A_{pp}(\kappa) & \mathcal A_{p\phi}\\
  \mathcal A_{\phi u} & \mathcal A_{\phi p} & \mathcal A_{\phi\phi}
 \end{bmatrix}
 \begin{bmatrix}
  \bm u^h\\ p^h\\ \phi^h
 \end{bmatrix}
 =
 \begin{bmatrix}
  \mathcal F_u\\ \mathcal F_p\\ \mathcal F_\phi
 \end{bmatrix}.
 \label{eq:parsel-complete-original-matrix}
\end{equation*}
The second stabilization adds $\delta_2\mathcal L_h$ only to the $p$ block.

To state the system after the decomposition~\eqref{eq:parsel-pressure-splitting}, let $a,b\in\{r,f\}$ and define the restricted pressure blocks by
\begin{equation*}
 \begin{aligned}
 \langle\mathcal A_{ab}(\kappa)p_b^h,q_a^h\rangle
 &:={}
 i\theta\lambda^{-1}(p_b^h,q_a^h)_\Omega
 +\frac{\kappa}{\mu_f\omega\alpha}
 (\nabla p_b^h,\nabla q_a^h)_\Omega,
 \\
 \langle\mathcal A_{\phi a}p_a^h,\xi^h\rangle
 &:={}
 -\alpha\lambda^{-1}(p_a^h,\xi^h)_\Omega,
 \\
 \langle\mathcal A_{a\phi}\phi^h,q_a^h\rangle
 &:={}
 -i\lambda^{-1}(\phi^h,q_a^h)_\Omega,
 \end{aligned}
\end{equation*}
for all $p_b^h\in P_h^b$, $q_a^h\in P_h^a$, and $\xi^h\in S_h$.  Let $\mathcal F_r$ and $\mathcal F_f$ denote the corresponding restrictions of $\mathcal F_p$.  By~\eqref{eq:parsel-LPS-fluctuation-only}, reordering the
unknowns as $\bm u^h,\phi^h,p_r^h,p_f^h$ gives
\begin{equation}
 \begin{bmatrix}
  \mathcal A_{uu} & \mathcal A_{u\phi} & 0 & 0\\
  \mathcal A_{\phi u} & \mathcal A_{\phi\phi} &
  \mathcal A_{\phi r} & \mathcal A_{\phi f}\\
  0 & \mathcal A_{r\phi} & \mathcal A_{rr}(\kappa) &
  \mathcal A_{rf}(\kappa)\\
  0 & \mathcal A_{f\phi} & \mathcal A_{fr}(\kappa) &
  \mathcal A_{ff}(\kappa)+\delta_2\mathcal L_f
 \end{bmatrix}
 \begin{bmatrix}
  \bm u^h\\ \phi^h\\ p_r^h\\ p_f^h
 \end{bmatrix}
 =
 \begin{bmatrix}
  \mathcal F_u\\ \mathcal F_\phi\\ \mathcal F_r\\ \mathcal F_f
 \end{bmatrix}.
 \label{eq:parsel-complete-transformed-matrix}
\end{equation}

For the elimination step, let the subscript $c$ denote the complementary variables $(\bm u^h,\phi^h,p_r^h)$ and define
\begin{equation*}
 \begin{aligned}
 \mathcal A_c(\kappa)
 &:={}
 \begin{bmatrix}
  \mathcal A_{uu} & \mathcal A_{u\phi} & 0\\
  \mathcal A_{\phi u} & \mathcal A_{\phi\phi} & \mathcal A_{\phi r}\\
  0 & \mathcal A_{r\phi} & \mathcal A_{rr}(\kappa)
 \end{bmatrix},
 \\
 \mathcal A_{cf}(\kappa)
 &:={}
 \begin{bmatrix}
  0\\ \mathcal A_{\phi f}\\ \mathcal A_{rf}(\kappa)
 \end{bmatrix},
 \qquad
 \mathcal A_{fc}(\kappa)
 :=
 \begin{bmatrix}
  0 & \mathcal A_{f\phi} & \mathcal A_{fr}(\kappa)
 \end{bmatrix}.
 \end{aligned}
\end{equation*}
Assume that $\mathcal A_c(\kappa)$ is invertible. This excludes a resonance or loss of solvability in the block containing the mechanics and resolved pressure; changing $\delta_2$ cannot repair such a loss of solvability. 
The first three rows of~\eqref{eq:parsel-complete-transformed-matrix} give
\begin{equation*}
 \begin{bmatrix}
  \bm u^h\\ \phi^h\\ p_r^h
 \end{bmatrix}
 =
 \mathcal A_c(\kappa)^{-1}
 \left(
 \begin{bmatrix}
  \mathcal F_u\\ \mathcal F_\phi\\ \mathcal F_r
 \end{bmatrix}
 -\mathcal A_{cf}(\kappa)p_f^h
 \right).
 \label{eq:parsel-eliminate-complement}
\end{equation*}
Substitution into the fourth row yields
\begin{equation}
 \left(\mathcal S_f^0(\kappa)+\delta_2\mathcal L_f\right)p_f^h
 =\mathcal R_f,
 \label{parsel-reduced-equation}
\end{equation}
where
\begin{align*}
 \mathcal S_f^0(\kappa)
 &:={}
 \mathcal A_{ff}(\kappa)
 -\mathcal A_{fc}(\kappa)
  \mathcal A_c(\kappa)^{-1}
  \mathcal A_{cf}(\kappa),
 \\
 \mathcal R_f
 &:={}
 \mathcal F_f
 -\mathcal A_{fc}(\kappa)\mathcal A_c(\kappa)^{-1}
 \begin{bmatrix}
  \mathcal F_u\\ \mathcal F_\phi\\ \mathcal F_r
 \end{bmatrix}.
\end{align*}
The superscript $0$ in $\mathcal S_f^0$ indicates that this Schur complement is formed before adding the second stabilization term, that is, with $\delta_2=0$. 

The stabilized fluctuation pressure operator is therefore
\begin{equation*}
 \mathcal S_f(\kappa,\delta_2)
 :=
 \mathcal S_f^0(\kappa)+\delta_2\mathcal L_f.
 \label{eq:parsel-Sf-delta}
\end{equation*}
This identity is exact: the physical diffusion in
$\mathcal A_{pp}(\kappa)$ and the LPS operator $\mathcal L_h$ remain distinct. The effects of $\mu_f$, $\omega$, $\alpha$, $\rho$, $\mu_e$, $\lambda$, storage, the boundary conditions, and $\delta_1$ are retained in $\mathcal S_f^0(\kappa)$.

\subsection{Hermitian fluctuation stability and selection of
\texorpdfstring{$\delta_2$}{delta2}}
\label{app:parsel-hermitian-criterion}

Since $\mathcal L_f$ is positive definite, it defines the norm
\begin{equation*}
 \|q_f^h\|_{\mathcal L_f}
 :=
 \ell_h(q_f^h,q_f^h)^{1/2}
 \qquad
 \forall q_f^h\in P_h^f.
 \label{eq:parsel-Lf-norm}
\end{equation*}
For $\mathcal G_f\in(P_h^f)^*$, define the corresponding dual norm by
\begin{equation*}
 \|\mathcal G_f\|_{\mathcal L_f^{-1}}
 :=
 \sup_{0\neq q_f^h\in P_h^f}
 \frac{|\langle\mathcal G_f,q_f^h\rangle|}
      {\|q_f^h\|_{\mathcal L_f}}.
 \label{eq:parsel-Lf-dual-norm}
\end{equation*}
Let $\beta_*>0$ be a prescribed dimensionless margin, chosen independently of $\delta_2$.  We require
\begin{equation}
 \operatorname{Re}
 \langle\mathcal S_f(\kappa,\delta_2)q_f^h,q_f^h\rangle
 \geq
 \beta_*\ell_h(q_f^h,q_f^h)
 \qquad
 \forall q_f^h\in P_h^f.
 \label{eq:parsel-hermitian-stability}
\end{equation}
This condition gives a load-independent stability estimate.  Indeed, duality and~\eqref{eq:parsel-hermitian-stability} imply
\begin{equation*}
 \beta_*\|q_f^h\|_{\mathcal L_f}^2
 \leq
 \operatorname{Re}
 \langle\mathcal S_f(\kappa,\delta_2)q_f^h,q_f^h\rangle
 \leq
 \|\mathcal S_f(\kappa,\delta_2)q_f^h\|_{\mathcal L_f^{-1}}
 \|q_f^h\|_{\mathcal L_f}.
 \label{eq:parsel-hermitian-inverse-bound}
\end{equation*}
Applying this estimate to~\eqref{parsel-reduced-equation} yields
\begin{equation}
 \|p_f^h\|_{\mathcal L_f}
 \leq
 \beta_*^{-1}\|\mathcal R_f\|_{\mathcal L_f^{-1}}.
 \label{eq:parsel-pressure-amplification-bound}
\end{equation}

Define the Hermitian part of the unstabilized Schur form by
\begin{equation*}
 s_{f,\mathrm H}^0(p_f^h,q_f^h)
 :=
 \frac{1}{2}
 \left(
  \langle\mathcal S_f^0(\kappa)p_f^h,q_f^h\rangle
  +\overline{
   \langle\mathcal S_f^0(\kappa)q_f^h,p_f^h\rangle
  }
 \right).
 \label{eq:parsel-Hermitian-form}
\end{equation*}
Its smallest generalized eigenvalue relative to $\mathcal L_f$ is
\begin{equation*}
 \lambda_f^0(\kappa)
 :=
 \min_{0\neq q_f^h\in P_h^f}
 \frac{s_{f,\mathrm H}^0(q_f^h,q_f^h)}
      {\ell_h(q_f^h,q_f^h)}
 =
 \min_{0\neq q_f^h\in P_h^f}
 \frac{
  \operatorname{Re}
  \langle\mathcal S_f^0(\kappa)q_f^h,q_f^h\rangle
 }{
  \ell_h(q_f^h,q_f^h)
 }.
 \label{eq:parsel-lambda-f-appendix}
\end{equation*}

For every nonzero $q_f^h\in P_h^f$,
\begin{equation*}
 \frac{
  \operatorname{Re}
  \langle\mathcal S_f(\kappa,\delta_2)q_f^h,q_f^h\rangle
 }{
  \ell_h(q_f^h,q_f^h)
 }
 =
 \frac{
  \operatorname{Re}
  \langle\mathcal S_f^0(\kappa)q_f^h,q_f^h\rangle
 }{
  \ell_h(q_f^h,q_f^h)
 }
 +\delta_2.
 \label{eq:parsel-Rayleigh-shift-step}
\end{equation*}
Taking the minimum gives the exact shift identity
\begin{equation*}
 \min_{0\neq q_f^h\in P_h^f}
 \frac{
  \operatorname{Re}
  \langle\mathcal S_f(\kappa,\delta_2)q_f^h,q_f^h\rangle
 }{
  \ell_h(q_f^h,q_f^h)
 }
 =
 \lambda_f^0(\kappa)+\delta_2.
 \label{eq:parsel-generalized-eigenvalue-shift}
\end{equation*}
Therefore,~\eqref{eq:parsel-hermitian-stability} holds if and only if
\begin{equation*}
 \lambda_f^0(\kappa)+\delta_2\geq\beta_*.
 \label{eq:parsel-scalar-H-condition}
\end{equation*}
Imposing $\delta_2\geq0$ yields
\begin{equation*}
 \delta_{2}(\kappa)
 =
 \left[\beta_*-\lambda_f^0(\kappa)\right]_+.
 \label{eq:parsel-delta2-appendix}
\end{equation*}
This value is minimal with respect to the Hermitian stability requirement. If $\lambda_f^0(\kappa)\geq\beta_*$, the unstabilized reduced operator already satisfies~\eqref{eq:parsel-hermitian-stability}, and $\delta_{2}=0$. Otherwise, $\delta_{2}=\beta_*-\lambda_f^0(\kappa)$, and every smaller nonnegative value violates the criterion along a minimizing generalized
eigenfunction.

By \eqref{eq:parsel-pressure-amplification-bound}, $\beta_*^{-1}$ is an upper bound for the amplification measured in the $\mathcal L_f$ norm. Since
\eqref{eq:parsel-hermitian-stability} imposes coercivity on the Hermitian part of a generally nonnormal complex operator, it is a sufficient and potentially conservative stability requirement. It controls the LPS fluctuation component in an energy norm and should not be interpreted as a nodal monotonicity or discrete maximum principle condition.

 
\subsection{Computational realization}
\label{app:parsel-computational-realization}

Algorithm~\ref{alg:parsel-computation} evaluates~\eqref{eq:parsel-delta2-main}
under the invertibility assumption in~\ref{app:parsel-schur-complement}.
Calligraphic symbols denote coefficient matrices:
$\mathcal A_h^0$ is assembled with $\delta_2=0$, and $\mathcal L_h$ acts on
one pressure component. For BP, $\mathcal L_h$ is obtained
from~\eqref{eq:parsel-Lh} by replacing $\mathcal K_h$ with $I$.
\texttt{BlockSplit} applies the same transformation
to trial and test pressure coordinates in both real and imaginary blocks.
The Schur and metric matrices use identical fluctuation coordinates;
ordinary transposes in this real representation correspond to complex
adjoints.

\begingroup
\RestyleAlgo{ruled}
\LinesNumbered
\DontPrintSemicolon
\SetAlgoVlined
\SetAlgoNoEnd
\SetFuncSty{texttt}
\SetCommentSty{textnormal}
\SetAlFnt{\small}
\SetAlCapFnt{\small}
\SetAlCapNameFnt{\small\bfseries}
\SetKwProg{Fn}{Function}{:}{end}
\SetKw{Stop}{stop}
\SetKwFunction{Assemble}{Assemble}
\SetKwFunction{Restrict}{Restrict}
\SetKwFunction{KernelBasis}{KernelBasis}
\SetKwFunction{Complement}{CoordinateComplement}
\SetKwFunction{CheckKernel}{CheckKernel}
\SetKwFunction{BlockSplit}{BlockSplit}
\SetKwFunction{SetUp}{SetUpSolvers}
\SetKwFunction{Solve}{Solve}
\SetKwFunction{ApplySchur}{ApplySchur}
\SetKwFunction{ApplyTranspose}{ApplyTranspose}
\SetKwFunction{Shell}{Shell}
\SetKwFunction{SmallestEigenpair}{SLEPcSmallest}
\SetKwFunction{Save}{Save}
\begin{algorithm}[H]
\caption{Computation of $\delta_2$}
\label{alg:parsel-computation}
\KwIn{Mesh $\mathcal T_h$, physical data, boundary conditions,
$k\geq1$, $\mathrm{method}\in\{\mathrm{BP},\mathrm{LPS}\}$,
$\delta_1$, $\beta_*$, $\mathrm{tol}_{\ker}$, and solver settings.}
\KwOut{$\delta_2$ on successful completion.}
\tcp{Assembly and boundary elimination}
$(\mathcal A_h^0,\mathcal L_h)\gets$
\Assemble{$\mathcal T_h,k,\delta_1,\mathrm{method};\ \delta_2=0$}\;
$(\mathcal A_h^0,\mathcal L_h)\gets$
\Restrict{$\mathcal A_h^0,\mathcal L_h$}\;
\tcp{Pressure kernel and coordinate complement}
$J_r\gets$ \KernelBasis{$P_h,k,\mathrm{method}$}\;
\CheckKernel{$\mathcal L_h,J_r,\mathrm{tol}_{\ker}$}\;
$J_f\gets$ \Complement{$J_r$}\;
\lIf{$J_f$ \textnormal{has no columns}}{\Return{$0$}}
$T_p\gets[J_r,J_f]$\;
$(\mathcal A_c,\mathcal A_{cf},\mathcal A_{fc},\mathcal A_{ff})\gets$
\BlockSplit{$\mathcal A_h^0,T_p$}\;
$\mathcal L_f\gets\operatorname{diag}
\bigl(J_f^\top\mathcal L_hJ_f,\,
      J_f^\top\mathcal L_hJ_f\bigr)$\;
\SetUp{$\mathcal A_c,\mathcal L_f$}\;
\tcp{Schur and transpose actions}
\Fn{\ApplySchur{$\bm x$}}{
  $\bm z\gets$ \Solve{$\mathcal A_c,\mathcal A_{cf}\bm x$}\;
  \Return{$\mathcal A_{ff}\bm x-\mathcal A_{fc}\bm z$}\;
}
\Fn{\ApplyTranspose{$\bm x$}}{
  $\bm z\gets$ \Solve{$\mathcal A_c^\top,\mathcal A_{fc}^\top\bm x$}\;
  \Return{$\mathcal A_{ff}^\top\bm x-\mathcal A_{cf}^\top\bm z$}\;
}
\tcp{Hermitian eigenproblem; $\mathrm{converged}\in\{\mathrm{true},\mathrm{false}\}$}
$\mathcal H_f^0\gets$ \Shell{$\bm x\mapsto
\tfrac12\bigl(\mbox{\normalfont\ttfamily ApplySchur}(\bm x)
+\mbox{\normalfont\ttfamily ApplyTranspose}(\bm x)\bigr)$}\;
$(\widehat\lambda_f^0,\mathrm{info},\mathrm{converged})\gets$
\SmallestEigenpair{$\mathcal H_f^0,\mathcal L_f$}\;
\lIf{$\mathrm{converged}=\mathrm{false}$}{\Stop}
\tcp{Parameter and diagnostic output}
$\delta_2\gets[\beta_*-\widehat\lambda_f^0]_+$\;
\Save{$\delta_2,\widehat\lambda_f^0,\mathrm{info}$}\;
\Return{$\delta_2$}\;
\end{algorithm}
\endgroup

For affine simplicial $P_k$ elements with the projection used here,
\texttt{KernelBasis} constructs $J_r$ by nodal interpolation of the
constrained continuous $P_{k-1}$ space into $P_h$.
The kernel is trivial for the present $P_1$--BP problems. 
\texttt{CoordinateComplement} selects coordinate vectors outside a maximal
set of linearly independent rows of $J_r$; for $P_2$--LPS these are the free edge
coordinates. This avoids an explicit $L^2$ orthogonalization: Schur elimination produces a pencil congruent to that on $P_h^f$.
The implementation used here covers $P_1$--BP and $P_2$--LPS.

PETSc reuses sparse factorizations or preconditioner setups; each
$\mathcal H_f^0$ action uses one solve with $\mathcal A_c$ and one with
$\mathcal A_c^\top$, and the generalized eigensolver also uses
$\mathcal L_f$ solves. SLEPc uses Krylov--Schur by default.
The record $\mathrm{info}$ contains the eigenpair residual, eigenvalue
iteration count, and Schur action counts.
The structural kernel check uses relative threshold
$\max\{100\,\mathrm{tol}_{\ker},10^{-8}\}$; failed checks or inner solves
stop the computation. No dense Schur matrix is formed.
The returned parameter approximates~\eqref{eq:parsel-delta2-main} at the
accuracy of the decomposition and solver tolerances; a small eigenpair
residual alone does not certify a lower spectral bound.

\section*{Data Availability}
The underlying computer codes are available from the corresponding author upon reasonable request.


\bibliographystyle{elsarticle-num}

\end{document}